# Partition of unity systems and B-splines

Daniel J. Greenhoe

*Abstract*: This paper presents the basic principles of partition of unity systems and B-splines. Analysis of these systems are performed using Fourier analysis, multi-resolution analysis, and wavelet analysis.

## Contents

# 1 Background: harmonic analysis

## 1.1 Families of functions

This paper is largely set in the space of *Lebesgue square-integrable functions* $\boldsymbol{L}^2_{\mathbb{R}}$ (Definition 1.2 page 3). The space $\boldsymbol{L}^2_{\mathbb{R}}$ is a subspace of the space $\mathbb{R}^{\mathbb{R}}$, the set of all functions with *domain* $\mathbb{R}$ (the set of real numbers) and *range* $\mathbb{R}$. The space $\mathbb{R}^{\mathbb{R}}$ is a subspace of the space $\mathbb{C}^{\mathbb{C}}$, the set of all functions with *domain* $\mathbb{C}$ (the set of complex numbers) and *range* $\mathbb{C}$. That is, $\boldsymbol{L}^2_{\mathbb{R}} \subseteq \mathbb{R}^{\mathbb{R}} \subseteq \mathbb{C}^{\mathbb{C}}$. In general, the notation $Y^X$ represents the set of all functions with domain $X$ and range $Y$ (Definition 1.1 page 3). Although this notation may seem curious, note that for finite $X$ and finite $Y$, the number of functions (elements) in $Y^X$ is $\left|Y^X\right| = |Y|^{|X|}$.

**Definition 1.1** Let $X$ and $Y$ be sets.
The space $Y^X$ represents the set of all functions with *domain* $X$ and *range* $Y$ such that
$$Y^X \triangleq \{\mathsf{f}(x) \,|\, \mathsf{f}(x) : X \to Y\}$$

**Definition 1.2** Let $\mathbb{R}$ be the set of real numbers, $\mathscr{B}$ the set of *Borel sets* on $\mathbb{R}$, and $\mu$ the standard *Borel measure* on $\mathscr{B}$. Let $\mathbb{R}^{\mathbb{R}}$ be as in Definition 1.1 page 3.

The space of **Lebesgue square-integrable functions** $\boldsymbol{L}^2_{(\mathbb{R},\mathscr{B},\mu)}$ (or $\boldsymbol{L}^2_{\mathbb{R}}$) is defined as
$$\boldsymbol{L}^2_{\mathbb{R}} \triangleq \boldsymbol{L}^2_{(\mathbb{R},\mathscr{B},\mu)} \triangleq \left\{ \mathsf{f} \in \mathbb{R}^{\mathbb{R}} \,\middle|\, \left(\int_{\mathbb{R}} |\mathsf{f}|^2\right)^{\frac{1}{2}} \mathrm{d}\mu < \infty \right\}.$$
The **standard inner product** $\langle \triangle \,|\, \nabla \rangle$ on $\boldsymbol{L}^2_{\mathbb{R}}$ is defined as
$$\langle \mathsf{f}(x) \,|\, \mathsf{g}(x) \rangle \triangleq \int_{\mathbb{R}} \mathsf{f}(x)\mathsf{g}^*(x)\,\mathrm{d}x.$$
The **standard norm** $\|\cdot\|$ on $\boldsymbol{L}^2_{\mathbb{R}}$ is defined as $\|\mathsf{f}(x)\| \triangleq \langle \mathsf{f}(x) \,|\, \mathsf{f}(x) \rangle^{\frac{1}{2}}$.

**Definition 1.3** [1] Let $X$ be a set. The **indicator function** $\mathbb{1} \in \{0,1\}^{2^X}$ is defined as
$$\mathbb{1}_A(x) = \begin{cases} 1 & \text{for } x \in A \qquad \forall x \in X,\ A \in 2^X \\ 0 & \text{for } x \notin A \qquad \forall x \in X,\ A \in 2^X \end{cases}$$
The indicator function $\mathbb{1}$ is also called the **characteristic function**.

## 1.2 Trigonometric functions

### 1.2.1 Definitions

**Lemma 1.4** [2] *Let $\boldsymbol{C}$ be the* SPACE OF ALL CONTINUOUSLY DIFFERENTIABLE REAL FUNCTIONS *and $\frac{\mathsf{d}}{\mathsf{dx}} \in \boldsymbol{C}^{\boldsymbol{C}}$ the differentiation operator.* $\frac{\mathsf{d}}{\mathsf{dx}}^2\mathsf{f} + \mathsf{f} = 0 \iff$
$$\left\{ \begin{array}{rcl} \mathsf{f}(x) &=& \underbrace{[\mathsf{f}](0)\sum_{n=0}^{\infty}(-1)^n\frac{x^{2n}}{(2n)!}}_{\textit{even terms}} + \underbrace{\left[\frac{\mathsf{d}}{\mathsf{dx}}\mathsf{f}\right](0)\sum_{n=0}^{\infty}(-1)^n\frac{x^{2n+1}}{(2n+1)!}}_{\textit{odd terms}} \\ &=& \left(\mathsf{f}(0) + \left[\frac{\mathsf{d}}{\mathsf{dx}}\mathsf{f}\right](0)x\right) - \left(\frac{\mathsf{f}(0)}{2!}x^2 + \frac{\left[\frac{\mathsf{d}}{\mathsf{dx}}\mathsf{f}\right](0)}{3!}x^3\right) + \left(\frac{\mathsf{f}(0)}{4!}x^4 + \frac{\left[\frac{\mathsf{d}}{\mathsf{dx}}\mathsf{f}\right](0)}{5!}x^5\right)\cdots \end{array} \right\}$$

[1] [4], page 126, [57], page 22, [109], page 440
[2] [97], page 156, [84]

**Definition 1.5** [3] Let $\boldsymbol{C}$ be the *space of all continuously differentiable real functions* and $\frac{\mathsf{d}}{\mathsf{dx}} \in \boldsymbol{C}^{\boldsymbol{C}}$ the differentiation operator.

The **cosine** function $\cos(x)$ is the function $\mathsf{f} \in \boldsymbol{C}$ that satisfies the following conditions:

$$\underbrace{\frac{\mathsf{d}}{\mathsf{dx}}^2\mathsf{f} + \mathsf{f} = 0}_{\text{2nd order homogeneous differential equation}} \qquad \underbrace{\mathsf{f}(0) = 1}_{\text{1st initial condition}} \qquad \underbrace{\left[\frac{\mathsf{d}}{\mathsf{dx}}\mathsf{f}\right](0) = 0}_{\text{2nd initial condition}}$$

The **sine** function $\sin(x)$ is the function $\mathsf{g} \in \boldsymbol{C}$ that satisfies the following conditions:

$$\underbrace{\frac{\mathsf{d}}{\mathsf{dx}}^2\mathsf{g} + \mathsf{g} = 0}_{\text{2nd order homogeneous differential equation}} \qquad \underbrace{\mathsf{g}(0) = 0}_{\text{1st initial condition}} \qquad \underbrace{\left[\frac{\mathsf{d}}{\mathsf{dx}}\mathsf{g}\right](0) = 1}_{\text{2nd initial condition}}$$

**Theorem 1.6** [4]

$$\begin{aligned}
\cos(x) &= \sum_{n=0}^{\infty}(-1)^n\frac{x^{2n}}{(2n)!} &&= 1 - \frac{x^2}{2} + \frac{x^4}{4!} - \frac{x^6}{6!} + \cdots && \forall x\in\mathbb{R}\\
\sin(x) &= \sum_{n=0}^{\infty}(-1)^n\frac{x^{2n+1}}{(2n+1)!} &&= x - \frac{x^3}{3!} + \frac{x^5}{5!} - \frac{x^7}{7!} + \cdots && \forall x\in\mathbb{R}
\end{aligned}$$

**Proposition 1.7** [5] *Let $\boldsymbol{C}$ be the space of all continuously differentiable real functions and $\frac{\mathsf{d}}{\mathsf{dx}} \in \boldsymbol{C}^{\boldsymbol{C}}$ the differentiation operator. Let $\mathsf{f}'(0) \triangleq \left[\frac{\mathsf{d}}{\mathsf{dx}}\mathsf{f}\right](0)$.*

$$\underbrace{\frac{\mathsf{d}}{\mathsf{dx}}^2\mathsf{f} + \mathsf{f} = 0}_{\text{2ND ORDER HOMOGENEOUS DIFFERENTIAL EQUATION}} \iff \mathsf{f}(x) = \mathsf{f}(0)\cos(x) + \mathsf{f}'(0)\sin(x) \qquad \forall \mathsf{f}\in\boldsymbol{C}, \forall x\in\mathbb{R}$$

**Theorem 1.8** [6] *Let $\frac{\mathsf{d}}{\mathsf{dx}} \in \boldsymbol{C}^{\boldsymbol{C}}$ be the differentiation operator.*

$$\begin{aligned}
\frac{\mathsf{d}}{\mathsf{dx}}\cos(x) &= -\sin(x) && \forall x\in\mathbb{R}\\
\frac{\mathsf{d}}{\mathsf{dx}}\sin(x) &= \cos(x) && \forall x\in\mathbb{R}
\end{aligned}$$

---

[3] 📖 [97], page 157, 📖 [38], pages 228–229

[4] 📖 [97], page 157

[5] 📖 [97], page 157. The general solution for the *non-homogeneous* equation $\frac{\mathsf{d}}{\mathsf{dx}}^2\mathsf{f}(x) + \mathsf{f}(x) = \mathsf{g}(x)$ with initial conditions $\mathsf{f}(a) = 1$ and $\mathsf{f}'(a) = \rho$ is

$$\mathsf{f}(x) = \cos(x) + \rho\sin(x) + \int_a^x \mathsf{g}(y)\sin(x-y)\,\mathsf{d}y.$$

This type of equation is called a *Volterra integral equation of the second type*. References: 📖 [39], page 371, 📖 [84]. Volterra equation references: 📖 [92], page 99, 📖 [78], 📖 [79].

[6] 📖 [97], page 157

### 1.2.2 The complex exponential

**Definition 1.9** The function $\mathsf{f} \in \mathbb{C}^{\mathbb{C}}$ is the **exponential function** $\exp(ix) \triangleq \mathsf{f}(x)$ if

1. $\frac{\mathsf{d}}{\mathsf{dx}}^2 \mathsf{f} + \mathsf{f} = 0$ (second order homogeneous differential equation) and
2. $\mathsf{f}(0) = 1$ (first initial condition) and
3. $\left[\frac{\mathsf{d}}{\mathsf{dx}} \mathsf{f}\right](0) = i$ (second initial condition).

**Theorem 1.10** (Euler's identity) [7]

$$e^{ix} = \cos(x) + i\sin(x) \qquad \forall x \in \mathbb{R}$$

**Corollary 1.11**

$$e^{ix} = \sum_{n\in\mathbb{W}} \frac{(ix)^n}{n!} \qquad \forall x \in \mathbb{R}$$

**Corollary 1.12** [8]

$$e^{i\pi} + 1 = 0$$

The exponential has two properties that makes it extremely special:

- The exponential is an eigenvalue of any LTI operator (Theorem 1.13 page 5).
- The exponential generates a continuous point spectrum for the differential operator.

**Theorem 1.13** [9] *Let* $\mathbf{L}$ *be an operator with kernel* $\mathrm{h}(t, \omega)$ *and*

$$\hat{\mathrm{h}}(s) \triangleq \left\langle \mathrm{h}(t,\omega) \mid e^{st} \right\rangle \qquad \text{(Laplace transform)}.$$

$$\left.\begin{array}{ll} 1. & \mathbf{L} \textit{ is linear and} \\ 2. & \mathbf{L} \textit{ is time-invariant} \end{array}\right\} \implies \mathbf{L}e^{st} = \underbrace{\hat{\mathrm{h}}^*(-s)}_{\textit{eigenvalue}} \; \underbrace{e^{st}}_{\textit{eigenvector}}$$

### 1.2.3 Trigonometric Identities

**Corollary 1.14** (Euler formulas) [10]

$$\begin{aligned} \cos(x) &= \Re\left(e^{ix}\right) = \frac{e^{ix} + e^{-ix}}{2} && \forall x \in \mathbb{R} \\ \sin(x) &= \Im\left(e^{ix}\right) = \frac{e^{ix} - e^{-ix}}{2i} && \forall x \in \mathbb{R} \end{aligned}$$

[7] [34], [16], page 12
[8] [34], [35], http://www.daviddarling.info/encyclopedia/E/Eulers_formula.html
[9] [87], page 2, …page 2 online: http://www.cmap.polytechnique.fr/~mallat/WTintro.pdf
[10] [34], [16], page 12

**Theorem 1.15** [11]

$$e^{(\alpha+\beta)} = e^{\alpha}\, e^{\beta} \qquad \forall \alpha,\beta\in\mathbb{C}$$

**Theorem 1.16** (shift identities)

$$\begin{array}{rclr|rclr}
\cos\left(x+\frac{\pi}{2}\right) &=& -\sin x & \forall x\in\mathbb{R} & \sin\left(x+\frac{\pi}{2}\right) &=& +\cos x & \forall x\in\mathbb{R}\\
\cos\left(x-\frac{\pi}{2}\right) &=& +\sin x & \forall x\in\mathbb{R} & \sin\left(x-\frac{\pi}{2}\right) &=& -\cos x & \forall x\in\mathbb{R}
\end{array}$$

**Theorem 1.17** (product identities)

$$\begin{array}{rcrlclr}
\cos x\cos y &=& & \frac{1}{2}\cos(x-y) &+& \frac{1}{2}\cos(x+y) & \forall x,y\in\mathbb{R}\\
\cos x\sin y &=& - & \frac{1}{2}\sin(x-y) &+& \frac{1}{2}\sin(x+y) & \forall x,y\in\mathbb{R}\\
\sin x\cos y &=& & \frac{1}{2}\sin(x-y) &+& \frac{1}{2}\sin(x+y) & \forall x,y\in\mathbb{R}\\
\sin x\sin y &=& & \frac{1}{2}\cos(x-y) &-& \frac{1}{2}\cos(x+y) & \forall x,y\in\mathbb{R}
\end{array}$$

**Theorem 1.18** (double angle formulas) [12]

$$\begin{array}{rclr}
\cos(x+y) &=& \cos x\cos y - \sin x\sin y & \forall x,y\in\mathbb{R}\\
\sin(x+y) &=& \sin x\cos y + \cos x\sin y & \forall x,y\in\mathbb{R}\\
\tan(x+y) &=& \dfrac{\tan x+\tan y}{1-\tan x\tan y} & \forall x,y\in\mathbb{R}
\end{array}$$

**Theorem 1.19** (squared identities)

$$\begin{array}{rclr}
\cos^2 x &=& \frac{1}{2}(1+\cos 2x) & \forall x\in\mathbb{R}\\
\sin^2 x &=& \frac{1}{2}(1-\cos 2x) & \forall x\in\mathbb{R}\\
\cos^2 x+\sin^2 x &=& 1 & \forall x\in\mathbb{R}
\end{array}$$

## 1.3 Fourier Series

The *Fourier Series* expansion of a periodic function is simply a complex trigonometric polynomial. In the special case that the periodic function is even, then the Fourier Series expansion is a cosine polynomial.

**Definition 1.20** [13] The **Fourier Series operator** $\hat{\mathbf{F}} : \boldsymbol{L}^2_{\mathbb{R}} \to \boldsymbol{\ell}^2_{\mathbb{R}}$ is defined as

$$\left[\hat{\mathbf{F}}\mathsf{f}\right](n) \triangleq \frac{1}{\sqrt{\tau}}\int_0^{\tau} \mathsf{f}(x)e^{-i\frac{2\pi}{\tau}nx}\,\mathrm{d}x \qquad \forall \mathsf{f}\in\left\{\mathsf{f}\in\boldsymbol{L}^2_{\mathbb{R}} \,\middle|\, \mathsf{f} \text{ is periodic with period } \tau\right\}$$

**Theorem 1.21** *Let $\hat{\mathbf{F}}$ be the Fourier Series operator.*
*The **inverse Fourier Series** operator $\hat{\mathbf{F}}^{-1}$ is given by*

$$\left[\hat{\mathbf{F}}^{-1}\left(\!\left(\tilde{x}_n\right)\!\right)_{n\in\mathbb{Z}}\right](x) \triangleq \frac{1}{\sqrt{\tau}}\sum_{n\in\mathbb{Z}}\tilde{x}_n e^{i\frac{2\pi}{\tau}nx} \qquad \forall (\!(\tilde{x}_n)\!)\in\boldsymbol{\ell}^2_{\mathbb{R}}$$

---

[11] [98], page 1

[12] Expressions for $\cos(\alpha+\beta)$, $\sin(\alpha+\beta)$, and $\sin^2 x$ appear in works as early as [95]. Reference: http://en.wikipedia.org/wiki/History_of_trigonometric_functions

[13] [70], page 3

**Theorem 1.22** *The **Fourier Series adjoint** operator* $\hat{\mathbf{F}}^*$ *is given by*

$$\hat{\mathbf{F}}^* = \hat{\mathbf{F}}^{-1}$$

✎PROOF:

$$\begin{aligned}
\left\langle \hat{\mathbf{F}}\mathsf{x}(x) \,\middle|\, \tilde{\mathsf{y}}(n) \right\rangle_{\mathbb{Z}} &= \left\langle \frac{1}{\sqrt{\tau}} \int_0^\tau \mathsf{x}(x) e^{-i\frac{2\pi}{\tau}nx} \, \mathrm{d}x \,\middle|\, \tilde{\mathsf{y}}(n) \right\rangle_{\mathbb{Z}} && \text{by definition of } \hat{\mathbf{F}} \text{ Definition 1.20 page 6} \\
&= \frac{1}{\sqrt{\tau}} \int_0^\tau \mathsf{x}(x) \left\langle e^{-i\frac{2\pi}{\tau}nx} \,\middle|\, \tilde{\mathsf{y}}(n) \right\rangle_{\mathbb{Z}} \mathrm{d}x && \text{by additivity property of } \langle \triangle \mid \nabla \rangle \\
&= \int_0^\tau \mathsf{x}(x) \frac{1}{\sqrt{\tau}} \left\langle \tilde{\mathsf{y}}(n) \mid e^{-i\frac{2\pi}{\tau}nx} \right\rangle_{\mathbb{Z}}^* \mathrm{d}x && \text{by property of } \langle \triangle \mid \nabla \rangle \\
&= \int_0^\tau \mathsf{x}(x) \left[ \hat{\mathbf{F}}^{-1} \tilde{\mathsf{y}}(n) \right]^* \mathrm{d}x && \text{by definition of } \hat{\mathbf{F}}^{-1} \text{ page 6} \\
&= \left\langle \mathsf{x}(x) \mid \underbrace{\hat{\mathbf{F}}^{-1}}_{\hat{\mathbf{F}}^*} \tilde{\mathsf{y}}(n) \right\rangle_{\mathbb{R}}
\end{aligned}$$

✎

The Fourier Series operator has several nice properties:

- $\hat{\mathbf{F}}$ is *unitary* (Corollary 1.23 page 7).
- Because $\hat{\mathbf{F}}$ is unitary, it automatically has several other nice properties such as being *isometric*, and satisfying *Parseval's equation*, satisfying *Plancheral's formula*, and more (Corollary 1.24 page 7).

**Corollary 1.23** *Let* $\mathbf{I}$ *be the identity operator and let* $\hat{\mathbf{F}}$ *be the Fourier Series operator with adjoint* $\hat{\mathbf{F}}^*$.

$$\hat{\mathbf{F}}\hat{\mathbf{F}}^* = \hat{\mathbf{F}}^*\hat{\mathbf{F}} = \mathbf{I} \qquad \text{(}\hat{\mathbf{F}} \textit{ is unitary...and thus also normal and isometric)}$$

✎PROOF: This follows directly from the fact that $\hat{\mathbf{F}}^* = \hat{\mathbf{F}}^{-1}$ (Theorem 1.22 (page 7)). ✎

**Corollary 1.24** *Let* $\hat{\mathbf{F}}$ *be the Fourier series operator,* $\hat{\mathbf{F}}^*$ *be its adjoint, and* $\hat{\mathbf{F}}^{-1}$ *be its inverse.*

$$\begin{array}{rclcll}
\mathcal{R}(\hat{\mathbf{F}}) &=& \mathcal{R}(\hat{\mathbf{F}}^{-1}) &=& \boldsymbol{L}^2_{\mathbb{R}} & \\
\left\|\left|\hat{\mathbf{F}}\right|\right\| &=& \left\|\left|\hat{\mathbf{F}}^{-1}\right|\right\| &=& 1 & \text{(UNITARY)} \\
\left\langle \hat{\mathbf{F}}\boldsymbol{x} \,\middle|\, \hat{\mathbf{F}}\boldsymbol{y} \right\rangle &=& \left\langle \hat{\mathbf{F}}^{-1}\boldsymbol{x} \,\middle|\, \hat{\mathbf{F}}^{-1}\boldsymbol{y} \right\rangle &=& \langle \boldsymbol{x} \mid \boldsymbol{y} \rangle & \text{(PARSEVAL'S EQUATION)} \\
\left\|\hat{\mathbf{F}}\boldsymbol{x}\right\| &=& \left\|\hat{\mathbf{F}}^{-1}\boldsymbol{x}\right\| &=& \|\boldsymbol{x}\| & \text{(PLANCHEREL'S FORMULA)} \\
\left\|\hat{\mathbf{F}}\boldsymbol{x} - \hat{\mathbf{F}}\boldsymbol{y}\right\| &=& \left\|\hat{\mathbf{F}}^{-1}\boldsymbol{x} - \hat{\mathbf{F}}^{-1}\boldsymbol{y}\right\| &=& \|\boldsymbol{x} - \boldsymbol{y}\| & \text{(ISOMETRIC)}
\end{array}$$

✎PROOF: These results follow directly from the fact that $\hat{\mathbf{F}}$ is unitary (Corollary 1.23 page 7) and from the properties of unitary operators. ✎

**Theorem 1.25** *The set*

$$\left\{ \frac{1}{\sqrt{\tau}} e^{i\frac{2\pi}{\tau}nx} \,\middle|\, n \in \mathbb{Z} \right\}$$

*is an* ORTHONORMAL BASIS *for all functions* $\mathsf{f}(x)$ *with support in* $[0, \tau]$.

## 1.4 Fourier Transform

**Definition 1.26** [14] The **Fourier Transform** operator $\tilde{\mathbf{F}}$ is defined as[15]

$$\left[\tilde{\mathbf{F}}\mathsf{f}\right](\omega) \triangleq \frac{1}{\sqrt{2\pi}}\int_{\mathbb{R}} \mathsf{f}(x)\,e^{-i\omega x}\,\mathrm{d}x \qquad \forall \mathsf{f}\in \boldsymbol{L}^2_{(\mathbb{R},\mathscr{B},\mu)}$$

This definition of the Fourier Transform is also called the **unitary Fourier Transform**.

**Remark 1.27** (**Fourier transform scaling factor**) [16] If the Fourier transform operator $\tilde{\mathbf{F}}$ and inverse Fourier transform operator $\tilde{\mathbf{F}}^{-1}$ are defined as

$$\tilde{\mathbf{F}}\mathsf{f}(x) \triangleq A\int_{\mathbb{R}} \mathsf{f}(x)e^{-i\omega x}\,\mathrm{d}x \quad\text{and}\quad \tilde{\mathbf{F}}^{-1}\tilde{\mathsf{f}}(\omega) \triangleq B\int_{\mathbb{R}} \mathsf{F}(\omega)e^{i\omega x}\,\mathrm{d}\omega \quad,$$

then $A$ and $B$ can be any constants as long as $AB = \frac{1}{2\pi}$. The Fourier transform is often defined with the scaling factor $A$ set equal to 1 such that $\left[\tilde{\mathbf{F}}\mathsf{f}(x)\right](\omega) \triangleq \int_{\mathbb{R}} \mathsf{f}(x)\,e^{-i\omega x}\,\mathrm{d}x$. In this case, the inverse Fourier transform operator $\tilde{\mathbf{F}}^{-1}$ is either defined as

- $\left[\tilde{\mathbf{F}}^{-1}\mathsf{f}(x)\right](f) \triangleq \int_{\mathbb{R}} \mathsf{f}(x)\,e^{i2\pi f x}\,\mathrm{d}x$ (using oscillatory frequency free variable $f$) or
- $\left[\tilde{\mathbf{F}}^{-1}\mathsf{f}(x)\right](\omega) \triangleq \frac{1}{2\pi}\int_{\mathbb{R}} \mathsf{f}(x)\,e^{i\omega x}\,\mathrm{d}x$ (using angular frequency free variable $\omega$).

In short, the $2\pi$ has to show up somewhere, either in the argument of the exponential ($e^{-i2\pi f t}$) or in front of the integral ($\frac{1}{2\pi}\int\cdots$). One could argue that it is unnecessary to burden the exponential argument with the $2\pi$ factor ($e^{-i2\pi f t}$), and thus could further argue in favor of using the angular frequency variable $\omega$ thus giving the inverse operator definition $\left[\tilde{\mathbf{F}}^{-1}\mathsf{f}(x)\right](\omega) \triangleq \frac{1}{2\pi}\int_{\mathbb{R}} \mathsf{f}(x)\,e^{-i\omega x}\,\mathrm{d}x$. But this causes a new problem. In this case, the Fourier operator $\tilde{\mathbf{F}}$ is not *unitary* (see Theorem 1.29 page 8)—in particular, $\tilde{\mathbf{F}}\tilde{\mathbf{F}}^* \neq \mathbf{I}$, where $\tilde{\mathbf{F}}^*$ is the *adjoint* of $\tilde{\mathbf{F}}$; but rather, $\tilde{\mathbf{F}}\left(\frac{1}{2\pi}\tilde{\mathbf{F}}^*\right) = \left(\frac{1}{2\pi}\tilde{\mathbf{F}}^*\right)\tilde{\mathbf{F}} = \mathbf{I}$. But if we define the operators $\tilde{\mathbf{F}}$ and $\tilde{\mathbf{F}}^{-1}$ to both have the scaling factor $\frac{1}{\sqrt{2\pi}}$, then $\tilde{\mathbf{F}}$ and $\tilde{\mathbf{F}}^{-1}$ are inverses *and* $\tilde{\mathbf{F}}$ is *unitary*—that is, $\tilde{\mathbf{F}}\tilde{\mathbf{F}}^* = \tilde{\mathbf{F}}^*\tilde{\mathbf{F}} = \mathbf{I}$.

**Theorem 1.28** (Inverse Fourier transform) [17] *Let $\tilde{\mathbf{F}}$ be the Fourier Transform operator (Definition 1.26 page 8). The inverse $\tilde{\mathbf{F}}^{-1}$ of $\tilde{\mathbf{F}}$ is*

$$\left[\tilde{\mathbf{F}}^{-1}\tilde{\mathsf{f}}\right](x) \triangleq \frac{1}{\sqrt{2\pi}}\int_{\mathbb{R}} \tilde{\mathsf{f}}(\omega)e^{i\omega x}\,\mathrm{d}\omega \qquad \forall \tilde{\mathsf{f}}\in \boldsymbol{L}^2_{(\mathbb{R},\mathscr{B},\mu)}$$

**Theorem 1.29** *Let $\tilde{\mathbf{F}}$ be the Fourier Transform operator with inverse $\tilde{\mathbf{F}}^{-1}$ and adjoint $\tilde{\mathbf{F}}^*$.*

$$\tilde{\mathbf{F}}^* = \tilde{\mathbf{F}}^{-1}$$

---

[14] [11], page 363, [20], page 13, [85], page 144, [73], pages 374–375, [41], [42], page 336?

[16] [50], page 274, ⟨Remark F.1⟩, [20], page 13, [66], pages xxxi–xxxii, [73], pages 374–375

[17] [20], page 13

✎PROOF:

$$\begin{aligned}
\left\langle \tilde{\mathbf{F}}\mathsf{f} \,\middle|\, \mathsf{g} \right\rangle &= \left\langle \frac{1}{\sqrt{2\pi}} \int_{\mathbb{R}} \mathsf{f}(x) e^{-i\omega x}\, \mathrm{d}x \,\middle|\, \mathsf{g}(\omega) \right\rangle && \text{by definition of } \tilde{\mathbf{F}} \text{ page 8} \\
&= \frac{1}{\sqrt{2\pi}} \int_{\mathbb{R}} \mathsf{f}(x) \left\langle e^{-i\omega x} \,\middle|\, \mathsf{g}(\omega) \right\rangle \mathrm{d}x && \text{by } \textit{additive property} \text{ of } \langle \triangle \mid \nabla \rangle \\
&= \int_{\mathbb{R}} \mathsf{f}(x) \frac{1}{\sqrt{2\pi}} \left\langle \mathsf{g}(\omega) \mid e^{-i\omega x} \right\rangle^* \mathrm{d}x && \text{by } \textit{conjugate symmetric property} \text{ of } \langle \triangle \mid \nabla \rangle \\
&= \left\langle \mathsf{f}(x) \,\middle|\, \frac{1}{\sqrt{2\pi}} \left\langle \mathsf{g}(\omega) \mid e^{-i\omega x} \right\rangle \right\rangle && \text{by definition of } \langle \triangle \mid \nabla \rangle \\
&= \left\langle \mathsf{f} \,\middle|\, \underbrace{\tilde{\mathbf{F}}^{-1}}_{\tilde{\mathbf{F}}^*} \mathsf{g} \right\rangle && \text{by Theorem 1.28 page 8}
\end{aligned}$$

✎

The Fourier Transform operator has several nice properties:

- $\tilde{\mathbf{F}}$ is *unitary* (Corollary 1.30—next corollary).
- Because $\tilde{\mathbf{F}}$ is unitary, it automatically has several other nice properties (Theorem 1.31 page 9).

**Corollary 1.30** *Let* $\mathbf{I}$ *be the identity operator and let* $\tilde{\mathbf{F}}$ *be the Fourier Transform operator with adjoint* $\tilde{\mathbf{F}}^*$ *and inverse* $\tilde{\mathbf{F}}^{-1}$*.*

$$\underbrace{\tilde{\mathbf{F}}\tilde{\mathbf{F}}^* = \tilde{\mathbf{F}}^*\tilde{\mathbf{F}} = \mathbf{I}}_{\tilde{\mathbf{F}}^* = \tilde{\mathbf{F}}^{-1}} \qquad (\tilde{\mathbf{F}} \textit{ is unitary})$$

✎PROOF: This follows directly from the fact that $\tilde{\mathbf{F}}^* = \tilde{\mathbf{F}}^{-1}$ (Theorem 1.29 page 8). ✎

**Theorem 1.31** *Let* $\tilde{\mathbf{F}}$ *be the Fourier transform operator with adjoint* $\tilde{\mathbf{F}}^*$ *and inverse* $\tilde{\mathbf{F}}$*. Let* $|||\cdot|||$ *be the operator norm with respect to the vector norm* $\|\cdot\|$ *with respect to the Hilbert space* $\left(\mathbb{C}^{\mathbb{R}}, \langle \triangle \mid \nabla \rangle\right)$*. Let* $\mathcal{R}(\mathbf{A})$ *be the range of an operator* $\mathbf{A}$*.*

$$\begin{array}{rclcll}
\mathcal{R}(\mathbf{F}\tau) &=& \mathcal{R}(\tilde{\mathbf{F}}^{-1}) &=& \boldsymbol{L}^2_{\mathbb{R}} & \\
|||\tilde{\mathbf{F}}||| &=& |||\tilde{\mathbf{F}}^{-1}||| &=& 1 & (\text{UNITARY}) \\
\left\langle \tilde{\mathbf{F}}\mathsf{f} \,\middle|\, \tilde{\mathbf{F}}\mathsf{g} \right\rangle &=& \left\langle \tilde{\mathbf{F}}^{-1}\mathsf{f} \,\middle|\, \tilde{\mathbf{F}}^{-1}\mathsf{g} \right\rangle &=& \langle \mathsf{f} \mid \mathsf{g} \rangle & (\text{PARSEVAL'S EQUATION}) \\
\left\|\tilde{\mathbf{F}}\mathsf{f}\right\| &=& \left\|\tilde{\mathbf{F}}^{-1}\mathsf{f}\right\| &=& \|\mathsf{f}\| & (\text{PLANCHEREL'S FORMULA}) \\
\left\|\tilde{\mathbf{F}}\mathsf{f} - \tilde{\mathbf{F}}\mathsf{g}\right\| &=& \left\|\tilde{\mathbf{F}}^{-1}\mathsf{f} - \tilde{\mathbf{F}}^{-1}\mathsf{g}\right\| &=& \|\mathsf{f} - \mathsf{g}\| & (\text{ISOMETRIC})
\end{array}$$

✎PROOF: These results follow directly from the fact that $\tilde{\mathbf{F}}$ is unitary (Corollary 1.30 page 9) and from the properties of unitary operators. ✎

**Theorem 1.32** (Shift relations) [18] *Let $\tilde{\mathbf{F}}$ be the Fourier transform operator.*

$$\begin{aligned}\tilde{\mathbf{F}}[\mathsf{f}(x-u)](\omega) &= e^{-i\omega u}\left[\tilde{\mathbf{F}}\mathsf{f}(x)\right](\omega)\\ \left[\tilde{\mathbf{F}}\left(e^{ivx}\mathsf{g}(x)\right)\right](\omega) &= \left[\tilde{\mathbf{F}}\mathsf{g}(x)\right](\omega-v)\end{aligned}$$

**Theorem 1.33** (Complex conjugate) [19] *Let $\tilde{\mathbf{F}}$ be the Fourier Transform operator and $*$ represent the complex conjugate operation on the set of complex numbers.*

$$\tilde{\mathbf{F}}\mathsf{f}^*(-x) = \left[\tilde{\mathbf{F}}\mathsf{f}(x)\right]^* \qquad \forall \mathsf{f}\in \boldsymbol{L}^2_{(\mathbb{R},\mathscr{B},\mu)}$$

**Definition 1.34** [20] The **convolution operation** is defined as

$$\left[\mathsf{f}\star\mathsf{g}\right](x) \triangleq \mathsf{f}(x)\star\mathsf{g}(x) \triangleq \int_{u\in\mathbb{R}} \mathsf{f}(u)\mathsf{g}(x-u)\,\mathrm{d}u \qquad \forall \mathsf{f},\mathsf{g}\in \boldsymbol{L}^2_{(\mathbb{R},\mathscr{B},\mu)}$$

Theorem 1.35 (next) demonstrates that multiplication in the "time domain" is equivalent to convolution in the "frequency domain" and vice-versa.

**Theorem 1.35** (convolution theorem) [21] *Let $\tilde{\mathbf{F}}$ be the Fourier Transform operator and $\star$ the convolution operator.*

$$\begin{aligned}\underbrace{\tilde{\mathbf{F}}\left[\mathsf{f}(x)\star\mathsf{g}(x)\right](\omega)}_{\text{convolution in "time domain"}} &= \underbrace{\sqrt{2\pi}\left[\tilde{\mathbf{F}}\mathsf{f}\right](\omega)\left[\tilde{\mathbf{F}}\mathsf{g}\right](\omega)}_{\text{multiplication in "frequency domain"}} && \forall \mathsf{f},\mathsf{g}\in \boldsymbol{L}^2_{(\mathbb{R},\mathscr{B},\mu)}\\ \underbrace{\tilde{\mathbf{F}}\left[\mathsf{f}(x)\mathsf{g}(x)\right](\omega)}_{\text{multiplication in "time domain"}} &= \underbrace{\frac{1}{\sqrt{2\pi}}\left[\tilde{\mathbf{F}}\mathsf{f}\right](\omega)\star\left[\tilde{\mathbf{F}}\mathsf{g}\right](\omega)}_{\text{convolution in "frequency domain"}} && \forall \mathsf{f},\mathsf{g}\in \boldsymbol{L}^2_{(\mathbb{R},\mathscr{B},\mu)}.\end{aligned}$$

✎PROOF:

$$\begin{aligned}\tilde{\mathbf{F}}\left[\mathsf{f}(x)\star\mathsf{g}(x)\right](\omega) &= \tilde{\mathbf{F}}\left[\int_{u\in\mathbb{R}}\mathsf{f}(u)\mathsf{g}(x-u)\,\mathrm{d}u\right](\omega) && \text{by def. of } \star \text{ (Definition 1.34 page 10)}\\ &= \int_{u\in\mathbb{R}}\mathsf{f}(u)\left[\tilde{\mathbf{F}}\mathsf{g}(x-u)\right](\omega)\,\mathrm{d}u\\ &= \int_{u\in\mathbb{R}}\mathsf{f}(u)e^{-i\omega u}\left[\tilde{\mathbf{F}}\mathsf{g}(x)\right](\omega)\,\mathrm{d}u && \text{by Theorem 1.32 page 10}\\ &= \sqrt{2\pi}\underbrace{\left(\frac{1}{\sqrt{2\pi}}\int_{u\in\mathbb{R}}\mathsf{f}(u)e^{-i\omega u}\,\mathrm{d}u\right)}_{\left[\tilde{\mathbf{F}}\mathsf{f}\right](\omega)}\left[\tilde{\mathbf{F}}\mathsf{g}\right](\omega)\end{aligned}$$

---

[18] 📚 [50], page 276, ⟨Theorem F.4⟩
[19] 📚 [50], page 276, ⟨Theorem F.5⟩
[20] 📚 [9], page 6
[21] 📚 [50], pages 277–278, ⟨Theorem F.6⟩, 📄 [51], ⟨Theorem 2.31⟩

$$
\begin{aligned}
&= \sqrt{2\pi}\left[\tilde{\mathbf{F}}\mathsf{f}\right](\omega)\ \left[\tilde{\mathbf{F}}\mathsf{g}\right](\omega) && \text{by definition of } \tilde{\mathbf{F}} \text{ (Definition 1.26 page 8)}\\
\tilde{\mathbf{F}}[\mathsf{f}(x)\mathsf{g}(x)](\omega) &= \tilde{\mathbf{F}}\left[\left(\tilde{\mathbf{F}}^{-1}\tilde{\mathbf{F}}\mathsf{f}(x)\right)\ \mathsf{g}(x)\right](\omega) && \text{by definition of operator inverse}\\
&= \tilde{\mathbf{F}}\left[\left(\frac{1}{\sqrt{2\pi}}\int_{v\in\mathbb{R}}\left[\tilde{\mathbf{F}}\mathsf{f}(x)\right](v)e^{ivx}\,\mathrm{d}v\right)\ \mathsf{g}(x)\right](\omega) && \text{by Theorem 1.28 page 8}\\
&= \frac{1}{\sqrt{2\pi}}\int_{v\in\mathbb{R}}\left[\tilde{\mathbf{F}}\mathsf{f}(x)\right](v)\left[\tilde{\mathbf{F}}\left(e^{ivx}\ \mathsf{g}(x)\right)\right](\omega,v)\,\mathrm{d}v\\
&= \frac{1}{\sqrt{2\pi}}\int_{v\in\mathbb{R}}\left[\tilde{\mathbf{F}}\mathsf{f}(x)\right](v)\left[\tilde{\mathbf{F}}\mathsf{g}(x)\right](\omega-v)\,\mathrm{d}v && \text{by Theorem 1.32 page 10}\\
&= \frac{1}{\sqrt{2\pi}}\left[\tilde{\mathbf{F}}\mathsf{f}\right](\omega)\star\left[\tilde{\mathbf{F}}\mathsf{g}\right](\omega) && \text{by def. of } \star \text{ (Definition 1.34 page 10)}
\end{aligned}
$$

✏

## 1.5 Z-transform

**Definition 1.36** [22] Let $X^Y$ be the set of all functions from a set $Y$ to a set $X$. Let $\mathbb{Z}$ be the set of integers. A function $\mathsf{f}$ in $X^Y$ is a **sequence** over $X$ if $Y = \mathbb{Z}$.
A sequence may be denoted in the form $(\!(x_n)\!)_{n\in\mathbb{Z}}$ or simply as $(\!(x_n)\!)$.

**Definition 1.37** [23] Let $(\mathbb{F}, +, \cdot)$ be a field. The **space of all absolutely square summable sequences** $\boldsymbol{\ell}^2_{\mathbb{F}}$ over $\mathbb{F}$ is defined as

$$\boldsymbol{\ell}^2_{\mathbb{F}} \triangleq \left\{ (\!(x_n)\!)_{n\in\mathbb{Z}} \left| \sum_{n\in\mathbb{Z}} |x_n|^2 < \infty \right. \right\}$$

The space $\boldsymbol{\ell}^2_{\mathbb{R}}$ is an example of a *separable Hilbert space*. In fact, $\boldsymbol{\ell}^2_{\mathbb{R}}$ is the *only* separable Hilbert space in the sense that all separable Hilbert spaces are isomorphically equivalent. For example, $\boldsymbol{\ell}^2_{\mathbb{R}}$ is isomorphic to $\boldsymbol{L}^2_{\mathbb{R}}$, the *space of all absolutely square Lebesgue integrable functions*.

**Definition 1.38** The **convolution** operation $\ast$ is defined as

$$(\!(x_n)\!) \ast (\!(y_n)\!) \triangleq \left(\!\!\left( \sum_{m\in\mathbb{Z}} x_m y_{n-m} \right)\!\!\right)_{n\in\mathbb{Z}} \qquad \forall (\!(x_n)\!)_{n\in\mathbb{Z}},(\!(y_n)\!)_{n\in\mathbb{Z}}\in\boldsymbol{\ell}^2_{\mathbb{R}}$$

**Definition 1.39** [24]
The **z-transform Z** of $(\!(x_n)\!)_{n\in\mathbb{Z}}$ is defined as

$$\left[\mathbf{Z}\,(\!(x_n)\!)\right](z) \triangleq \underbrace{\sum_{n\in\mathbb{Z}} x_n z^{-n}}_{\text{Laurent series}} \qquad \forall (\!(x_n)\!)\in\boldsymbol{\ell}^2_{\mathbb{R}}$$

[22] [17], page 1, [108], page 23, ⟨Definition 2.1⟩, [67], page 31
[23] [77], page 347, ⟨Example 5.K⟩
[24] *Laurent series*: [1], page 49

**Proposition 1.40** [25] *Let $\star$ be the* CONVOLUTION OPERATOR *(Definition 1.38 page 11).*

$$\left(\!\left(x_n\right)\!\right) \star \left(\!\left(y_n\right)\!\right) = \left(\!\left(y_n\right)\!\right) \star \left(\!\left(x_n\right)\!\right) \qquad \forall \left(\!\left(x_n\right)\!\right)_{n\in\mathbb{Z}}, \left(\!\left(y_n\right)\!\right)_{n\in\mathbb{Z}} \in \ell^2_{\mathbb{R}} \qquad (\star \text{ is COMMUTATIVE})$$

**Theorem 1.41** [26] *Let $\star$ be the convolution operator (Definition 1.38 page 11).*

$$\underbrace{\mathbf{Z}\left(\left(\!\left(x_n\right)\!\right) \star \left(\!\left(y_n\right)\!\right)\right)}_{\text{sequence convolution}} = \underbrace{\left(\mathbf{Z}\left(\!\left(x_n\right)\!\right)\right)\left(\mathbf{Z}\left(\!\left(y_n\right)\!\right)\right)}_{\text{series multiplication}} \qquad \forall \left(\!\left(x_n\right)\!\right)_{n\in\mathbb{Z}}, \left(\!\left(y_n\right)\!\right)_{n\in\mathbb{Z}} \in \ell^2_{\mathbb{R}}$$

## 1.6 Discrete Time Fourier Transform

**Definition 1.42** The **discrete-time Fourier transform $\check{\mathbf{F}}$** of $\left(\!\left(x_n\right)\!\right)_{n\in\mathbb{Z}}$ is defined as

$$\left[\check{\mathbf{F}}\left(\!\left(x_n\right)\!\right)\right](\omega) \triangleq \sum_{n\in\mathbb{Z}} x_n e^{-i\omega n} \qquad \forall \left(\!\left(x_n\right)\!\right)_{n\in\mathbb{Z}} \in \ell^2_{\mathbb{R}}$$

If we compare the definition of the *Discrete Time Fourier Transform* (Definition 1.42 page 12) to the definition of the Z-transform (Definition 1.39 page 11), we see that the DTFT is just a special case of the more general Z-Transform, with $z = e^{i\omega}$. If we imagine $z \in \mathbb{C}$ as a complex plane, then $e^{i\omega}$ is a unit circle in this plane. The "frequency" $\omega$ in the DTFT is the unit circle in the much larger z-plane as illustrated in Figure 1 (page 12).

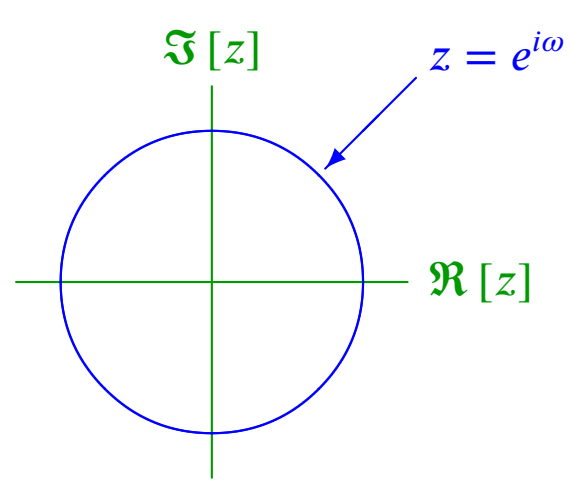

Figure 1: Unit circle in complex-z plane

**Proposition 1.43** [27] *Let $\breve{x}(\omega) \triangleq \check{\mathbf{F}}\left[\left(\!\left(x_n\right)\!\right)\right](\omega)$ be the* DISCRETE-TIME FOURIER TRANSFORM *(Definition 1.42 page 12) of a sequence $\left(\!\left(x_n\right)\!\right)_{n\in\mathbb{Z}}$ in $\ell^2_{\mathbb{R}}$.*

$$\underbrace{\breve{x}(\omega) = \breve{x}(\omega + 2\pi n)}_{\text{PERIODIC with period } 2\pi} \qquad \forall n \in \mathbb{Z}$$

[25] 📚 [50], page 344, ⟨Proposition J.1⟩
[26] 📚 [50], pages 344–345, ⟨Theorem J.1⟩
[27] 📚 [50], pages 348–349, ⟨Proposition J.2⟩

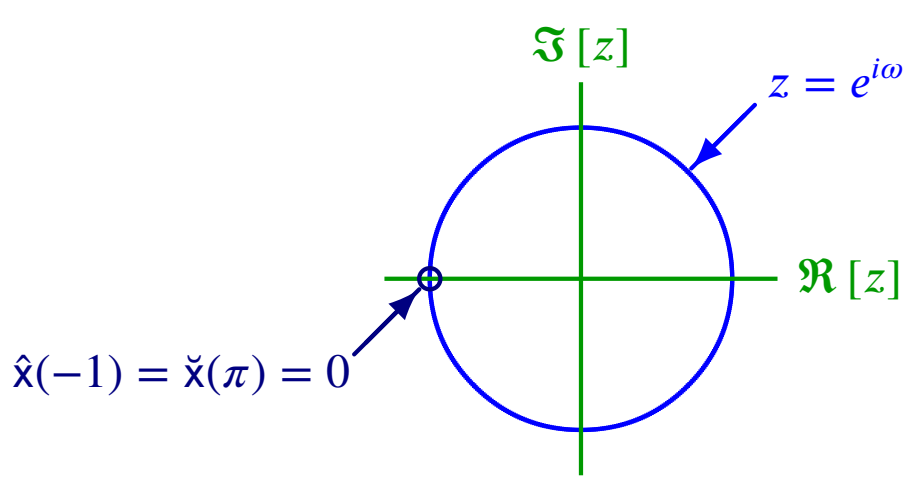

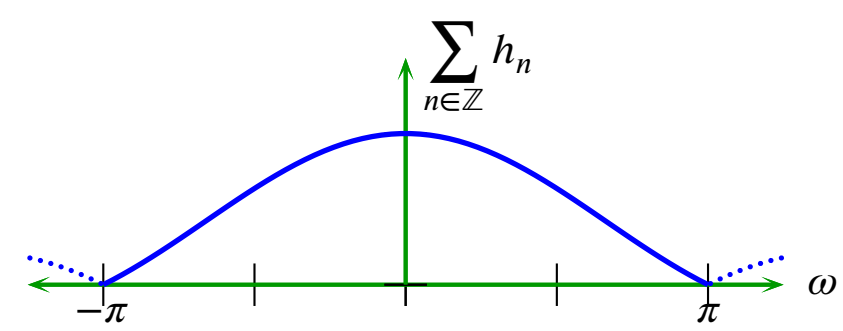

**Proposition 1.44** [28] *Let $\hat{x}(z)$ be the* Z-TRANSFORM (Definition 1.39 page 11) *and* $\breve{x}(\omega)$ *the* DISCRETE-TIME FOURIER TRANSFORM (Definition 1.42 page 12) *of* $(\!(x_n)\!)$.

$$\underbrace{\left\{\sum_{n\in\mathbb{Z}} x_n = c\right\}}_{\text{(1) time domain}} \iff \underbrace{\left\{\hat{x}(z)\Big|_{z=1} = c\right\}}_{\text{(2) z domain}} \iff \underbrace{\left\{\breve{x}(\omega)\Big|_{\omega=0} = c\right\}}_{\text{(3) frequency domain}} \qquad \forall (\!(x_n)\!)_{n\in\mathbb{Z}} \in \ell^2_{\mathbb{R}}, c\in\mathbb{R}$$

**Proposition 1.45** [29]

$$\underbrace{\sum_{n\in\mathbb{Z}} (-1)^n x_n = c}_{\text{(1) in "time"}} \iff \underbrace{\hat{x}(z)|_{z=-1} = c}_{\text{(2) in "z domain"}} \iff \underbrace{\breve{x}(\omega)|_{\omega=\pi} = c}_{\text{(3) in "frequency"}}$$

$$\iff \underbrace{\left(\sum_{n\in\mathbb{Z}} h_{2n}, \sum_{n\in\mathbb{Z}} h_{2n+1}\right) = \left(\frac{1}{2}\left(\sum_{n\in\mathbb{Z}} h_n + c\right), \frac{1}{2}\left(\sum_{n\in\mathbb{Z}} h_n - c\right)\right)}_{\text{(4) sum of even, sum of odd}}$$

$\forall c\in\mathbb{R}, (\!(x_n)\!)_{n\in\mathbb{Z}}, (\!(y_n)\!)_{n\in\mathbb{Z}} \in \ell^2_{\mathbb{R}}$

**Lemma 1.46** [30] *Let $\tilde{f}(\omega)$ be the* DTFT (Definition 1.42 page 12) *of a sequence* $(\!(x_n)\!)_{n\in\mathbb{Z}}$.

$$\underbrace{(\!(x_n \in \mathbb{R})\!)_{n\in\mathbb{Z}}}_{\text{REAL-VALUED sequence}} \implies \underbrace{|\breve{x}(\omega)|^2 = |\breve{x}(-\omega)|^2}_{\text{EVEN}} \qquad \forall (\!(x_n)\!)_{n\in\mathbb{Z}} \in \ell^2_{\mathbb{R}}$$

**Theorem 1.47** (inverse DTFT) [31] *Let $\breve{x}(\omega)$ be the* DISCRETE-TIME FOURIER TRANSFORM (Definition 1.42 page 12) *of a sequence* $(\!(x_n)\!)_{n\in\mathbb{Z}} \in \ell^2_{\mathbb{R}}$. *Let $\tilde{x}^{-1}$ be the inverse of $\breve{x}$.*

$$\underbrace{\left\{\breve{x}(\omega) \triangleq \sum_{n\in\mathbb{Z}} x_n e^{-i\omega n}\right\}}_{\breve{x}(\omega) \triangleq \breve{\mathbf{F}}(\!(x_n)\!)} \implies \underbrace{\left\{x_n = \frac{1}{2\pi}\int_{\alpha-\pi}^{\alpha+\pi} \breve{x}(\omega) e^{i\omega n}\, d\omega \quad \forall \alpha\in\mathbb{R}\right\}}_{(\!(x_n)\!) = \breve{\mathbf{F}}^{-1}\breve{\mathbf{F}}(\!(x_n)\!)} \qquad \forall (\!(x_n)\!)_{n\in\mathbb{Z}} \in \ell^2_{\mathbb{R}}$$

---

[28] [50], pages 349–350, ⟨Proposition J.3⟩
[29] [25], page 123
[30] [50], pages 352–353, ⟨Lemma J.2⟩
[31] [68], page 3-95, ⟨(3.6.2)⟩

**Theorem 1.48** (orthonormal quadrature conditions) [32] *Let $\check{x}(\omega)$ be the* DISCRETE-TIME FOURIER TRANSFORM *(Definition 1.42 page 12) of a sequence $(\!(x_n)\!)_{n\in\mathbb{Z}} \in \ell^2_{\mathbb{R}}$. Let $\bar{\delta}_n$ be the* KRONECKER DELTA FUNCTION *at $n$ (Definition 6.1 page 58).*

$$\begin{array}{rclcrcll} \sum_{m\in\mathbb{Z}} x_m y^*_{m-2n} &=& 0 & \iff & \check{x}(\omega)\check{y}^*(\omega) + \check{x}(\omega+\pi)\check{y}^*(\omega+\pi) &=& 0 & \forall n\in\mathbb{Z}, \forall (\!(x_n)\!),(\!(y_n)\!)\in\ell^2_{\mathbb{R}} \\ \sum_{m\in\mathbb{Z}} x_m x^*_{m-2n} &=& \bar{\delta}_n & \iff & |\check{x}(\omega)|^2 + |\check{x}(\omega+\pi)|^2 &=& 2 & \forall n\in\mathbb{Z}, \forall (\!(x_n)\!),(\!(y_n)\!)\in\ell^2_{\mathbb{R}} \end{array}$$

# 2 Background: transversal operators

## 2.1 Definitions

Much of B-spline and wavelet theory can be constructed with the help of the **translation operator T** and the **dilation operator D** (next).

**Definition 2.1** [33]

1. **T** is the **translation operator** on $\mathbb{C}^{\mathbb{C}}$ defined as
$$\mathbf{T}_\tau \mathsf{f}(x) \triangleq \mathsf{f}(x-\tau) \quad \text{and} \quad \mathbf{T} \triangleq \mathbf{T}_1 \qquad \forall \mathsf{f}\in\mathbb{C}^{\mathbb{C}}$$
2. **D** is the **dilation operator** on $\mathbb{C}^{\mathbb{C}}$ defined as
$$\mathbf{D}_\alpha \mathsf{f}(x) \triangleq \mathsf{f}(\alpha x) \quad \text{and} \quad \mathbf{D} \triangleq \sqrt{2}\mathbf{D}_2 \qquad \forall \mathsf{f}\in\mathbb{C}^{\mathbb{C}}$$

**Example 2.2**

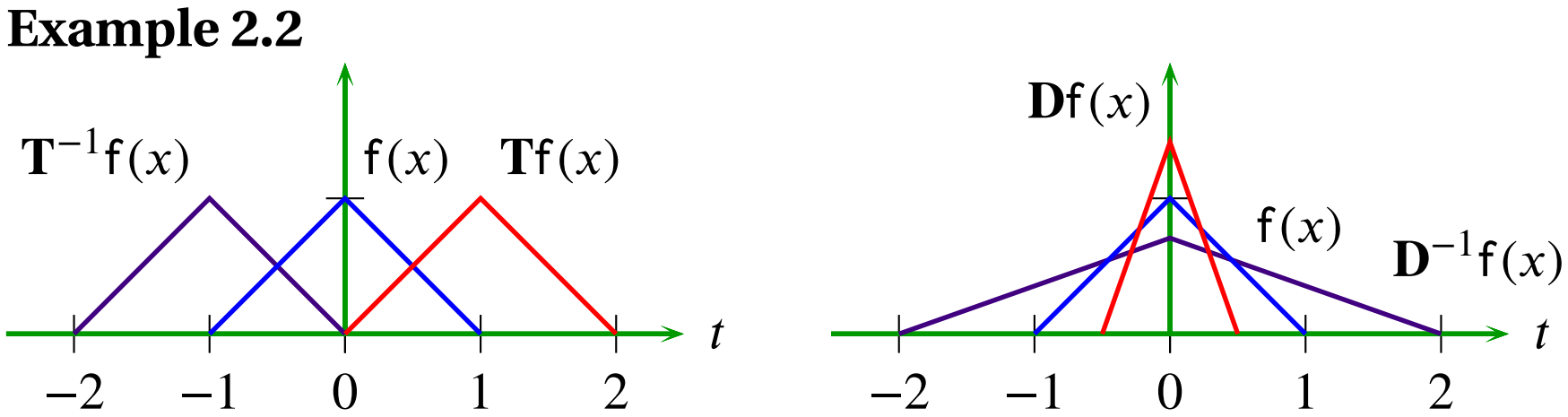

---

[32] [30], pages 132–137, ⟨(5.1.20),(5.1.39)⟩

[33] [112], pages 79–80, ⟨Definition 3.39⟩, [21], pages 41–42, [116], page 18, ⟨Definitions 2.3,2.4⟩, [69], page A-21, [11], page 473, [90], page 260, [14], page , [58], page 250, ⟨Notation 9.4⟩, [19], page 74, [45], page 639, [29], page 81, [28], page 2, [50], page 2

## 2.2 Properties

### 2.2.1 Algebraic properties

**Proposition 2.3** [34] *Let* $\mathbf{T}$ *be the* TRANSLATION OPERATOR *(Definition 2.1 page 14).*

$$\sum_{n\in\mathbb{Z}} \mathbf{T}^n \mathsf{f}(x) = \sum_{n\in\mathbb{Z}} \mathbf{T}^n \mathsf{f}(x+1) \qquad \forall \mathsf{f}\in\mathbb{R}^{\mathbb{R}} \qquad \left(\sum_{n\in\mathbb{Z}} \mathbf{T}^n \mathsf{f}(x) \text{ is PERIODIC with period } 1\right)$$

In a linear space, every operator has an *inverse*. Although the inverse always exists as a relation, it may not exist as a function or as an operator. But in some cases the inverse of an operator is itself an operator. The inverses of the operators $\mathbf{T}$ and $\mathbf{D}$ both exists as operators, as demonstrated by Proposition 2.4 (next).

**Proposition 2.4** [35] *Let* $\mathbf{T}$ *and* $\mathbf{D}$ *be as defined in Definition 2.1 page 14.*

$\mathbf{T}$ *has an inverse* $\mathbf{T}^{-1}$ *in* $\mathbb{C}^{\mathbb{C}}$ *expressed by the relation*

$$\mathbf{T}^{-1}\mathsf{f}(x) = \mathsf{f}(x+1) \qquad \forall \mathsf{f}\in\mathbb{C}^{\mathbb{C}} \qquad \text{(translation operator inverse)}.$$

$\mathbf{D}$ *has an inverse* $\mathbf{D}^{-1}$ *in* $\mathbb{C}^{\mathbb{C}}$ *expressed by the relation*

$$\mathbf{D}^{-1}\mathsf{f}(x) = \frac{\sqrt{2}}{2}\mathsf{f}\left(\frac{1}{2}x\right) \qquad \forall \mathsf{f}\in\mathbb{C}^{\mathbb{C}} \qquad \text{(dilation operator inverse)}.$$

**Proposition 2.5** [36] *Let* $\mathbf{T}$ *and* $\mathbf{D}$ *be as defined in Definition 2.1 page 14. Let* $\mathbf{D}^0 = \mathbf{T}^0 \triangleq \mathbf{I}$ *be the* IDENTITY OPERATOR.

$$\mathbf{D}^j\mathbf{T}^n\mathsf{f}(x) = 2^{j/2}\mathsf{f}\left(2^j x - n\right) \qquad \forall j,n\in\mathbb{Z},\ \mathsf{f}\in\mathbb{C}^{\mathbb{C}}$$

### 2.2.2 Linear space properties

**Definition 2.6** [37] Let $+$ be an addition operator on a tuple $\left(x_n\right)_m^N$.

The **summation** of $\left(x_n\right)$ from index $m$ to index $N$ with respect to $+$ is

$$\sum_{n=m}^{N} x_n \triangleq \begin{cases} 0 & \text{for } N < m \\ \left(\sum_{n=m}^{N-1} x_n\right) + x_N & \text{for } N \geq m \end{cases}$$

An infinite summation $\sum_{n=1}^{\infty}\phi_n$ is meaningless outside some topological space (e.g. metric space, normed space, etc.). The sum $\sum_{n=1}^{\infty}\phi_n$ is an abbreviation for $\lim_{N\to\infty}\sum_{n=1}^{N}\phi_n$ (the limit of partial sums). And the concept of limit is also itself meaningless outside of a topological space.

[34] [50], page 3
[35] [50], page 3
[36] [50], page 4
[37] [15], page 8, ⟨Definition I.3.1⟩, [40], page 280, ⟨"$\sum$" notation⟩

**Definition 2.7** [38] Let $(X, \boldsymbol{T})$ be a topological space and $\lim$ be the limit induced by the topology $\boldsymbol{T}$.

$$\sum_{n=1}^{\infty} \boldsymbol{x}_n \triangleq \sum_{n\in\mathbb{N}} \boldsymbol{x}_n \triangleq \lim_{N\to\infty} \sum_{n=1}^{N} \boldsymbol{x}_n$$

$$\sum_{n=-\infty}^{\infty} \boldsymbol{x}_n \triangleq \sum_{n\in\mathbb{Z}} \boldsymbol{x}_n \triangleq \lim_{N\to\infty} \left( \sum_{n=0}^{N} \boldsymbol{x}_n \right) + \left( \lim_{N\to-\infty} \sum_{n=-1}^{N} \boldsymbol{x}_n \right)$$

In general the operators $\mathbf{T}$ and $\mathbf{D}$ are *noncommutative* ($\mathbf{TD} \neq \mathbf{DT}$), as demonstrated by Proposition 2.9 and by the following illustration.

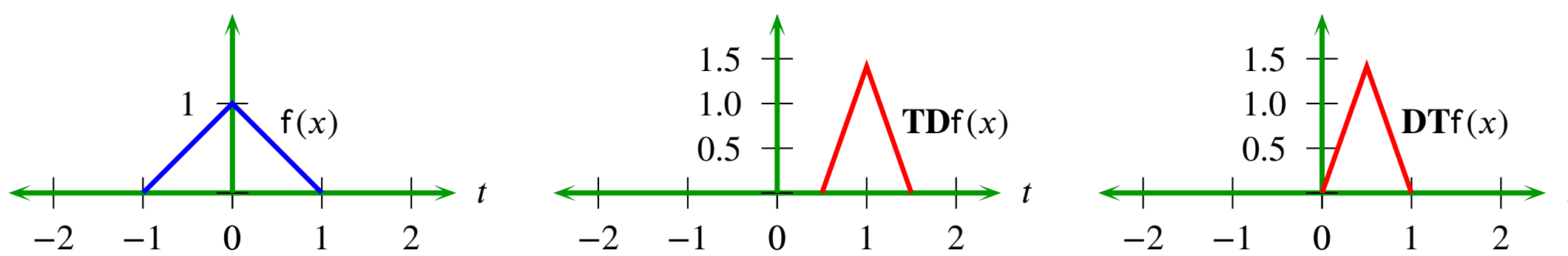

**Proposition 2.8** *Let* $\mathbf{T}$ *and* $\mathbf{D}$ *be as in Definition 2.1 page 14.*

$$\mathbf{D}^j\mathbf{T}^n\left[\mathsf{f}\mathsf{g}\right] = 2^{-j/2}\left[\mathbf{D}^j\mathbf{T}^n\mathsf{f}\right]\left[\mathbf{D}^j\mathbf{T}^n\mathsf{g}\right] \qquad \forall j,n\in\mathbb{Z},\ \mathsf{f}\in\mathbb{C}^{\mathbb{C}}$$

**Proposition 2.9** (commutator relation) [39] *Let* $\mathbf{T}$ *and* $\mathbf{D}$ *be as in Definition 2.1 page 14.*

$$\begin{aligned} \mathbf{D}^j\mathbf{T}^n &= \mathbf{T}^{2^{-j/2}n}\mathbf{D}^j && \forall j,n\in\mathbb{Z} \\ \mathbf{T}^n\mathbf{D}^j &= \mathbf{D}^j\mathbf{T}^{2^j n} && \forall n,j\in\mathbb{Z} \end{aligned}$$

### 2.2.3 Inner-product space properties

In an inner product space, every operator has an *adjoint* and this adjoint is always itself an operator. In the case where the adjoint and inverse of an operator $\mathbf{U}$ coincide, then $\mathbf{U}$ is said to be *unitary*. And in this case, $\mathbf{U}$ has several nice properties (see Proposition 2.15 and Theorem 2.16 page 18). Proposition 2.10 (next) gives the adjoints of $\mathbf{D}$ and $\mathbf{T}$, and Proposition 2.11 (page 17) demonstrates that both $\mathbf{D}$ and $\mathbf{T}$ are unitary. Other examples of unitary operators include the *Fourier Transform operator* $\tilde{\mathbf{F}}$ and the *rotation matrix operator*.

**Proposition 2.10** *Let* $\mathbf{T}$ *be the translation operator* (Definition 2.1 page 14) *with adjoint* $\mathbf{T}^*$ *and* $\mathbf{D}$ *the dilation operator with adjoint* $\mathbf{D}^*$.

$$\begin{aligned} \mathbf{T}^*\mathsf{f}(x) &= \mathsf{f}(x+1) && \forall \mathsf{f}\in\boldsymbol{L}^2_{\mathbb{R}} && \text{(translation operator adjoint)} \\ \mathbf{D}^*\mathsf{f}(x) &= \frac{\sqrt{2}}{2}\mathsf{f}\left(\frac{1}{2}x\right) && \forall \mathsf{f}\in\boldsymbol{L}^2_{\mathbb{R}} && \text{(dilation operator adjoint)} \end{aligned}$$

---
[38] [72], page 4, [76], page 43, [10], pages 3–4
[39] [21], page 42, ⟨equation (2.9)⟩, [28], page 21, [45], page 641, [46], page 110

**Proposition 2.11** [40] *Let* $\mathbf{T}$ *and* $\mathbf{D}$ *be as in Definition 2.1 page 14. Let* $\mathbf{T}^{-1}$ *and* $\mathbf{D}^{-1}$ *be as in Proposition 2.4 page 15.*

$$\begin{array}{ll} \mathbf{T} \text{ is UNITARY in } \boldsymbol{L}^2_{\mathbb{R}} & (\mathbf{T}^{-1} = \mathbf{T}^* \text{ in } \boldsymbol{L}^2_{\mathbb{R}}). \\ \mathbf{D} \text{ is UNITARY in } \boldsymbol{L}^2_{\mathbb{R}} & (\mathbf{D}^{-1} = \mathbf{D}^* \text{ in } \boldsymbol{L}^2_{\mathbb{R}}). \end{array}$$

### 2.2.4 Normed linear space properties

**Proposition 2.12** *Let* $\mathbf{D}$ *be the* DILATION OPERATOR *(Definition 2.1 page 14).*

$$\left\{\begin{array}{lll} (1). & \mathbf{D}\mathsf{f}(x) = \sqrt{2}\mathsf{f}(x) & \text{and} \\ (2). & \mathsf{f}(x) \text{ is CONTINUOUS} & \end{array}\right\} \iff \{\mathsf{f}(x) \text{ is a CONSTANT}\} \qquad \forall \mathsf{f} \in \boldsymbol{L}^2_{\mathbb{R}}$$

Note that in Proposition 2.12, it is not possible to remove the *continuous* constraint outright (next two counterexamples).

**Counterexample 2.13** Let $\mathsf{f}(x)$ be a function in $\mathbb{R}^{\mathbb{R}}$.

Let $\mathsf{f}(x) \triangleq \begin{cases} 0 & \text{for } x = 0 \\ 1 & \text{otherwise.} \end{cases}$

Then $\mathbf{D}\mathsf{f}(x) \triangleq \sqrt{2}\mathsf{f}(2x) = \sqrt{2}\mathsf{f}(x)$, but $\mathsf{f}(x)$ is *not constant*.

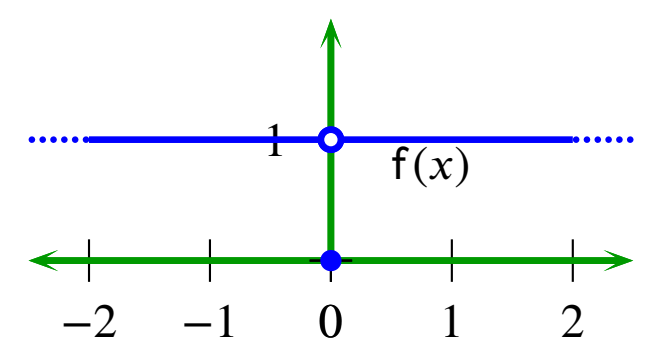

**Counterexample 2.14** Let $\mathsf{f}(x)$ be a function in $\mathbb{R}^{\mathbb{R}}$. Let $\mathbb{Q}$ be the set of *rational numbers* and $\mathbb{R}\backslash\mathbb{Q}$ the set of *irrational numbers*.

Let $\mathsf{f}(x) \triangleq \begin{cases} 1 & \text{for } x \in \mathbb{Q} \\ -1 & \text{for } x \in \mathbb{R}\backslash\mathbb{Q}. \end{cases}$

Then $\mathbf{D}\mathsf{f}(x) \triangleq \sqrt{2}\mathsf{f}(2x) = \sqrt{2}\mathsf{f}(x)$, but $\mathsf{f}(x)$ is *not constant*.

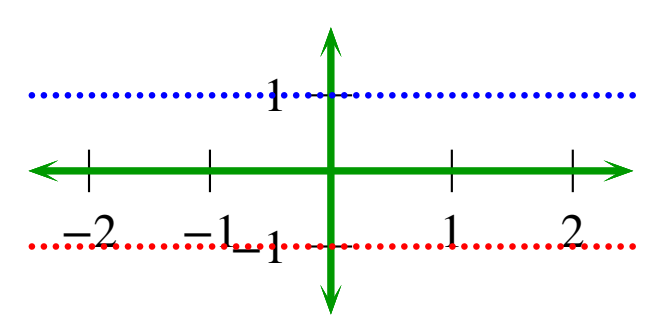

**Proposition 2.15** (Operator norm) *Let* $\mathbf{T}$ *and* $\mathbf{D}$ *be as in Definition 2.1 page 14. Let* $\mathbf{T}^{-1}$ *and* $\mathbf{D}^{-1}$ *be as in Proposition 2.4 page 15. Let* $\mathbf{T}^*$ *and* $\mathbf{D}^*$ *be as in Proposition 2.10 page 16. Let* $\|\cdot\|$ *and* $\langle \triangle \mid \nabla \rangle$ *be as in Definition 1.2 page 3. Let* $|||\cdot|||$ *be the operator norm induced by* $\|\cdot\|$.

$$|||\mathbf{T}||| = |||\mathbf{D}||| = \left|\left|\left|\mathbf{T}^*\right|\right|\right| = \left|\left|\left|\mathbf{D}^*\right|\right|\right| = \left|\left|\left|\mathbf{T}^{-1}\right|\right|\right| = \left|\left|\left|\mathbf{D}^{-1}\right|\right|\right| = 1$$

✎PROOF: These results follow directly from the fact that $\mathbf{T}$ and $\mathbf{D}$ are *unitary* and from properties of unitary operators. ✏

[40] [21], page 41, ⟨Lemma 2.5.1⟩, [116], page 18, ⟨Lemma 2.5⟩

**Theorem 2.16** *Let* $\mathbf{T}$ *and* $\mathbf{D}$ *be as in Definition 2.1 page 14. Let* $\mathbf{T}^{-1}$ *and* $\mathbf{D}^{-1}$ *be as in Proposition 2.4 page 15. Let* $\|\cdot\|$ *and* $\langle \triangle \mid \nabla \rangle$ *be as in Definition 1.2 page 3.*

1. $\|\mathbf{T}f\| = \|\mathbf{D}f\| = \|f\| \quad \forall f \in L^2_{\mathbb{R}}$ (ISOMETRIC IN LENGTH)
2. $\|\mathbf{T}f - \mathbf{T}g\| = \|\mathbf{D}f - \mathbf{D}g\| = \|f - g\| \quad \forall f,g \in L^2_{\mathbb{R}}$ (ISOMETRIC IN DISTANCE)
3. $\left\|\mathbf{T}^{-1}f - \mathbf{T}^{-1}g\right\| = \left\|\mathbf{D}^{-1}f - \mathbf{D}^{-1}g\right\| = \|f - g\| \quad \forall f,g \in L^2_{\mathbb{R}}$ (ISOMETRIC IN DISTANCE)
4. $\langle \mathbf{T}f \mid \mathbf{T}g \rangle = \langle \mathbf{D}f \mid \mathbf{D}g \rangle = \langle f \mid g \rangle \quad \forall f,g \in L^2_{\mathbb{R}}$ (SURJECTIVE)
5. $\left\langle \mathbf{T}^{-1}f \,\middle|\, \mathbf{T}^{-1}g \right\rangle = \left\langle \mathbf{D}^{-1}f \,\middle|\, \mathbf{D}^{-1}g \right\rangle = \langle f \mid g \rangle \quad \forall f,g \in L^2_{\mathbb{R}}$ (SURJECTIVE)

✎PROOF: These results follow directly from the fact that $\mathbf{T}$ and $\mathbf{D}$ are *unitary* (Proposition 2.11 page 17) and from properties of unitary operators. ✏

**Proposition 2.17** *Let* $\mathbf{T}$ *be as in Definition 2.1 page 14. Let* $\mathbf{A}^*$ *be the adjoint of an operator* $\mathbf{A}$.

$$\left(\sum_{n\in\mathbb{Z}} \mathbf{T}^n\right) = \left(\sum_{n\in\mathbb{Z}} \mathbf{T}^n\right)^* \qquad \left(\text{The operator } \left[\sum_{n\in\mathbb{Z}} \mathbf{T}^n\right] \text{ is SELF-ADJOINT}\right)$$

### 2.2.5 Fourier transform properties

**Proposition 2.18** *Let* $\mathbf{T}$ *and* $\mathbf{D}$ *be as in Definition 2.1 page 14. Let* $\mathbf{B}$ *be the* TWO-SIDED LAPLACE TRANSFORM *defined as*

$$[\mathbf{B}f](s) \triangleq \frac{1}{\sqrt{2\pi}} \int_{\mathbb{R}} f(x) e^{-sx} \, dx .$$

1. $\mathbf{B}\mathbf{T}^n = e^{-sn}\mathbf{B} \quad \forall n\in\mathbb{Z}$
2. $\mathbf{B}\mathbf{D}^j = \mathbf{D}^{-j}\mathbf{B} \quad \forall j\in\mathbb{Z}$
3. $\mathbf{D}\mathbf{B} = \mathbf{B}\mathbf{D}^{-1} \quad \forall n\in\mathbb{Z}$
4. $\mathbf{B}\mathbf{D}^{-1}\mathbf{B}^{-1} = \mathbf{B}^{-1}\mathbf{D}^{-1}\mathbf{B} = \mathbf{D} \quad \forall n\in\mathbb{Z}$ ($\mathbf{D}^{-1}$ *is* SIMILAR *to* $\mathbf{D}$)
5. $\mathbf{D}\mathbf{B}\mathbf{D} = \mathbf{D}^{-1}\mathbf{B}\mathbf{D}^{-1} = \mathbf{B} \quad \forall n\in\mathbb{Z}$

**Corollary 2.19** *Let* $\mathbf{T}$ *and* $\mathbf{D}$ *be as in Definition 2.1 page 14. Let* $\tilde{f}(\omega) \triangleq \tilde{\mathbf{F}}f(x)$ *be the* FOURIER TRANSFORM *(Definition 1.26 page 8) of some function* $f \in L^2_{\mathbb{R}}$ *(Definition 1.2 page 3).*

1. $\tilde{\mathbf{F}}\mathbf{T}^n = e^{-i\omega n}\tilde{\mathbf{F}}$
2. $\tilde{\mathbf{F}}\mathbf{D}^j = \mathbf{D}^{-j}\tilde{\mathbf{F}}$
3. $\mathbf{D}\tilde{\mathbf{F}} = \tilde{\mathbf{F}}\mathbf{D}^{-1}$
4. $\mathbf{D} = \tilde{\mathbf{F}}\mathbf{D}^{-1}\tilde{\mathbf{F}}^{-1} = \tilde{\mathbf{F}}^{-1}\mathbf{D}^{-1}\tilde{\mathbf{F}}$
5. $\tilde{\mathbf{F}} = \mathbf{D}\tilde{\mathbf{F}}\mathbf{D} = \mathbf{D}^{-1}\tilde{\mathbf{F}}\mathbf{D}^{-1}$

✎PROOF: These results follow directly from Proposition 2.18 page 18. ✏

**Proposition 2.20** *Let* $\mathbf{T}$ *and* $\mathbf{D}$ *be as in Definition 2.1 page 14. Let* $\tilde{\mathsf{f}}(\omega) \triangleq \tilde{\mathbf{F}}\mathsf{f}(x)$ *be the* FOURIER TRANSFORM (Definition 1.26 page 8) *of some function* $\mathsf{f} \in \boldsymbol{L}^2_{\mathbb{R}}$ (Definition 1.2 page 3).

$$\tilde{\mathbf{F}}\mathbf{D}^j\mathbf{T}^n\mathsf{f}(x) = \frac{1}{2^{j/2}} e^{-i\frac{\omega}{2^j}n}\tilde{\mathsf{f}}\left(\frac{\omega}{2^j}\right)$$

**Proposition 2.21** *Let* $\mathbf{T}$ *be the translation operator* (Definition 2.1 page 14). *Let* $\tilde{\mathsf{f}}(\omega) \triangleq \tilde{\mathbf{F}}\mathsf{f}(x)$ *be the* FOURIER TRANSFORM (Definition 1.26 page 8) *of a function* $\mathsf{f} \in \boldsymbol{L}^2_{\mathbb{R}}$. *Let* $\breve{\mathsf{a}}(\omega)$ *be the* DTFT (Definition 1.42 page 12) *of a sequence* $(\!(a_n)\!)_{n\in\mathbb{Z}} \in \boldsymbol{\ell}^2_{\mathbb{R}}$ (Definition 1.37 page 11).

$$\tilde{\mathbf{F}}\sum_{n\in\mathbb{Z}} a_n\mathbf{T}^n\phi(x) = \breve{\mathsf{a}}(\omega)\tilde{\phi}(\omega) \qquad \forall (\!(a_n)\!)\in\boldsymbol{\ell}^2_{\mathbb{R}}, \phi(x)\in\boldsymbol{L}^2_{\mathbb{R}}$$

**Theorem 2.22** (Poisson Summation Formula—PSF) [41] *Let* $\tilde{\mathsf{f}}(\omega)$ *be the* FOURIER TRANSFORM (Definition 1.26 page 8) *of a function* $\mathsf{f}(x) \in \boldsymbol{L}^2_{\mathbb{R}}$.

$$\underbrace{\sum_{n\in\mathbb{Z}}\mathbf{T}^n_\tau\mathsf{f}(x) = \sum_{n\in\mathbb{Z}}\mathsf{f}(x+n\tau)}_{\text{summation in “time”}} = \underbrace{\sqrt{\frac{2\pi}{\tau}}\,\hat{\mathbf{F}}^{-1}\mathbf{S}\tilde{\mathbf{F}}[\mathsf{f}(x)]}_{\text{operator notation}} = \underbrace{\frac{\sqrt{2\pi}}{\tau}\sum_{n\in\mathbb{Z}}\tilde{\mathsf{f}}\left(\frac{2\pi}{\tau}n\right)e^{i\frac{2\pi}{\tau}nx}}_{\text{summation in “frequency”}}$$

*where* $\mathbf{S} \in \boldsymbol{\ell}^{2\,\boldsymbol{L}^2_{\mathbb{R}}}_{\mathbb{R}}$ *is the* SAMPLING OPERATOR *defined as*

$$[\mathbf{S}\mathsf{f}(x)](n) \triangleq \mathsf{f}\left(\tfrac{2\pi}{\tau}n\right) \qquad \forall \mathsf{f}\in\boldsymbol{L}^2_{(\mathbb{R},\mathscr{B},\mu)},\ \tau\in\mathbb{R}^+$$

**Theorem 2.23** (Inverse Poisson Summation Formula—IPSF) [42] *Let* $\tilde{\mathsf{f}}(\omega)$ *be the* FOURIER TRANSFORM (Definition 1.26 page 8) *of a function* $\mathsf{f}(x) \in \boldsymbol{L}^2_{\mathbb{R}}$.

$$\underbrace{\sum_{n\in\mathbb{Z}}\mathbf{T}^n_{2\pi/\tau}\tilde{\mathsf{f}}(\omega) \triangleq \sum_{n\in\mathbb{Z}}\tilde{\mathsf{f}}\left(\omega-\frac{2\pi}{\tau}n\right)}_{\text{summation in “frequency”}} = \underbrace{\frac{\tau}{\sqrt{2\pi}}\sum_{n\in\mathbb{Z}}\mathsf{f}(n\tau)e^{-i\omega n\tau}}_{\text{summation in “time”}}$$

**Remark 2.24** The left hand side of the *Poisson Summation Formula* (Theorem 2.22 page 19) is very similar to the *Zak Transform* $\mathbf{Z}$: [43]

$$(\mathbf{Z}\mathsf{f})(t,\omega) \triangleq \sum_{n\in\mathbb{Z}}\mathsf{f}(x+n\tau)e^{i2\pi n\omega}$$

**Remark 2.25** A generalization of the *Poisson Summation Formula* (Theorem 2.22 page 19) is the **Selberg Trace Formula**.[44]

[41] 📚 [6], page 624, 📚 [73], page 389, 📚 [80], page 254, 📚 [98], pages 194–195, 📚 [39], page 337, 📚 [50], pages 12–13, ⟨Theorem 1.2⟩
[42] 📚 [50], pages 14–15, ⟨Theorem 1.3⟩, 📚 [44], page 88,
[43] 📚 [64], page 24, 📚 [117], page 482
[44] 📚 [81], page 349, 📚 [100], 📚 [107]

**Lemma 2.26** [45] *Let* $\boldsymbol{\Omega} \triangleq \left(X, +, \cdot, (\mathbb{F}, \dot{+}, \dot{\times}), \boldsymbol{T}\right)$ *be a topological linear space. Let* $\operatorname{span} A$ *be the* SPAN *of a set* $A$. *Let* $\tilde{\mathsf{f}}(\omega)$ *and* $\tilde{\mathsf{g}}(\omega)$ *be the* FOURIER TRANSFORMS (Definition 1.26 page 8) *of the functions* $\mathsf{f}(x)$ *and* $\mathsf{g}(x)$, *respectively, in* $\boldsymbol{L}^2_{\mathbb{R}}$ (Definition 1.2 page 3). *Let* $\breve{a}(\omega)$ *be the* DTFT (Definition 1.42 page 12) *of a sequence* $(\!(a_n)\!)_{n\in\mathbb{Z}}$ *in* $\boldsymbol{\ell}^2_{\mathbb{R}}$ (Definition 1.37 page 11).

$$\left\{\begin{array}{ll} (1). & \left\{\mathbf{T}^n\mathsf{f}|_{n\in\mathbb{Z}}\right\} \text{ is a SCHAUDER BASIS for } \boldsymbol{\Omega} \quad \text{and} \\ (2). & \left\{\mathbf{T}^n\mathsf{g}|_{n\in\mathbb{Z}}\right\} \text{ is a SCHAUDER BASIS for } \boldsymbol{\Omega} \end{array}\right\} \implies \left\{\begin{array}{l} \exists\,(\!(a_n)\!)_{n\in\mathbb{Z}} \quad \text{such that} \\ \tilde{\mathsf{f}}(\omega) = \breve{a}(\omega)\tilde{\mathsf{g}}(\omega) \end{array}\right\}$$

**Theorem 2.27** (Battle-Lemarié orthogonalization) [46] *Let* $\tilde{\mathsf{f}}(\omega)$ *be the* FOURIER TRANSFORM (Definition 1.26 page 8) *of a function* $\mathsf{f} \in \boldsymbol{L}^2_{\mathbb{R}}$.

$$\left\{\begin{array}{ll} 1. & \left\{\mathbf{T}^n\mathsf{g}|_{n\in\mathbb{Z}}\right\} \text{ is a RIESZ BASIS for } \boldsymbol{L}^2_{\mathbb{R}} \quad \text{and} \\ 2. & \tilde{\mathsf{f}}(\omega) \triangleq \dfrac{\tilde{\mathsf{g}}(\omega)}{\sqrt{2\pi\sum_{n\in\mathbb{Z}}|\tilde{\mathsf{g}}(\omega+2\pi n)|^2}} \end{array}\right\} \implies \left\{\begin{array}{l} \left\{\mathbf{T}^n\mathsf{f}|_{n\in\mathbb{Z}}\right\} \\ \text{is an ORTHONORMAL BASIS for } \boldsymbol{L}^2_{\mathbb{R}} \end{array}\right\}$$

# 3 Background: MRA-wavelet analysis

## 3.1 Orientation

In Fourier analysis, *continuous* dilations (Definition 2.1 page 14) of the *complex exponential* form a *basis* for the *space of square integrable functions* $\boldsymbol{L}^2_{\mathbb{R}}$ such that

$$\boldsymbol{L}^2_{\mathbb{R}} = \operatorname{span}\left\{\mathbf{D}_\omega e^{ix}\,|_{\omega\in\mathbb{R}}\right\}.$$

In Fourier series analysis , *discrete* dilations of the complex exponential form a basis for $\boldsymbol{L}^2_{\mathbb{R}}(0,\, 2\pi)$ such that

$$\boldsymbol{L}^2_{\mathbb{R}}(0,\, 2\pi) = \operatorname{span}\left\{\mathbf{D}_j e^{ix}\big|_{j\in\mathbb{Z}}\right\}.$$

In Wavelet analysis, for some *mother wavelet* (Definition 3.13 page 25) $\psi(x)$,

$$\boldsymbol{L}^2_{\mathbb{R}} = \operatorname{span}\left\{\mathbf{D}_\omega\mathbf{T}_\tau\psi(x)\,|\,\omega, \tau \in \mathbb{R}\right\}.$$

However, the ranges of parameters $\omega$ and $\tau$ can be much reduced to the countable set $\mathbb{Z}$ resulting in a *dyadic* wavelet basis such that for some mother wavelet $\psi(x)$,

$$\boldsymbol{L}^2_{\mathbb{R}} = \operatorname{span}\left\{\mathbf{D}^j\mathbf{T}^n\psi(x)\,|\,j, n \in \mathbb{Z}\right\}.$$

Wavelets that are both *dyadic* and *compactly supported* have the attractive feature that they can be easily implemented in hardware or software by use of the *Fast Wavelet Transform*.

[45] [30], page 140, [50], pages 22–23, ⟨Lemma 1.1⟩,
[46] [116], page 25, ⟨Remark 2.4⟩, [111], page 71, [86], page 72, [87], page 225, [30], page 140, ⟨(5.3.3)⟩, [50], pages 23–24, ⟨Theorem 1.7⟩

In 1989, Stéphane G. Mallat introduced the *Multiresolution Analysis* (MRA, Definition 3.1 page 21) method for wavelet construction. The MRA has since become the dominate wavelet construction method. Moreover, P.G. Lemarié has proved that all wavelets with *compact support* are generated by an MRA.[47]

## 3.2 Multiresolution analysis

### 3.2.1 Definition

A multiresolution analysis provides "coarse" approximations of a function in a linear space $\boldsymbol{L}^2_{\mathbb{R}}$ at multiple "scales" or "resolutions". Key to this process is a sequence of *scaling functions*. Most traditional transforms feature a single *scaling function* $\phi(x)$ set equal to one $(\phi(x) = 1)$. This allows for convenient representation of the most basic functions, such as constants.[48] A multiresolution system, on the other hand, uses a generalized form of the scaling concept:[49]

(1) Instead of the scaling function simply being set *equal to unity* $(\phi(x) = 1)$, a multiresolution analysis (Definition 3.1 page 21) is often constructed in such a way that the scaling function $\phi(x)$ forms a *partition of unity* such that $\sum_{n\in\mathbb{Z}} \mathbf{T}^n\phi(x) = 1$.

(2) Instead of there being *just one* scaling function, there is an entire sequence of scaling functions $\left(\!\left(\mathbf{D}^j\phi(x)\right)\!\right)_{j\in\mathbb{Z}}$, each corresponding to a different "*resolution*".

**Definition 3.1** [50] Let $\left(\!\left(\boldsymbol{V}_j\right)\!\right)_{j\in\mathbb{Z}}$ be a sequence of subspaces on $\boldsymbol{L}^2_{\mathbb{R}}$. Let $A^-$ be the *closure* of a set $A$. The sequence $\left(\!\left(\boldsymbol{V}_j\right)\!\right)_{j\in\mathbb{Z}}$ is a **multiresolution analysis** on $\boldsymbol{L}^2_{\mathbb{R}}$ if

1. $\boldsymbol{V}_j = \boldsymbol{V}_j^-$ $\quad\forall j\in\mathbb{Z}$ (*closed*) and
2. $\boldsymbol{V}_j \subset \boldsymbol{V}_{j+1}$ $\quad\forall j\in\mathbb{Z}$ (*linearly ordered*) and
3. $\left(\bigcup_{j\in\mathbb{Z}} \boldsymbol{V}_j\right)^- = \boldsymbol{L}^2_{\mathbb{R}}$ (*dense* in $\boldsymbol{L}^2_{\mathbb{R}}$) and
4. $\mathsf{f} \in \boldsymbol{V}_j \iff \mathbf{D}\mathsf{f} \in \boldsymbol{V}_{j+1}$ $\quad\forall j\in\mathbb{Z}, \mathsf{f}\in\boldsymbol{L}^2_{\mathbb{R}}$ (*self-similar*) and
5. $\exists\phi$ such that $\left\{\mathbf{T}^n\phi|_{n\in\mathbb{Z}}\right\}$ is a *Riesz basis* for $\boldsymbol{V}_0$.

[47] [82], [87], page 240

[48] [65], page 8

[49] The concept of a scaling space was perhaps first introduced by Taizo Iijima in 1959 in Japan, and later as the *Gaussian Pyramid* by Burt and Adelson in the 1980s in the West. [86], page 70, [62], [18], [2], [83], [5], [53], [113], ⟨historical survey⟩

[50] [59], page 44, [87], page 221, ⟨Definition 7.1⟩ , [86], page 70, [88], page 21, ⟨Definition 2.2.1⟩, [21], page 284, ⟨Definition 13.1.1⟩, [11], pages 451–452, ⟨Definition 7.7.6⟩, [112], pages 300–301, ⟨Definition 10.16⟩, [30], pages 129–140, ⟨Riesz basis: page 139⟩

A *multiresolution analysis* is also called an **MRA**. An element $V_j$ of $(\!(V_j)\!)_{j\in\mathbb{Z}}$ is a **scaling subspace** of the space $L^2_{\mathbb{R}}$. The pair $(L^2_{\mathbb{R}}, (\!(V_j)\!))$ is a **multiresolution analysis space**, or **MRA space**. The function $\phi$ is the **scaling function** of the *MRA space*.

The traditional definition of the *MRA* also includes the following:

6. $f \in V_j \iff \mathbf{T}^n f \in V_j \quad \forall n,j\in\mathbb{Z}, f\in L^2_{\mathbb{R}}$ (*translation invariant*)
7. $\bigcap_{j\in\mathbb{Z}} V_j = \{\mathbb{0}\}$ (*greatest lower bound* is **0**)

However, it can be shown that these follow from the *MRA* as defined in Definition 3.1.[51]

The MRA (Definition 3.1 page 21) is more than just an interesting mathematical toy. Under some very "reasonable" conditions (next proposition), as $j \to \infty$, the *scaling subspace* $V_j$ is *dense* in $L^2_{\mathbb{R}}$ ...meaning that with the MRA we can represent any "reasonable" function to within an arbitrary accuracy.

**Proposition 3.2** [52]

$$\left\{\begin{array}{lll} (1). & (\!(\mathbf{T}^n\phi)\!) \text{ is a } \textsc{Riesz sequence} & \text{and} \\ (2). & \tilde{\phi}(\omega) \text{ is } \textsc{continuous} \text{ at } 0 & \text{and} \\ (3). & \tilde{\phi}(0) \neq 0 & \end{array}\right\} \implies \left\{ \left(\bigcup_{j\in\mathbb{Z}} V_j\right)^- = L^2_{\mathbb{R}} \quad (\textsc{dense} \text{ in } L^2_{\mathbb{R}}) \right\}$$

### 3.2.2 Dilation equation

The scaling function $\phi(x)$ (Definition 3.1 page 21) exhibits a kind of *self-similar* property. By Definition 3.1 page 21, the dilation $\mathbf{D}f$ of each vector $f$ in $V_0$ is in $V_1$. If $\{\mathbf{T}^n\phi|_{n\in\mathbb{Z}}\}$ is a basis for $V_0$, then $\{\mathbf{D}\mathbf{T}^n\phi|_{n\in\mathbb{Z}}\}$ is a basis for $V_1$, $\{\mathbf{D}^2\mathbf{T}^n\phi|_{n\in\mathbb{Z}}\}$ is a basis for $V_2$, ...; and in general $\{\mathbf{D}^j\mathbf{T}^m\phi \,|\, j\in\mathbb{Z}\}$ is a basis for $V_j$. Also, if $\phi$ is in $V_0$, then it is also in $V_1$ (because $V_0 \subset V_1$). And because $\phi$ is in $V_1$ and because $\{\mathbf{D}\mathbf{T}^n\phi|_{n\in\mathbb{Z}}\}$ is a basis for $V_1$, $\phi$ is a linear combination of the elements in $\{\mathbf{D}\mathbf{T}^n\phi|_{n\in\mathbb{Z}}\}$. That is, $\phi$ can be represented as a linear combination of translated and dilated versions of itself. The resulting equation is called the *dilation equation* (Definition 3.3, next).[53]

**Definition 3.3** [54] Let $(L^2_{\mathbb{R}}, (\!(V_j)\!))$ be a *multiresolution analysis space* with scaling function $\phi$ (Definition 3.1 page 21). Let $(\!(h_n)\!)_{n\in\mathbb{Z}}$ be a *sequence* (Definition 1.36 page 11) in $\ell^2_{\mathbb{R}}$ (Definition 1.37 page 11).

---

[51] 📖 [59], page 45, ⟨Theorem 1.6⟩, 📖 [116], pages 19–28, ⟨Proposition 2.14⟩, 📖 [93], pages 313–314, ⟨Lemma 6.4.28⟩, 📖 [50], pages 32–35, ⟨Propositions 2.1, 2.2⟩

[52] 📖 [116], pages 28–31, ⟨Proposition 2.15⟩, 📖 [50], pages 35–37, ⟨Proposition 2.3⟩

[53] The property of *translation invariance* is of particular significance in the theory of *normed linear spaces* (a Hilbert space is a complete normed linear space equipped with an inner product).

[54] 📖 [65], page 7

The equation

$$\phi(x) = \sum_{n\in\mathbb{Z}} h_n \mathbf{D}\mathbf{T}^n \phi(x) \qquad \forall x\in\mathbb{R}$$

is called the **dilation equation**. It is also called the **refinement equation**, **two-scale difference equation**, and **two-scale relation**.

**Theorem 3.4** (dilation equation) [55] *Let an* MRA SPACE *and* SCALING FUNCTION *be as defined in Definition 3.1 page 21.*

$$\left\{\begin{array}{l}\left(\boldsymbol{L}^2_{\mathbb{R}}, \left(\!\left(\boldsymbol{V}_j\right)\!\right)\right) \text{ is an MRA SPACE}\\ \text{with SCALING FUNCTION } \phi\end{array}\right\} \implies \underbrace{\left\{\begin{array}{l}\exists \left(\!\left(h_n\right)\!\right)_{n\in\mathbb{Z}} \text{ such that}\\ \phi(x) = \sum_{n\in\mathbb{Z}} h_n \mathbf{D}\mathbf{T}^n\phi(x) \qquad \forall x\in\mathbb{R}\end{array}\right\}}_{\text{DILATION EQUATION IN “TIME”}}$$

**Lemma 3.5** [56] *Let $\phi(x)$ be a function in $\boldsymbol{L}^2_{\mathbb{R}}$ (Definition 1.2 page 3). Let $\tilde{\phi}(\omega)$ be the* FOURIER TRANSFORM *of $\phi(x)$. Let $\breve{\mathrm{h}}(\omega)$ be the* DISCRETE TIME FOURIER TRANSFORM *of a sequence $\left(\!\left(h_n\right)\!\right)_{n\in\mathbb{Z}}$.*

$$(A) \quad \phi(x) = \sum_{n\in\mathbb{Z}} h_n \mathbf{D}\mathbf{T}^n\phi(x) \quad \forall x\in\mathbb{R} \iff \tilde{\phi}(\omega) = \frac{\sqrt{2}}{2}\breve{\mathrm{h}}\left(\frac{\omega}{2}\right)\tilde{\phi}\left(\frac{\omega}{2}\right) \qquad \forall\omega\in\mathbb{R} \qquad (1)$$

$$\iff \tilde{\phi}(\omega) = \tilde{\phi}\left(\frac{\omega}{2^N}\right)\prod_{n=1}^{N}\frac{\sqrt{2}}{2}\breve{\mathrm{h}}\left(\frac{\omega}{2^n}\right) \qquad \forall n\in\mathbb{N}, \omega\in\mathbb{R} \qquad (2)$$

Definition 3.6 (next) formally defines the coefficients that appear in Theorem 3.4 (page 23).

**Definition 3.6** Let $\left(\boldsymbol{L}^2_{\mathbb{R}}, \left(\!\left(\boldsymbol{V}_j\right)\!\right)\right)$ be a multiresolution analysis space with scaling function $\phi$. Let $\left(\!\left(h_n\right)\!\right)_{n\in\mathbb{Z}}$ be a sequence of coefficients such that $\phi = \sum_{n\in\mathbb{Z}} h_n\mathbf{D}\mathbf{T}^n\phi$. A **multiresolution system** is the tuple $\left(\boldsymbol{L}^2_{\mathbb{R}}, \left(\!\left(\boldsymbol{V}_j\right)\!\right), \phi, \left(\!\left(h_n\right)\!\right)\right)$. The sequence $\left(\!\left(h_n\right)\!\right)_{n\in\mathbb{Z}}$ is the **scaling coefficient sequence**. A multiresolution system is also called an **MRA system**. An *MRA system* is an **orthonormal MRA system** if $\left\{\mathbf{T}^n\phi|n\in\mathbb{Z}\right\}$ is *orthonormal*.

**Definition 3.7** Let $\left(\boldsymbol{L}^2_{\mathbb{R}}, \left(\!\left(\boldsymbol{V}_j\right)\!\right), \phi, \left(\!\left(h_n\right)\!\right)\right)$ be a multiresolution system, and $\mathbf{D}$ the dilation operator. The **normalization coefficient at resolution** $n$ is the quantity

$$\left\|\mathbf{D}^j\phi\right\|.$$

**Theorem 3.8** [57] *Let $\left(\boldsymbol{L}^2_{\mathbb{R}}, \left(\!\left(\boldsymbol{V}_j\right)\!\right), \phi, \left(\!\left(h_n\right)\!\right)\right)$ be an* MRA SYSTEM *(Definition 3.6 page 23). Let $\operatorname{span} A$ be the* LINEAR SPAN *of a set $A$.*

$$\underbrace{\operatorname{span}\left\{\mathbf{T}^n\phi|n\in\mathbb{Z}\right\} = \boldsymbol{V}_0}_{\left\{\mathbf{T}^n\phi|n\in\mathbb{Z}\right\} \text{ is a BASIS for } \boldsymbol{V}_0} \implies \underbrace{\operatorname{span}\left\{\mathbf{D}^j\mathbf{T}^n\phi\big|n\in\mathbb{Z}\right\} = \boldsymbol{V}_j \quad \forall j\in\mathbb{W}}_{\left\{\mathbf{D}^j\mathbf{T}^n\phi|n\in\mathbb{Z}\right\} \text{ is a BASIS for } \boldsymbol{V}_j}$$

---

[55] [50], page 39, ⟨Theorem 2.1⟩
[56] [87], page 228, [50], pages 39–41, ⟨Lemma 2.1⟩
[57] [50], page 43, ⟨Theorem 2.2⟩

### 3.2.3 Necessary Conditions

**Theorem 3.9** (admissibility condition) [58] *Let* $\hat{\mathrm{h}}(z)$ *be the* Z-TRANSFORM *(Definition 1.39 page 11) and* $\check{\mathrm{h}}(\omega)$ *the* DISCRETE-TIME FOURIER TRANSFORM *(Definition 1.42 page 12) of a sequence* $(\!(h_n)\!)_{n\in\mathbb{Z}}$. $\left\{\left(\boldsymbol{L}^2_{\mathbb{R}}, (\!(\boldsymbol{V}_j)\!), \phi, (\!(h_n)\!)\right)\right.$ *is an* MRA SYSTEM *(Definition 3.6 page 23)*$\left.\right\}$

$$\underset{\not\Leftarrow}{\Longrightarrow} \underbrace{\left\{\sum_{n\in\mathbb{Z}} h_n = \sqrt{2}\right\}}_{\text{(1) ADMISSIBILITY in "time"}} \iff \underbrace{\left\{\hat{\mathrm{h}}(z)\Big|_{z=1} = \sqrt{2}\right\}}_{\text{(2) ADMISSIBILITY in "z domain"}} \iff \underbrace{\left\{\check{\mathrm{h}}(\omega)\Big|_{\omega=0} = \sqrt{2}\right\}}_{\text{(3) ADMISSIBILITY in "frequency"}}$$

**Counterexample 3.10** Let $\left(\boldsymbol{L}^2_{\mathbb{R}}, (\!(\boldsymbol{V}_j)\!), \phi, (\!(h_n)\!)\right)$ be an *MRA system* (Definition 3.6 page 23).

$$\left\{ (\!(h_n)\!) \triangleq \sqrt{2}\bar{\delta}_{n-1} \triangleq \begin{cases} \sqrt{2} & \text{for } n=1 \\ 0 & \text{otherwise.} \end{cases} \right\} \Longrightarrow \{\phi(x) = 0\}$$

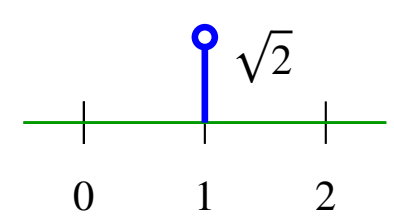

which means

$$\left\{\sum_{n\in\mathbb{Z}} h_n = \sqrt{2}\right\} \quad\not\Longrightarrow\quad \left\{\left(\boldsymbol{L}^2_{\mathbb{R}}, (\!(\boldsymbol{V}_j)\!), \phi, (\!(h_n)\!)\right) \text{ is an MRA system for } \boldsymbol{L}^2_{\mathbb{R}}.\right\}$$

✎PROOF:

$$\begin{aligned}
\phi(x) &= \sum_{n\in\mathbb{Z}} h_n \mathbf{D}\mathbf{T}^n \phi(x) && \text{by } \textit{dilation equation} \text{ (Theorem 3.4 page 23)} \\
&= \sum_{n\in\mathbb{Z}} h_n \phi(2x-n) && \text{by definitions of } \mathbf{D} \text{ and } \mathbf{T} \text{ (Definition 2.1 page 14)} \\
&= \sum_{n\in\mathbb{Z}} \underbrace{\sqrt{2}\bar{\delta}_{n-1}}_{(\!(h_n)\!)} \phi(2x-n) && \text{by definitions of } (\!(h_n)\!) \\
&= \sqrt{2}\phi(2x-1) && \text{by definition of } \phi(x) \\
\Longrightarrow \phi(x) &= 0
\end{aligned}$$

This implies $\phi(x) = 0$, which implies that $\left(\boldsymbol{L}^2_{\mathbb{R}}, (\!(\boldsymbol{V}_j)\!), \phi, (\!(h_n)\!)\right)$ is *not* an *MRA system* for $\boldsymbol{L}^2_{\mathbb{R}}$ because

$$\left(\bigcup_{j\in\mathbb{Z}} \boldsymbol{V}_j\right)^{-} = \left(\bigcup_{j\in\mathbb{Z}} \operatorname{span}\left\{\mathbf{D}^j\mathbf{T}^n\phi\,|_{n\in\mathbb{Z}}\right\}\right)^{-} \neq \boldsymbol{L}^2_{\mathbb{R}}$$

(the *least upper bound* is *not* $\boldsymbol{L}^2_{\mathbb{R}}$). ✎

[58] 📚 [50], pages 36–37, ⟨Theorem 2.3⟩

**Theorem 3.11** (Quadrature condition in "time") [59] *Let* $\left(L^2_{\mathbb{R}}, (\!(V_j)\!), \phi, (\!(h_n)\!)\right)$ *be an* MRA SYSTEM *(Definition 3.6 page 23)*.

$$\sum_{m\in\mathbb{Z}} h_m \sum_{k\in\mathbb{Z}} h_k^* \left\langle \phi \mid \mathbf{T}^{2n-m+k}\phi \right\rangle = \langle \phi \mid \mathbf{T}^n \phi \rangle \qquad \forall n\in\mathbb{Z}$$

### 3.2.4 Sufficient conditions

Theorem 3.12 (next) gives a set of *sufficient* conditions on the *scaling function* (Definition 3.1 page 21) $\phi$ to generate an *MRA*.

**Theorem 3.12** [60] *Let an* MRA *be defined as in Definition 3.1 page 21. Let* $V_j \triangleq \operatorname{span}\left\{ \mathbf{T}\phi(x) | n\in\mathbb{Z} \right\}$.

$$\left\{\begin{array}{lll} (1). & (\!(\mathbf{T}^n\phi)\!) \textit{ is a } \text{RIESZ SEQUENCE} & \textit{and} \\ (2). & \exists\, (\!(h_n)\!) \text{ such that } \phi(x) = \sum_{n\in\mathbb{Z}} h_n \mathbf{D}\mathbf{T}^n\phi(x) & \textit{and} \\ (3). & \tilde{\phi}(\omega) \textit{ is } \text{CONTINUOUS} \textit{ at } 0 & \textit{and} \\ (4). & \tilde{\phi}(0) \neq 0 & \end{array}\right\} \implies \left\{ (\!(V_j)\!)_{j\in\mathbb{Z}} \textit{ is an } \text{MRA} \right\}$$

## 3.3 Wavelet analysis

### 3.3.1 Definition

The term "wavelet" comes from the French word "*ondelette*", meaning "small wave". And in essence, wavelets are "small waves" (as opposed to the "long waves" of Fourier analysis) that form a basis for the Hilbert space $L^2_{\mathbb{R}}$.[61]

**Definition 3.13** [62] Let $\mathbf{T}$ and $\mathbf{D}$ be as defined in Definition 2.1 page 14. A function $\psi(x)$ in $L^2_{\mathbb{R}}$ is a **wavelet function** for $L^2_{\mathbb{R}}$ if

$\left\{ \mathbf{D}^j\mathbf{T}^n\psi \,|\, j,n\in\mathbb{Z} \right\}$ is a *Riesz basis* for $L^2_{\mathbb{R}}$.

In this case, $\psi$ is also called the **mother wavelet** of the basis $\left\{ \mathbf{D}^j\mathbf{T}^n\psi \,|\, j,n\in\mathbb{Z} \right\}$. The sequence of subspaces $(\!(W_j)\!)_{j\in\mathbb{Z}}$ is the **wavelet analysis** induced by $\psi$, where each subspace $W_j$ is defined as

$$W_j \triangleq \operatorname{span}\left\{ \mathbf{D}^j\mathbf{T}^n\psi \middle| n\in\mathbb{Z} \right\}.$$

---

[59] [50], page 48, ⟨Theorem 2.4⟩

[60] [116], page 28, ⟨Theorem 2.13⟩, [93], page 313, ⟨Theorem 6.4.27⟩, [50], pages 49–50, ⟨Theorem 2.6⟩

[61] [105], page ix, [8], page 191

[62] [116], page 17, ⟨Definition 2.1⟩, [50], page 50, ⟨Definition 2.4⟩

A *wavelet analysis* $(\!( \boldsymbol{W}_j )\!)$ is often constructed from a *multiresolution anaysis* (Definition 3.1 page 21) $(\!( \boldsymbol{V}_j )\!)$ under the relationship

$$\boldsymbol{V}_{j+1} = \boldsymbol{V}_j \dotplus \boldsymbol{W}_j , \qquad \text{where } \dotplus \text{ is subspace addition (\textit{Minkowski addition}).}$$

By this relationship alone, $(\!( \boldsymbol{W}_j )\!)$ is in no way uniquely defined in terms of a multiresolution analysis $(\!( \boldsymbol{V}_j )\!)$. In general there are many possible complements of a subspace $\boldsymbol{V}_j$. To uniquely define such a wavelet subspace, one or more additional constraints are required. One of the most common additional constraints is *orthogonality*, such that $\boldsymbol{V}_j$ and $\boldsymbol{W}_j$ are orthogonal to each other.

### 3.3.2 Dilation equation

Suppose $(\!(\mathbf{T}^n\psi)\!)_{n\in\mathbb{Z}}$ is a basis for $\boldsymbol{W}_0$. By Definition 3.13 page 25, the wavelet subspace $\boldsymbol{W}_0$ is contained in the scaling subspace $\boldsymbol{V}_1$. By Definition 3.1 page 21, the sequence $(\!(\mathbf{DT}^n\phi)\!)_{n\in\mathbb{Z}}$ is a basis for $\boldsymbol{V}_1$. Because $\boldsymbol{W}_0$ is contained in $\boldsymbol{V}_1$, the sequence $(\!(\mathbf{DT}^n\phi)\!)_{n\in\mathbb{Z}}$ is also a basis for $\boldsymbol{W}_0$.

**Theorem 3.14** [63] *Let* $\left(\boldsymbol{L}^2_{\mathbb{R}}, (\!( \boldsymbol{V}_j )\!), \phi, (\!(h_n)\!)\right)$ *be a multiresolution system and* $(\!( \boldsymbol{W}_j )\!)_{j\in\mathbb{Z}}$ *a wavelet analysis with respect to* $\left(\boldsymbol{L}^2_{\mathbb{R}}, (\!( \boldsymbol{V}_j )\!), \phi, (\!(h_n)\!)\right)$ *and with wavelet function* $\psi$.

$$\exists\, (\!(g_n)\!)_{n\in\mathbb{Z}} \quad \text{such that} \quad \psi = \sum_{n\in\mathbb{Z}} g_n \mathbf{DT}^n\phi$$

A *wavelet system* (next definition) consists of two subspace sequences:

- A **multiresolution analysis** $(\!( \boldsymbol{V}_j )\!)$ (Definition 3.1 page 21) provides "coarse" approximations of a function in $\boldsymbol{L}^2_{\mathbb{R}}$ at different "scales" or resolutions.
- A **wavelet analysis** $(\!( \boldsymbol{W}_j )\!)$ provides the "detail" of the function missing from the approximation provided by a given scaling subspace (Definition 3.13 page 25).

**Definition 3.15** Let $\left(\boldsymbol{L}^2_{\mathbb{R}}, (\!( \boldsymbol{V}_j )\!), \phi, (\!(h_n)\!)\right)$ be a multiresolution system (Definition 3.1 page 21) and $(\!( \boldsymbol{W}_j )\!)_{j\in\mathbb{Z}}$ a wavelet analysis (Definition 3.13 page 25) with respect to $(\!( \boldsymbol{V}_j )\!)_{j\in\mathbb{Z}}$. Let $(\!(g_n)\!)_{n\in\mathbb{Z}}$ be a sequence of coefficients such that $\psi = \sum_{n\in\mathbb{Z}} g_n \mathbf{DT}^n\phi$. A **wavelet system** is the tuple
$\left(\boldsymbol{L}^2_{\mathbb{R}}, (\!( \boldsymbol{V}_j )\!), (\!( \boldsymbol{W}_j )\!), \phi, \psi, (\!(h_n)\!), (\!(g_n)\!)\right)$
and the sequence $(\!(g_n)\!)_{n\in\mathbb{Z}}$ is the **wavelet coefficient sequence**.

[63] [50], page 51, ⟨Theorem 2.6⟩

### 3.3.3 Necessary conditions

**Theorem 3.16** (quadrature conditions in "time") [64] *Let* $\left(L^2_{\mathbb{R}}, (\!(V_j)\!), (\!(W_j)\!), \phi, \psi, (\!(h_n)\!), (\!(g_n)\!)\right)$ *be a wavelet system (Definition 3.15 page 26).*

1. $\sum_{m\in\mathbb{Z}} h_m \sum_{k\in\mathbb{Z}} h_k^* \left\langle \phi \,|\, \mathbf{T}^{2n-m+k}\phi \right\rangle = \langle \phi \,|\, \mathbf{T}^n\phi \rangle \quad \forall n\in\mathbb{Z}$
2. $\sum_{m\in\mathbb{Z}} g_m \sum_{k\in\mathbb{Z}} g_k^* \left\langle \phi \,|\, \mathbf{T}^{2n-m+k}\phi \right\rangle = \langle \psi \,|\, \mathbf{T}^n\psi \rangle \quad \forall n\in\mathbb{Z}$
3. $\sum_{m\in\mathbb{Z}} h_m \sum_{k\in\mathbb{Z}} g_k^* \left\langle \phi \,|\, \mathbf{T}^{2n-m+k}\phi \right\rangle = \langle \phi \,|\, \mathbf{T}^n\psi \rangle \quad \forall n\in\mathbb{Z}$

**Proposition 3.17** [65] *Let* $\left(L^2_{\mathbb{R}}, (\!(V_j)\!), (\!(W_j)\!), \phi, \psi, (\!(h_n)\!), (\!(g_n)\!)\right)$ *be a wavelet system. Let* $\tilde{\phi}(\omega)$ *and* $\tilde{\psi}(\omega)$ *be the* FOURIER TRANSFORM*s of* $\phi(x)$ *and* $\psi(x)$*, respectively. Let* $\breve{g}(\omega)$ *be the* DISCRETE TIME FOURIER TRANSFORM *of* $(\!(g_n)\!)$.

$$\tilde{\psi}(\omega) = \frac{\sqrt{2}}{2}\, \breve{g}\left(\frac{\omega}{2}\right) \tilde{\phi}\left(\frac{\omega}{2}\right)$$

### 3.3.4 Sufficient condition

In this text, an often used sufficient condition for designing the *wavelet coefficient sequence* $(\!(g_n)\!)$ (Definition 3.15 page 26) is the *conjugate quadrature filter condition.* It expresses the sequence $(\!(g_n)\!)$ in terms of the *scaling coefficient sequence* (Definition 3.6 page 23) and a "shift" integer $N$ as $g_n = \pm(-1)^n h^*_{N-n}$.

**Theorem 3.18** [66] *Let* $\left(L^2_{\mathbb{R}}, (\!(V_j)\!), (\!(W_j)\!), \phi, \psi, (\!(h_n)\!), (\!(g_n)\!)\right)$ *be a* WAVELET SYSTEM *(Definition 3.15 page 26). Let* $\breve{g}(\omega)$ *be the* DTFT *(Definition 1.42 page 12) and* $\hat{g}(z)$ *the* Z-TRANSFORM *(Definition 1.39 page 11) of* $(\!(g_n)\!)$.

$$\underbrace{g_n = \pm(-1)^n h^*_{N-n},\ N\in\mathbb{Z}}_{\text{CONJUGATE QUADRATURE FILTER}} \iff \breve{g}(\omega) = \pm(-1)^N e^{-i\omega N}\breve{h}^*(\omega+\pi)\Big|_{\omega=\pi} \quad (1)$$

$$\implies \sum_{n\in\mathbb{Z}}(-1)^n g_n = \sqrt{2} \quad (2)$$

$$\iff \hat{g}(z)\Big|_{z=-1} = \sqrt{2} \quad (3)$$

$$\iff \breve{g}(\omega)\Big|_{\omega=\pi} = \sqrt{2} \quad (4)$$

[64] [50], pages 55–56, ⟨Theorem 2.9⟩
[65] [50], page 56, ⟨Proposition 2.7⟩
[66] [50], pages 58–59, ⟨Theorem 2.11⟩

### 3.4 Support size

The *support* of a function is what it's non-zero part "sits" on. If the support of the scaling coefficients $\langle\!\langle h_n \rangle\!\rangle$ goes from say $[0,3]$ in $\mathbb{Z}$, what is the support of the scaling function $\phi(x)$? The answer is $[0,3]$ in $\mathbb{R}$—essentially the same as the support of $\langle\!\langle h_n \rangle\!\rangle$ except that the two functions have different domains ($\mathbb{Z}$ versus $\mathbb{R}$). This concept is defined in Definition 3.19 (next definition) and proven in Theorem 3.20 (next theorem).

**Definition 3.19** Let $\left(\boldsymbol{L}^2_{\mathbb{R}}, \left(\!\left(\boldsymbol{V}_j\right)\!\right), \left(\!\left(\boldsymbol{W}_j\right)\!\right), \phi, \psi, \langle\!\langle h_n \rangle\!\rangle, \langle\!\langle g_n \rangle\!\rangle\right)$ be a wavelet system. Let $X^-$ represent the closure of a set $X$ in $\boldsymbol{L}^2_{\mathbb{R}}$, $\vee X$ the *least upper bound* of an ordered set $(X, \leq)$, $\wedge X$ the *greatest lower bound* of an ordered set $(X, \leq)$, and

$$\begin{array}{rcll} \lfloor x \rfloor & \triangleq & \bigvee \{n \in \mathbb{Z} \,|\, n \leq x\} \quad \forall x \in \mathbb{R} & (\textit{floor} \text{ of } x) \\ \lceil x \rceil & \triangleq & \bigwedge \{n \in \mathbb{Z} \,|\, n \geq x\} \quad \forall x \in \mathbb{R} & (\textit{ceiling} \text{ of } x). \end{array}$$

The **support** $\mathsf{S}\mathrm{f}$ of a function $\mathrm{f} \in Y^X$ is defined as

$$\mathsf{S}\mathrm{f} \triangleq \begin{cases} \{x \in \mathbb{R} \,|\, \mathrm{f}(x) \neq 0\}^- & \text{for } X = \mathbb{R} \quad (\text{domain of f is } \mathbb{R}) \\ \left\{x \in \mathbb{R} \,\middle|\, \mathrm{f}(\lfloor x \rfloor) \neq 0 \text{ and } \mathrm{f}(\lceil x \rceil) \neq 0\right\}^- & \text{for } X = \mathbb{Z} \quad (\text{domain of f is } \mathbb{Z}) \end{cases}$$

**Theorem 3.20** (support size) [67] *Let* $\left(\boldsymbol{L}^2_{\mathbb{R}}, \left(\!\left(\boldsymbol{V}_j\right)\!\right), \left(\!\left(\boldsymbol{W}_j\right)\!\right), \phi, \psi, \langle\!\langle h_n \rangle\!\rangle, \langle\!\langle g_n \rangle\!\rangle\right)$ *be a wavelet system. Let* $\operatorname{supp} \mathrm{f}$ *be the support of a function* f (Definition 3.19 page 28).

$$\operatorname{supp} \phi = \operatorname{supp} \mathrm{h}$$

# 4 Background: binomial relations

### 4.1 Factorials

**Definition 4.1** (factorial) The **factorial** $n!$ is defined as

$$n! \triangleq \begin{cases} n(n-1)(n-2)\cdots 1 & \text{for } n \in \mathbb{Z},\ n \geq 1 \\ 1 & \text{for } n \in \mathbb{Z},\ n = 0 \\ 0 & \text{for } n \in \mathbb{Z},\ n \leq -1 \end{cases}$$

**Definition 4.2** [68] The quantities "$x$ **to the** $m$ **falling**", "$x$ **to the** $m$ **rising**", "$x$ **to the** $m$ **central**" are defined as follows:

[67] [87], pages 243–244, [50], pages 60–61, ⟨Theorem 2.12⟩

[68] [47], pages 47–48, ⟨equations (2.43), (2.44)⟩, [75], page 414, ⟨(2.11), (2.12)⟩, [3], page 10, [102], page 8, ⟨descending, ascending, and central factorials⟩, [101], page 8, ⟨descending, ascending, and central factorials⟩

$$x^{\underline{m}} \triangleq \begin{cases} \underbrace{x(x-1)(x-2)\cdots(x-m+1)}_{m \text{ factors}} & \forall x\in\mathbb{C}, m\in\mathbb{N} \\ 1 & \forall x\in\mathbb{C}, m=0 \end{cases} \qquad \text{("}x\textbf{ to the }m\textbf{ falling}\text{")}$$

$$x^{\overline{m}} \triangleq \begin{cases} \underbrace{x(x+1)(x+2)\cdots(x+m-1)}_{m \text{ factors}} & \forall x\in\mathbb{C}, m\in\mathbb{N} \\ 1 & \forall x\in\mathbb{C}, m=0 \end{cases} \qquad \text{("}x\textbf{ to the }m\textbf{ rising}\text{")}$$

$$x^{\overline{\underline{m}}} \triangleq \begin{cases} \underbrace{x\left(x+\frac{m}{2}-1\right)\left(x+\frac{m}{2}-2\right)\cdots\left(x-\frac{m}{2}+1\right)}_{m \text{ factors}} & \forall x\in\mathbb{C}, m\in\mathbb{N} \\ 1 & \forall x\in\mathbb{C}, m=0 \end{cases} \qquad \text{("}x\textbf{ to the }m\textbf{ central}\text{")}$$

The rising and central expressions may be expressed in terms of the falling expression (next).

**Proposition 4.3** [69]

$$x^{\overline{m}} = (-1)^m x^{\underline{m}} \qquad\qquad x^{\overline{\underline{m}}} = x\left(x+\frac{m}{2}-1\right)^{\underline{(m-1)}}$$

✎PROOF:

$$\begin{aligned}
(-1)^m(-x)^{\underline{m}} &= (-1)^m[(-x)(-x-1)(-x-2)\cdots(-x-m+1)] && \text{by Definition 4.2 page 28}\\
&= (-1)^m(-1)^m[(x)(x+1)(x+2)\cdots(x+m-1)]\\
&= x^{\overline{m}} && \text{by Definition 4.2 page 28}
\end{aligned}$$

$$\begin{aligned}
x\left(x+\frac{m}{2}-1\right)^{\underline{(m-1)}} &= x\left(x+\frac{m}{2}-1\right)\left(x+\frac{m}{2}-1-1\right)\cdots\left(x+\frac{m}{2}-1-(m-1)+1\right) && \text{by Definition 4.2 page 28}\\
&= x\left(x+\frac{m}{2}-1\right)\left(x+\frac{m}{2}-2\right)\cdots\left(x-\frac{m}{2}+1\right)\\
&= x^{\overline{\underline{m}}}
\end{aligned}$$

✎

## 4.2 Binomial identities

**Definition 4.4** (Binomial coefficient) [70] Let $\mathbb{C}$ be the set of complex numbers and $\mathbb{Z}$ the set of integers. Let $x^{\underline{m}}$ represent "$x$ to the $m$ falling" (Definition 4.2). Let $n!$ represent "$n$ factorial" (Definition 4.1). The **binomial coefficient** $\binom{x}{k}$ is defined as

[69] [102], page 8, ⟨(3)⟩
[70] [47], page 154, ⟨equation (5.1)⟩, [3], page 10, ⟨(1)⟩, [26], pages 149–150, [103]

$$\binom{x}{k} \triangleq \begin{cases} \frac{x^{\underline{k}}}{k!} & \forall x \in \mathbb{C} \quad k \in \mathbb{W} \quad (k = 0, 1, 2, 3, \ldots) \\ 0 & \forall x \in \mathbb{C} \quad k \in \mathbb{Z}^- \quad (k = -1, -2, -3, \ldots) \end{cases}$$

The value $x$ is called the **upper index** and the value $k$ is called the **lower index**.

**Proposition 4.5** *Let* $\binom{n}{k}$ *be the* BINOMIAL COEFFICIENT *(Definition 4.4 page 29).*

1. $\binom{x}{0} = 1 \quad \forall x \in \mathbb{C}$
2. $\binom{n}{n} = 1 \quad \forall n \in \mathbb{W}$
3. $\binom{x}{1} = x \quad \forall x \in \mathbb{C}$
4. $\binom{x}{k} = 0 \quad \forall x \in \mathbb{C},\ x < k$

✎PROOF:

(1) Proof that $\binom{x}{0} = 1$:

$$\begin{aligned} \binom{x}{0} &= \frac{x^{\underline{0}}}{0!} && \text{by Definition 4.4 page 29} \\ &= \frac{x^{\underline{0}}}{1} && \text{by Definition 4.1 page 28} \\ &= 1 && \text{by Definition 4.2 page 28} \end{aligned}$$

(2) Proof that $\binom{n}{n} = 1$:

$$\begin{aligned} \binom{n}{n} &= \frac{n^{\underline{n}}}{n!} && \text{by Definition 4.4 page 29} \\ &= \frac{n(n-1)\cdots(n-n+1)}{n!} && \text{by Definition 4.2 page 28} \\ &= \frac{n(n-1)\cdots(1)}{n(n-1)\cdots(1)} && \text{by Definition 4.1 page 28} \\ &= 1 \end{aligned}$$

(3) Proof that $\binom{x}{1} = x$:

$$\begin{aligned} \binom{x}{1} &= \frac{x^{\underline{1}}}{1!} && \text{by Definition 4.4 page 29} \\ &= \frac{x^{\underline{1}}}{1} && \text{by Definition 4.1 page 28} \\ &= x && \text{by Definition 4.2 page 28} \end{aligned}$$

(4) Proof that $\binom{x}{k} = 0, \forall x < k$:

$$\begin{aligned} \binom{x}{k} &= \frac{x^{\underline{k}}}{k!} && \text{by Definition 4.4 page 29} \\ &= \frac{x(x-1)\cdots(0)\cdots(x-k+1)}{k!} && \text{by Definition 4.2 page 28} \\ &= 0 \end{aligned}$$

✏

**Theorem 4.6** [71] *Let* $\binom{n}{k}$ *be the* BINOMIAL COEFFICIENT *(Definition 4.4 page 29).*

| | | | |
|---|---|---|---|
| 1. | $\binom{n}{k} = \frac{n!}{k!\,(n-k)!}$ | $\forall n,k\in\mathbb{Z},\ n\geq k\geq 0$ | (FACTORIAL EXPANSION) |
| 2. | $\binom{n}{k} = \binom{n}{n-k}$ | $\forall n,k\in\mathbb{Z},\ n\geq 0$ | (SYMMETRY) |
| 3. | $\binom{n+x+1}{n} = \binom{n+x}{n} + \binom{n+x}{n-1}$ | $\forall n\in\mathbb{Z},\ x\in\mathbb{C}$ | (PASCAL'S RULE) |
| 4. | $\binom{x+1}{k+1} = \binom{x}{k+1} + \binom{x}{k}$ | $\forall k\in\mathbb{Z},\ x\in\mathbb{C}$ | (PASCAL'S IDENTITY / STIFEL FORMULA) |
| 5. | $\binom{x}{m}\binom{m}{k} = \binom{x}{k}\binom{x-k}{m-k}$ | $\forall k,m\in\mathbb{Z},\ x\in\mathbb{C}$ | (TRINOMIAL REVISION) |
| 6. | $\binom{x}{k} = \frac{x}{k}\binom{x-1}{k-1}$ | $\forall k\in\mathbb{Z},\ x\in\mathbb{C}$ | (ABSORPTION IDENTITY) |
| 7. | $\binom{x}{k} = (-1)^k\binom{k-x-1}{k}$ | $\forall k\in\mathbb{Z},\ x\in\mathbb{C}$ | (UPPER NEGATION) |
| 8. | $\binom{x}{k} = \binom{x-2}{k-2} + 2\binom{x-2}{k-1} + \binom{x-2}{k}$ | $\forall k\in\mathbb{Z},\ x\in\mathbb{C}$ | (SECOND-ORDER PASCAL'S IDENTITY) |
| 9. | $\binom{x-1}{k-1}\binom{x}{k+1}\binom{x+1}{k} = \binom{x-1}{k}\binom{x}{k-1}\binom{x+1}{k+1}$ | $\forall k\in\mathbb{Z},\ x\in\mathbb{C}$ | (HEXAGON IDENTITY) |

✎PROOF:

(1) Proof for *factorial expansion*:

$$\begin{aligned}\binom{n}{k} &\triangleq \frac{n^{\underline{k}}}{k!} && \forall n,k\in\mathbb{Z},\ n\geq k\geq 0 \quad \text{by Definition 4.4}\\ &= \frac{n(n-1)(n-2)\cdot(n-k+1)}{k!} && \forall n,k\in\mathbb{Z},\ n\geq k\geq 0 \quad \text{by Definition 4.2}\end{aligned}$$

[71] 📚 [47], page 174, ⟨Table 174⟩, 📚 [43], page 221, 📚 [52], page 227, ⟨Table 4.1.2⟩, 📚 [26], pages 149–150, 📚 [103], 📚 [12], page 43, ⟨Pascal's Rule⟩, 📚 [56], page 143, ⟨*hexagon identity*, (2.15)⟩, 📚 [36], page 216, ⟨*second-order pascal identity*⟩

$$= \frac{n(n-1)(n-2)\cdot(n-k+1)(n-k)(n-k-1)\cdots 1}{k!\,(n-k)!} \qquad \forall n,k\in\mathbb{Z},\ n\ge k\ge 0 \quad \text{by Definition 4.2}$$

$$= \frac{n!}{k!\,(n-k)!} \qquad \forall n,k\in\mathbb{Z},\ n\ge k\ge 0 \quad \text{by Definition 4.1}$$

(2) Proof for *symmetry* property:

(a) Proof for $n,k\in\mathbb{Z},\ n\ge k\ge 0$: (use item 1 page 31)

$$\begin{aligned}
\binom{n}{n-k} &= \frac{n!}{(n-k)!\,(n-(n-k)!} && \forall n,k\in\mathbb{Z},\ n\ge k\ge 0 && \text{by item 1 page 31}\\
&= \frac{n!}{k!\,(n-k)!} && \forall n,k\in\mathbb{Z},\ n\ge k\ge 0 &&\\
&= \binom{n}{k} && \forall n,k\in\mathbb{Z},\ n\ge k\ge 0 && \text{by item 1 page 31}
\end{aligned}$$

(b) Proof for $n,k\in\mathbb{Z},\ n\ge 0>k$:

$$\begin{aligned}
\binom{n}{n-k} &= \frac{n^{\underline{n-k}}}{(n-k)!} && \forall n,k\in\mathbb{Z},\ n\ge 0>k && \text{by Definition 4.4}\\
&= \frac{n(n-1)(n-2)\cdots 0\cdots(n-n+k+1)}{(n-k)!} && \forall n,k\in\mathbb{Z},\ n\ge 0>k && \text{by Definition 4.2}\\
&= 0 &&&&\\
&= \binom{n}{k} && \forall n,k\in\mathbb{Z},\ n\ge 0>k && \text{by Definition 4.4}
\end{aligned}$$

(c) Proof for $n,k\in\mathbb{Z},\ n\ge 0>k$:

$$\begin{aligned}
\binom{n}{k} &= \frac{n^{\underline{k}}}{k!} && \forall n,k\in\mathbb{Z},\ k>n\ge 0 && \text{by Definition 4.4 page 29}\\
&= \frac{n(n-1)(n-2)\cdots 0\cdots(n-k+1)}{(n-k)!} && \forall n,k\in\mathbb{Z},\ k>n\ge 0 && \text{by Definition 4.2 page 28}\\
&= 0 &&&&\\
&= \binom{n}{n-k} && \forall n,k\in\mathbb{Z},\ k>n\ge 0 && \text{by Definition 4.4 page 29}
\end{aligned}$$

(3) Proof for *Pascal's Rule*:

(a) Proof for $n<0,\ x\in\mathbb{C}$:

$$\begin{aligned}
\binom{n+x}{n}+\binom{n+x}{n-1} &= 0+0 && \text{by Definition 4.4 page 29}\\
&= \binom{n+x+1}{n} && \text{by Definition 4.4 page 29}
\end{aligned}$$

(b) Proof for $n = 0,\ x \in \mathbb{C}$:

$$\begin{aligned}
\binom{n+x}{n} + \binom{n+x}{n-1} &= \binom{n+x}{0} + \binom{n+x}{-1} && \text{by } n = 0 \text{ hypothesis} \\
&= 1 + 0 && \text{by Definition 4.4 page 29} \\
&= \binom{n+x+1}{0} && \text{by Definition 4.4 page 29} \\
&= \binom{n+x+1}{n} && \text{by } n = 0 \text{ hypothesis}
\end{aligned}$$

(c) Proof for $n > 0,\ x \in \mathbb{C}$:

$$\begin{aligned}
&\binom{n+x}{n} + \binom{n+x}{n-1} \\
&\triangleq \frac{n+x^{\underline{n}}}{n!} + \frac{n+x^{\underline{n-1}}}{(n-1)!} && \text{by Definition 4.4 page 29} \\
&\triangleq \frac{(n+x)(n+x-1)\cdots(n+x-n+1)}{n!} \\
&\quad + \frac{(n+x)(n+x-1)\cdots(n+x-n+1+1)}{(n-1)!} && \text{by Definition 4.2 page 28} \\
&= \frac{[(n+x)(n+x-1)\cdots(x+1)] + [(n+x)(n+x-1)\cdots(x+2)n]}{n!} \\
&= \frac{[(x+1)+n][(n+x)(n+x-1)\cdots(x+2)]}{n!} \\
&= \frac{(n+x+1)(n+x)(n+x-1)\cdots(x+2)}{n!} \\
&\triangleq \frac{(n+x+1)^{\underline{n}}}{n!} && \text{by Definition 4.2 page 28} \\
&\triangleq \binom{n+x+1}{n} && \text{by Definition 4.4 page 29}
\end{aligned}$$

(4) Proof for *Pascal's Identity*:

$$\begin{aligned}
\binom{x+1}{k+1} &= \binom{k+y+1}{k+1} && \text{where } y \triangleq x-k \implies x = y+k \\
&= \binom{y+k}{k+1} + \binom{y+k}{k} && \text{by Pascal's Rule (item 3)} \\
&= \binom{x}{k+1} + \binom{x}{k} && \text{by definition of } m
\end{aligned}$$

(5) Proof for *Trinomial revision*:

(a) Proof for $k < 0$ case:

$$\begin{aligned}
\binom{x}{m}\binom{m}{k} &= \binom{x}{m}0 && \text{by } k < 0 \text{ hypothesis and Definition 4.4 page 29} \\
&= \cancelto{0}{\binom{x}{k}}\binom{x-k}{m-k} && \text{by } k < 0 \text{ hypothesis and Definition 4.4 page 29}
\end{aligned}$$

(b) Proof for $k \geq 0,\ m < 0$ case:

$$\begin{aligned}
\binom{x}{m}\binom{m}{k} &= 0\binom{m}{k} && \text{by } m < 0 \text{ hypothesis and Definition 4.4 page 29}\\
&= \binom{x}{k}\binom{x-k}{m-k} && \text{by } k \geq 0, m < 0 \text{ hypothesis and Definition 4.4 page 29}
\end{aligned}$$

(c) Proof for $m < k$ case:

$$\begin{aligned}
\binom{x}{m}\binom{m}{k} &= \binom{x}{m}0 && \text{by Proposition 4.5 page 30}\\
&= \binom{x}{k}\binom{x-k}{m-k} && \text{by } m < k \text{ hypothesis and Definition 4.4 page 29}
\end{aligned}$$

(d) Proof for remaining cases:

$$\begin{aligned}
&\binom{x}{m}\binom{m}{k}\\
&= \frac{x^{\underline{m}}}{m!}\frac{m^{\underline{k}}}{k!} && \text{by Definition 4.4}\\
&= \frac{x(x-1)\cdots(x-m+1)}{m!}\frac{m(m-1)\cdots(m-k+1)}{k!} && \text{by Definition 4.2}\\
&= \frac{x(x-1)\cdots(x-m+1)}{(m-k)!}\frac{1}{k!}\\
&= \frac{x(x-1)\cdots(x-k+1)}{k!}\frac{(x-k)(x-k-1)\cdots(x-m+1)}{(m-k)!}\\
&= \frac{x(x-1)\cdots(x-k+1)}{k!}\frac{(x-k)(x-k-1)\cdots((x-k)-(m-k)+1)}{(m-k)!}\\
&\triangleq \frac{x^{\underline{k}}}{k!}\frac{(x-k)^{\underline{m-k}}}{(m-k)!} && \text{by Definition 4.2}\\
&\triangleq \binom{x}{k}\binom{x-k}{m-k} && \text{by Definition 4.4}
\end{aligned}$$

(6) Proof for *Absorption identity*:

$$\begin{aligned}
\frac{x}{k}\binom{x-1}{k-1} &= \frac{1}{k}\binom{x}{1}\binom{n-1}{k-1} && \text{by Proposition 4.5 page 30}\\
&= \frac{1}{k}\binom{x}{k}\binom{k}{1} && \text{by Trinomial revision (item 5)}\\
&= \frac{1}{k}\binom{x}{k}k && \text{by Proposition 4.5 page 30}\\
&= \binom{x}{k}
\end{aligned}$$

(7) Proof for *Upper Negation*:

$$(-1)^k \binom{k-x-1}{k}$$

$$\triangleq (-1)^k \frac{(k-x-1)^{\underline{k}}}{k!}$$ by Definition 4.4 page 29

$$\triangleq (-1)^k \frac{(k-x-1)(k-x-2)(k-x-3)\cdots(k-x-1-k+1)}{k!}$$ by Definition 4.2 page 28

$$= (-1)^k \frac{(k-x-1)(k-x-2)(k-x-3)\cdots(-x)}{k!}$$

$$= (-1)^k(-1)^k \frac{(x)(x-1)\cdots(x(x-k+3)(x-k+2)(x-k+1)}{k!}$$

$$\triangleq \frac{x^{\underline{k}}}{k!}$$ by Definition 4.2 page 28

$$\triangleq \binom{x}{k}$$ by Definition 4.4 page 29

(8) Proof for *2nd Order Pascal's Identity*:

$$\binom{n-2}{k-2} + 2\binom{n-2}{k-1} + \binom{n-2}{k}$$

$$\triangleq \frac{(x-2)^{\underline{(k-2)}}}{(k-2)!} + \frac{(x-2)^{\underline{(k-1)}}}{(k-1)!} + \frac{(x-2)^{\underline{k}}}{k!}$$

$$\triangleq \frac{(x-2)(x-1)\cdots(x-k+2+1)}{(k-2)!} + 2\frac{(x-2)(x-1)\cdots(x-k+1+1)}{(k-1)!} + \frac{(x-2)\cdots(x-k+1)}{k!}$$

$$= \frac{(x-2)\cdots(x-2-k+2+1)k(k-1) + 2(x-2)\cdots(x-2-k+1+1)k + (x-2)\cdots(x-k-1)}{k!}$$

$$= \frac{(x-2)(x-1)\cdots(x-k+1)k(k-1) + 2(x-2)(x-1)\cdots(x-k)k + (n-2)(n-1)\cdots(x-k-1)}{k!}$$

$$= \frac{[(x-2)(x-1)\cdots(x-k+1)][k(k-1) + 2(x-k)k + (x-k)(x-k-1)]}{k!}$$

$$= \frac{[(x-2)(x-1)\cdots(x-k+1)][k(k-1) + 2(x-k)k - (x-k)k + (x-k)(x-1)]}{k!}$$

$$= \frac{[(x-2)(x-1)\cdots(x-k+1)][k(k-1) + (x-k)k + (x-k)(x-1)]}{k!}$$

$$= \frac{[(x-2)(x-1)\cdots(x-k+1)]\left[k^2 - k + kx - k^2 + x^2 - x - kx + k\right]}{k!}$$

$$= \frac{[(x-2)(x-1)\cdots(x-k+1)]\left[x^2 - x\right]}{k!}$$

$$= \frac{x(x-1)(x-2)(x-1)\cdots(x-k+1)}{k!}$$

$$\triangleq \frac{n^{\underline{k}}}{k!}$$

$$\triangleq \binom{n}{k}$$

(9) Proof for *Hexagon Identity*:

$$
\begin{aligned}
&\binom{x-1}{k-1}\binom{x}{k+1}\binom{x+1}{k} \\
&\triangleq \left[\frac{(x-1)^{\underline{(k-1)}}}{(k-1)!}\right]\left[\frac{x^{\underline{(k+1)}}}{(k+1)!}\right]\left[\frac{(x+1)^{\underline{k}}}{k!}\right] \\
&\triangleq \left[\frac{(x-1)\cdots(x-1-k+1+1)}{(k-1)!}\right]\left[\frac{x(x-1)\cdots(x-k-1+1)}{(k+1)!}\right]\left[\frac{(x+1)(x)\cdots(x+1-k+1)}{k!}\right] \\
&= \left[\frac{(x-1)\cdots(x-k+2)(x-k+1)}{(k-1)!}\right]\left[\frac{x(x-1)\cdots(x-k)}{(k+1)!}\right]\left[\frac{(x+1)(x)(x-1)\cdots(x-k+2)}{k!}\right] \\
&= \left[\frac{(x)(x-1)\cdots(x-k+2)}{(k-1)!}\right]\left[\frac{(x+1)x(x-1)\cdots(x-k)(x-k+1)}{(k+1)!}\right]\left[\frac{(x-1)\cdots(x-k)}{k!}\right] \\
&\triangleq \left[\frac{x^{\underline{(k-1)}}}{(k-1)!}\right]\left[\frac{(x+1)^{\underline{(k+1)}}}{(k+1)!}\right]\left[\frac{(x-1)^{\underline{k}}}{k!}\right] \\
&\triangleq \binom{x}{k-1}\binom{x+1}{k+1}\binom{x-1}{k}
\end{aligned}
$$

✏

From Pascal's Recursion we can construct *Pascal's Triangle*:[72]

$$
\begin{array}{ccccccccc}
 & & & & \binom{0}{0} & & & & \\
 & & & \binom{1}{0} & & \binom{1}{1} & & & \\
 & & \binom{2}{0} & & \binom{2}{1} & & \binom{2}{2} & & \\
 & \binom{3}{0} & & \binom{3}{1} & & \binom{3}{2} & & \binom{3}{3} & \\
\binom{4}{0} & & \binom{4}{1} & & \binom{4}{2} & & \binom{4}{3} & & \binom{4}{4} \\
 & & & & \vdots & & & &
\end{array}
=
\begin{array}{ccccccccc}
 & & & & 1 & & & & \\
 & & & 1 & & 1 & & & \\
 & & 1 & & 2 & & 1 & & \\
 & 1 & & 3 & & 3 & & 1 & \\
1 & & 4 & & 6 & & 4 & & 1 \\
 & & & & \vdots & & & &
\end{array}
$$

## 4.3 Binomial summations

**Theorem 4.7** [73] *Let* $(\!|x_n|\!)_1^N$ *and* $(\!|y_n|\!)_1^N$ *be sequences over a ring* $(\mathbb{X}, +, \times)$.

$$
\left(\sum_{n=0}^{p} x_n\right)\left(\sum_{m=0}^{q} y_m\right) = \sum_{n=0}^{p+q} \underbrace{\left(\sum_{k=\max(0,n-q)}^{\min(n,p)} x_k y_{n-k}\right)}_{\textit{Cauchy product}}
$$

✎PROOF:

[72] [91], [48], [49], [33], [55], pages 320–321, ⟨article 393⟩
[73] [7], page 237

(1)

$$\begin{aligned}
\left(\sum_{n=0}^{p} x_n\right)\left(\sum_{m=0}^{q} y_m\right) &= \sum_{n=0}^{p}\sum_{m=0}^{q} x_n y_m z^{n+m} \\
&= \sum_{n=0}^{p}\sum_{k=n}^{q+n} x_n y_{k-n} && k = n+m \qquad m = k-n \\
&\vdots \\
&= \sum_{n=0}^{p+q}\left(\sum_{k=0}^{n} x_k y_{n-k}\right)
\end{aligned}$$

(2) Perhaps the easiest way to see the relationship is by illustration with a matrix of product terms:

| | $y_0$ | $y_1$ | $y_2$ | $y_3$ | $\cdots$ | $y_q$ |
|---|---|---|---|---|---|---|
| $x_0$ | $x_0y_0$ | $x_0y_1$ | $x_0y_2$ | $x_0y_3$ | $\cdots$ | $x_0y_q$ |
| $x_1$ | $x_1y_0$ | $x_1y_1$ | $x_1y_2$ | $x_1y_3$ | $\cdots$ | $x_1y_q$ |
| $x_2$ | $x_2y_0$ | $x_2y_1$ | $x_2y_2$ | $x_2y_3$ | $\cdots$ | $x_2y_q$ |
| $x_3$ | $x_3y_0$ | $x_3y_1$ | $x_3y_2$ | $x_3y_3$ | $\cdots$ | $x_3y_q$ |
| $\vdots$ | $\vdots$ | $\vdots$ | $\vdots$ | $\vdots$ | $\ddots$ | $\vdots$ |
| $x_p$ | $x_py_0$ | $x_py_1$ | $x_py_2$ | $x_py_3$ | $\cdots$ | $x_py_q$ |

(a) The expression $\sum_{n=0}^{p}\sum_{m=0}^{q} x_n y_m z^{n+m}$ is equivalent to adding *horizontally* from left to right, from the first row to the last.

(b) If we switched the order of summation to $\sum_{m=0}^{q}\sum_{n=0}^{p} x_n y_m z^{n+m}$, then it would be equivalent to adding *vertically* from top to bottom, from the first column to the last.

(c) However the final result expression $\sum_{n=0}^{p+q}\left(\sum_{k=0}^{n} x_k y_{n-k}\right)$ is equivalent to adding *diagonally* starting from the upper left corner and proceding to the lower right.

(d) Upper limit on inner summation: Looking at the $x_k$ terms, we see that there are two constraints on $k$:

$$\left.\begin{array}{ccc} k & \le & n \\ k & \le & p \end{array}\right\} \implies k \le \min(n,p)$$

(e) Lower limit on inner summation: Looking at the $x_k$ terms, we see that there are two constraints on $k$:

$$\left.\begin{array}{ccc} k & \ge & 0 \\ k & \ge & n-q \end{array}\right\} \implies k \ge \max(0, n-q)$$

✎

**Theorem 4.8** [74] *Let* $\binom{n}{k}$ *be the* BINOMIAL COEFFICIENT *(Definition 4.4 page 29).*

[74] [47], page 169, ⟨Table 169⟩, [43], pages 218–223, [52], page 227, ⟨Table 4.1.2⟩, [56], pages 137–142, [74], [110], [118]

$$\sum_{k=0}^{n}\binom{n}{k} = 2^n \qquad \text{(row sum)}$$

$$\sum_{k=m}^{n}\binom{k}{m} = \binom{n+1}{m+1} \qquad \text{(upper sum / column sum)}$$

$$\sum_{k=0}^{n}\binom{m+k}{k} = \binom{n+m+1}{n} \qquad \text{(parallel summation formula/ southeast diagonal)}$$

$$\sum_{k=0}^{m}\binom{n-k}{m-k} = \binom{n+1}{m} \qquad \text{(northwest diagonal)}$$

$$\sum_{j=0}^{n}\binom{m}{j}\binom{n}{k-j} = \binom{m+n}{k} \qquad \text{(Vandermonde's convolution)}$$

$$\sum_{i=-j}^{n-j}\binom{m}{j+i}\binom{n}{k-i} = \binom{m+n}{j+k} \qquad \text{(alternate Vandermonde's convolution)}$$

$$\sum_{k=0}^{n}\binom{n}{k}^2 = \binom{2n}{n}$$

✎PROOF:

(1) Proof for *row sum* relation:

$$\begin{aligned}
\sum_{k=0}^{n}\binom{n}{k}x^k &= \left.\sum_{k=0}^{n}\binom{n}{k}x^k\right|_{x=1} \\
&= (1+x)^n\big|_{x=1} && \text{by Binomial Theorem} \\
&= (1+1)^n \\
&= 2^n
\end{aligned}$$

(2) Proof for *upper sum* relation (proof by induction):

(a) Proof for $(n,m) = (0,0)$ case:

$$\sum_{0}^{0}\binom{k}{m} = \binom{0}{0} = 1 = \binom{0+1}{0+1}$$

(b) Proof for $(n,m) = (1,0)$ case:

$$\sum_{0}^{1}\binom{k}{m} = \binom{1}{0} + \binom{1}{1} = 2 = \binom{1+1}{0+1}$$

(c) Proof for $(n,m) = (1,1)$ case:

$$\sum_{0}^{1}\binom{k}{m} = \binom{1}{1} = 1 = \binom{1+1}{1+1}$$

(d) Proof that $n$ case $\implies$ $n+1$ case:

$$\begin{aligned}
\sum_{k=m}^{n+1}\binom{k}{m} &= \binom{n+1}{m} + \sum_{k=m}^{n} \\
&= \binom{n+1}{m} + \binom{n+1}{m+1} && \text{by left hypothesis} \\
&= \binom{n+2}{m+1} && \text{by Pascal's recursion (Theorem 4.6 page 31)}
\end{aligned}$$

(3) Proof for *Parallel summation formula* (Proof by induction):

(a) Proof that $\sum_{k=0}^{n}\binom{m+k}{k} = \binom{n+m+1}{n}$ is true for $n=0$:

$$\begin{aligned}
\left.\sum_{k=0}^{n}\binom{m+k}{k}\right|_{n=0} &= \binom{m+0}{0} \\
&= \frac{(m+0)!}{(m-0)!\;0!} && \text{by Definition 4.4 page 29} \\
&= \frac{(m+1)!}{(m+1-0)!\;0!} \\
&= \binom{m+1}{0} && \text{by Definition 4.4 page 29} \\
&= \left.\binom{n+m+1}{n}\right|_{n=0}
\end{aligned}$$

(b) Proof that $\sum_{k=0}^{n}\binom{m+k}{k} = \binom{n+m+1}{n}$ is true for $n=1$:

$$\begin{aligned}
\left.\sum_{k=0}^{n}\binom{m+k}{k}\right|_{n=1} &= \binom{m+0}{0} + \binom{m+1}{1} \\
&= \binom{m+1}{0} + \binom{m+1}{1} \\
&= \binom{m+1+1}{1} && \text{by Pascal's Rule page 31} \\
&= \left.\binom{n+m+1}{n}\right|_{n=1}
\end{aligned}$$

(c) Proof that $\sum_{k=0}^{n}\binom{m+k}{k} = \binom{n+m+1}{n} \implies \sum_{k=0}^{n+1}\binom{m+k}{k} = \binom{(n+1)+m+1}{n+1}$:

$$\begin{aligned}
\sum_{k=0}^{n+1}\binom{m+k}{k} &= \binom{m}{0} + \sum_{k=1}^{n+1}\binom{m+k}{k} \\
&= \binom{m}{0} + \sum_{k=0}^{n}\binom{m+k+1}{k+1}
\end{aligned}$$

$$= \binom{m}{0} + \sum_{k=0}^{n} \binom{m+k}{k} - \binom{m}{0} + \binom{m+n+1}{n+1}$$

$$= \binom{n+m+1}{n} + \binom{m+n+1}{n+1} \qquad \text{by left hypothesis}$$

$$= \binom{n+m+2}{n+1} \qquad \text{by Pascal's Rule page 31}$$

$$= \binom{(n+1)+m+1}{(n+1)}$$

(4) Proof for *Vandermonde's convolution*:

$$\sum_{k=0}^{m+n} \binom{m+n}{k} x^k = (1+x)^{m+n} \qquad \text{by Binomial Theorem)}$$

$$= \left[\sum_{k=0}^{m} \binom{m}{k} x^k\right] \left[\sum_{j=0}^{n} \binom{n}{j} x^j\right] \qquad \text{by Binomial Theorem)}$$

$$= \sum_{k=0}^{m} \sum_{j=0}^{n} \binom{m}{k}\binom{n}{j} x^k x^j$$

$$= \sum_{k=0}^{m+n} \left[\sum_{j=0}^{n} \binom{m}{j}\binom{n}{k-j}\right] x^k \qquad \text{by Theorem 4.7 page 36}$$

$$\implies \binom{m+n}{k} = \sum_{j=0}^{n} \binom{m}{j}\binom{n}{k-j}$$

(5) Proof for *alternate Vandermonde's convolution*:

$$\binom{m+n}{j+k} = \binom{m+n}{u} \qquad \text{where } u \triangleq j+k \implies k = u-j$$

$$= \sum_{v=0}^{n} \binom{m}{v}\binom{n}{u-v}$$

$$= \sum_{v=0}^{n} \binom{m}{v}\binom{n}{j+k-v}$$

$$= \sum_{i+j=0}^{i+j=n} \binom{m}{j+i}\binom{n}{k-i} \qquad \text{where } i \triangleq v-j \implies v = i+j$$

$$= \sum_{i=-j}^{i=n-j} \binom{m}{j+i}\binom{n}{k-i}$$

(6) Proof that $\sum_{k=0}^{n} \binom{n}{k}^2 = \binom{2n}{n}$:

$$\binom{2n}{n} = \binom{n+n}{n}$$

$$= \sum_{k=0} n \binom{n}{k} \binom{n}{n-k}$$ by Vandermonde's convolution (item 4 page 40)

$$= \sum_{k=0} n \binom{n}{k} \binom{n}{k}$$ by item 2

$$= \sum_{k=0}^{n} \binom{n}{k}^2$$

✏

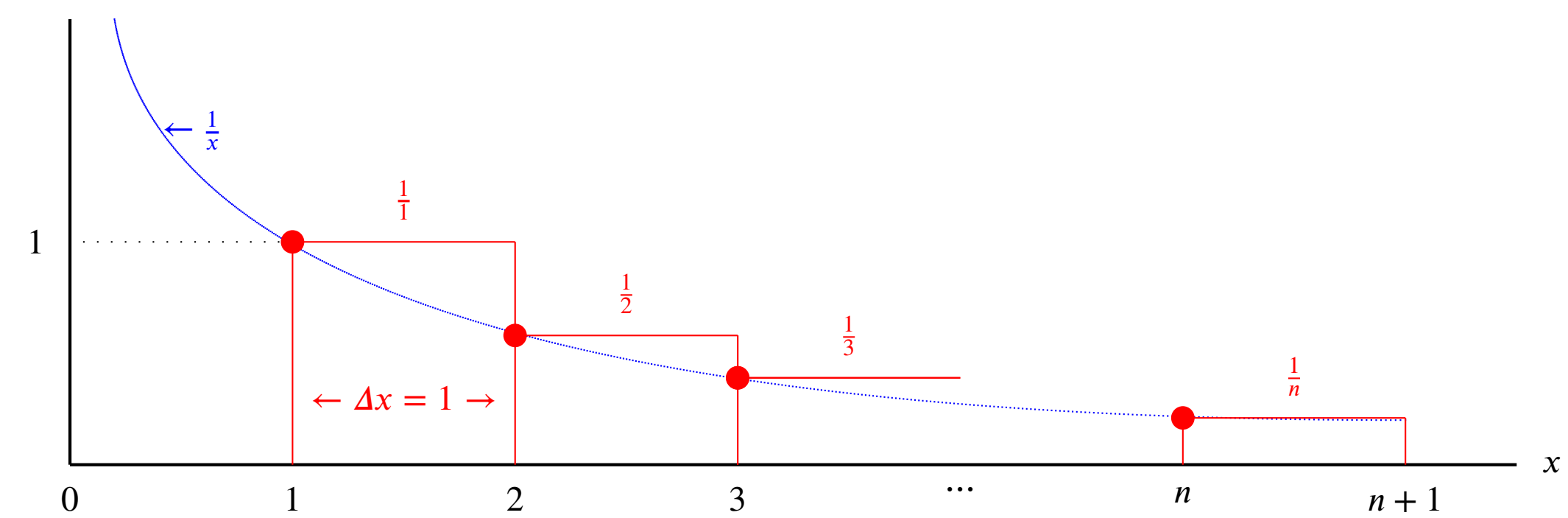

Figure 2: $\ln(n+1)$

**Theorem 4.9** [75]

$$\sum_{k=1}^{n} \frac{1}{k+1} < \ln(n+1) < \sum_{k=1}^{n} \frac{1}{k}$$

✎PROOF: The summations are simply lower and upper bounds of the integral of $\frac{1}{x}$ in the range $[1, n+1]$. This is illustrated in Figure 2.

1. Proof that $\ln(n+1) < \sum_{k=1}^{n} \frac{1}{k}$:

$$\begin{aligned}\sum_{k=1}^{n} \frac{1}{k} &> \int_1^{n+1} \frac{1}{x}\,dx \\ &= \ln x \Big|_1^{n+1} \\ &= \ln(n+1) - \cancelto{0}{\ln(1)} \\ &= \ln(n+1)\end{aligned}$$

2. Proof that $\sum_{k=1}^{n} \frac{1}{k+1} < \ln(n+1)$:

$$\begin{aligned}\sum_{k=1}^{n} \frac{1}{k+1} &< \int_1^{n+1} \frac{1}{x}\,dx \\ &= \ln(n+1) - \cancelto{0}{\ln(1)} \\ &= \ln(n+1)\end{aligned}$$

75 📘 [96], page 60

✏

# 5 B-splines

## 5.1 Definition

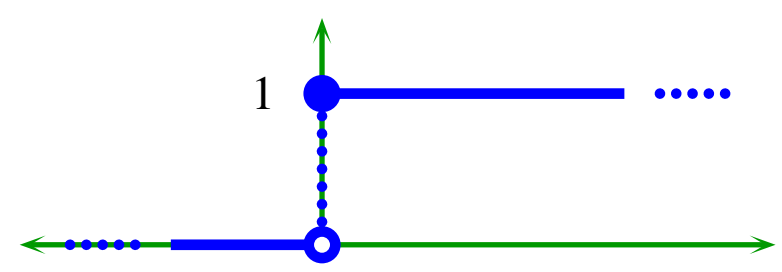

**Definition 5.1** Let $X$ be a set. The **step function** $\sigma \in \mathbb{R}^{\mathbb{R}}$ is defined as

$$\sigma(x) \triangleq \mathbb{1}_{[0,\infty)}(x) \qquad \forall x \in \mathbb{R}.$$

**Definition 5.2** [76] Let $\mathbb{1}$ be the *set indicator* function (Definition 1.3 page 3). Let $f(x) \star g(x)$ represent the *convolution* operation (Definition 1.34 page 10). The $n$**th order cardinal B-spline** $N_n$ for $n \in \mathbb{W}$ is defined as

$$N_n(x) \triangleq \begin{cases} \mathbb{1}_{[0,1]}(x) & \text{for } n = 0 \\ N_{n-1}(x) \star N_1(x) & \text{for } n \in \mathbb{W}\backslash 0 \end{cases}$$

**Lemma 5.3** [77]

$$N_n(x) = \int_0^1 N_{n-1}(x-\tau)\,d\tau \qquad \forall n \in \mathbb{W}\backslash 0$$

✎PROOF:

$$\begin{aligned} N_n(x) &\triangleq \int_{\mathbb{R}} N_{n-1}(x-\tau) N_1(\tau)\,d\tau && \text{by definition of } N_n \text{ (Definition 5.2 page 42)} \\ &= \int_0^1 N_{n-1}(x-\tau)\,d\tau && \text{by definition of } N_1 \text{ (Definition 5.2 page 42)} \end{aligned}$$

✏

**Lemma 5.4** [78] *Let $\mathbb{1}$ be the* SET INDICATOR *function (Definition 1.3 page 3). Let $\sigma(x)$ be the* STEP FUNCTION *(Definition 5.1 page 42).*

---

76 [25], page 85, ⟨(4.2.1)⟩, [22], page 140, [24], page 1
77 [22], page 140, [25], page 85, ⟨(4.2.1)⟩, [24], page 1
78 [22], page 148, ⟨Exercise 6.2⟩, [23], page 212, ⟨Exercise 10.2⟩, [99], page 136, ⟨Table 1⟩

$$\begin{aligned}
\mathsf{N}_0(x) &= \sigma(x) - \sigma(x-1) && \forall x\in\mathbb{R}\\
&= \begin{cases} 1 & for\, x \in [0,\ 1) \\ 0 & for\, x \in \mathbb{R}\backslash[0,\ 1) \end{cases}
\end{aligned}$$

$$\begin{aligned}
\mathsf{N}_1(x) &= x\sigma(x) - 2(x-1)\sigma(x-1) + (x-2)\sigma(x-2) && \forall x\in\mathbb{R}\\
&= \begin{cases} x & for\, x \in [0,\ 1] \\ -x+2 & for\, x \in [1,\ 2] \\ 0 & for\, x \in \mathbb{R}\backslash[0,\ 2] \end{cases}
\end{aligned}$$

$$\begin{aligned}
\mathsf{N}_2(x) &= \tfrac{1}{2}x^2\sigma(x) + \left[-\tfrac{3}{2}x^2 + 3x - \tfrac{3}{2}\right]\sigma(x-1) + \left[\tfrac{3}{2}x^2 - 6x + 6\right]\sigma(x-2)\\
&\quad + \left[-\tfrac{1}{2}x^2 + 3x - \tfrac{9}{2}\right]\sigma(x-3) && \forall x\in\mathbb{R}\\
&= \begin{cases} \tfrac{1}{2}x^2 & for\, x \in [0,\ 1] \\ -x^2 + 3x - \tfrac{3}{2} & for\, x \in [1,\ 2] \\ \tfrac{1}{2}x^2 - 3x + \tfrac{9}{2} & for\, x \in [2,\ 3] \\ 0 & for\, x \in \mathbb{R}\backslash[0,\ 3] \end{cases}
\end{aligned}$$

✎PROOF:

$$\begin{aligned}
\mathsf{N}_0(x) &= \mathbb{1}_{[0,\ 1]}(x) && \text{by definition of } \mathsf{N}_n \text{ (page 42)}\\
\mathsf{N}_1(x) &= \int_0^1 \mathsf{N}_0(x-\tau)\,\mathrm{d}\tau && \text{by Lemma 5.3 page 42}\\
&= \int_0^1 \mathbb{1}_{[0,\ 1]}(x-\tau)\,\mathrm{d}\tau && \text{by definition of } \mathsf{N}_1 \text{ (page 42)}\\
&= \int_{x-u=0}^{x-u=1} \mathbb{1}_{[0,\ 1]}(u)(-1)\,\mathrm{d}u && \text{where } u \triangleq x-\tau \implies \tau = x-u\\
&= \int_{u=x-1}^{u=t} \mathbb{1}_{[0,\ 1]}(u)\,\mathrm{d}u\\
&= u\sigma(u) - (u-1)\sigma(u-1) + a\Big|_{u=x-1}^{u=t}\\
&= \underbrace{\{x\sigma(x) - (x-1)\sigma(x-1) + a\}}_{u=t}\\
&\qquad - \underbrace{\{(x-1)\sigma(x-1) - (x-2)\sigma(x-2) + a\}}_{u=x-1}\\
&= x\sigma(x) - 2(x-1)\sigma(x-1) + (x-2)\sigma(x-2)\\
&= \begin{cases} t & \text{for } x \in [0,\ 1] \\ -x+2 & \text{for } x \in [1,\ 2] \\ 0 & \text{for } x \in \mathbb{R}\backslash[0,\ 2] \end{cases}
\end{aligned}$$

$$
\begin{aligned}
&\mathsf{N}_2(x)\\
&= \int_0^1 \mathsf{N}_1(x-\tau)\,\mathrm{d}\tau \qquad \text{by Lemma 5.3 page 42}\\
&= \int_0^1 (x-\tau)\sigma(x-\tau) - 2(x-\tau-1)\sigma(x-\tau-1) + (x-\tau-2)\sigma(x-\tau-2)\,\mathrm{d}\tau \qquad \text{by result for } \mathsf{N}_2\\
&= \int_{x-u=0}^{x-u=1} u\sigma(u) - 2(u-1)\sigma(u-1) + (u-2)\sigma(u-2)(-1)\,\mathrm{d}u \qquad \text{where } u \triangleq x-\tau \implies \tau = x-u\\
&= \int_{u=x-1}^{u=t} u\sigma(u)\,\mathrm{d}u + \int_{u=x-1}^{u=t} (-2u+2)\sigma(u-1)\,\mathrm{d}u + \int_{u=x-1}^{u=t} (u-2)\sigma(u-2)\,\mathrm{d}u\\
&= \left[\frac{1}{2}+a\right]u^2\sigma(u) + \left[-u^2+2u+b\right]\sigma(u-1) + \left[\frac{1}{2}u^2-2u+c\right]\sigma(u-2)\Big|_{u=x-1}^{u=t}\\
&= \underbrace{\left\{\left[\tfrac{1}{2}x^2+a\right]\sigma(x) + \left[-x^2+2x+b\right]\sigma(x-1) + \left[\tfrac{1}{2}x^2-2x+c\right]\sigma(x-2)\right\}}_{u=t}\\
&- \underbrace{\left\{\left[\tfrac{1}{2}(x-1)^2+a\right]\sigma(x-1) + \left[-(x-1)^2+2(x-1)+b\right]\sigma(x-2) + \left[\tfrac{1}{2}(x-1)^2-2(x-1)+c\right]\sigma(x-3)\right\}}_{u=x-1}\\
&= \left[\frac{1}{2}x^2+a\right]\sigma(x) + \left[-x^2+2x+b-\frac{1}{2}x^2+x-\frac{1}{2}-a\right]\sigma(x-1)\\
&\qquad + \left[\frac{1}{2}x^2-2x+c+x^2-2x+1-2x+2-b\right]\sigma(x-2) + \left[-\frac{1}{2}x^2+x-\frac{1}{2}+2x-2-c\right]\sigma(x-3)\\
&= \left[\frac{1}{2}x^2+a\right]\sigma(x) + \left[-\frac{3}{2}x^2+3t-\frac{1}{2}+b-a\right]\sigma(x-1) + \left[\frac{3}{2}x^2-6t+3+c-b\right]\sigma(x-2)\\
&\qquad + \left[-\frac{1}{2}x^2+3t-\frac{5}{2}-c\right]\sigma(x-3)\\
&= \begin{cases}
\frac{1}{2}x^2+a & \text{for } x\in[0,\,1]\\
-x^2+3x-\frac{1}{2}+b & \text{for } x\in[1,\,2]\\
\frac{1}{2}x^2-3x+\frac{5}{2}+c & \text{for } x\in[2,\,3]\\
0 & \text{for } x\in\mathbb{R}\backslash[0,\,3]
\end{cases}
\end{aligned}
$$

The B-spline $N_3(x)$ is continuous. Therefore, at each point $n$ where $\sigma(x-n)$ jumps from $0$ to $1$, the factor $\mathsf{f}_n(x)$ in $\mathsf{f}_n(x)\sigma(x-n)$ must be $0$. We can use this to compute the boundary conditions $a$, $b$, and $c$:

$$
\begin{aligned}
\frac{1}{2}x^2+a\Big|_{t=0} &= 0 &&\implies 0+a=0 &&\implies a=0\\
-\frac{3}{2}x^2+3x-\frac{1}{2}+b-a\Big|_{t=1} &= 0 &&\implies -\frac{3}{2}+3-\frac{1}{2}+b-0=0 &&\implies b=-1\\
\frac{3}{2}x^2-6x+3+c-b\Big|_{t=2} &= 0 &&\implies \frac{12}{2}-12+3+c+1=0 &&\implies c=2\\
-\frac{1}{2}x^2+3x-\frac{5}{2}-c\Big|_{t=3} &= 0 &&\implies -\frac{9}{2}+9-\frac{5}{2}-c=0 &&\implies c=2
\end{aligned}
$$

✎

## 5.2 Properties

**Theorem 5.5** [79] *Let* $\mathbb{1}$ *be the* SET INDICATOR *function* *(Definition 1.3 page 3)*. *Let* $\sigma(x)$ *be the* STEP FUNCTION *(Definition 5.1 page 42)*.

$$\mathsf{N}_n(x) = \frac{1}{(n)!}\sum_{k=0}^{n+1}(-1)^k\binom{n+1}{k}(x-k)^n\sigma(x-k) \qquad \forall n \in \mathbb{W}\setminus 0$$

✎PROOF: Proof by induction:

(1) Proof for $n = 1$ case:

$$\begin{aligned}\mathsf{N}_1(x) &= x\sigma(x) - 2(x-1)\sigma(x-1) + (x-2)\sigma(x-2) && \text{by Lemma 5.4 page 42}\\ &= \frac{1}{(2-1)!}\sum_{k=0}^{2}(-1)^k\binom{2}{k}(x-k)^{2-1}\sigma(x-k)\end{aligned}$$

(2) Proof that $n$ case $\implies$ $n+1$ case:

$$\begin{aligned}&\mathsf{N}_{n+1}(x)\\ &= \int_0^1 \mathsf{N}_n(x-\tau)\,\mathrm{d}\tau && \text{by Lemma 5.3 page 42}\\ &= \int_0^1 \frac{1}{(n-1)!}\sum_{k=0}^{n}(-1)^k\binom{n}{k}(x-\tau-k)^{n-1}\sigma(x-\tau-k) && \text{by left hypothesis}\\ &= \frac{1}{(n-1)!}\sum_{k=0}^{n}(-1)^k\binom{n}{k}\left(\frac{-1}{n}\right)(x-\tau-k)^n\sigma(x-\tau-k)\bigg|_0^1\\ &= \frac{1}{(n)!}\sum_{k=0}^{n}(-1)^{k+1}\binom{n}{k}(x-\tau-k)^n\sigma(x-\tau-k)\bigg|_0^1\\ &= \underbrace{\left\{\frac{1}{(n)!}\sum_{k=0}^{n}(-1)^{k+1}\binom{n}{k}(x-k-1)^n\sigma(x-k-1)\right\}}_{\tau=1} - \underbrace{\left\{\frac{1}{(n)!}\sum_{k=0}^{n}(-1)^{k+1}\binom{n}{k}(x-k)^n\sigma(x-k)\right\}}_{\tau=0}\\ &= \left\{\frac{1}{(n)!}\sum_{m=1}^{n}(-1)^m\binom{n}{m-1}(x-m)^n\sigma(x-m)\right\} - \left\{\frac{1}{(n)!}\sum_{k=0}^{n}(-1)^{k+1}\binom{n}{k}(x-k)^n\sigma(x-k)\right\}\\ &\qquad \text{where } m \triangleq k+1 \implies k = m-1\\ &= \left\{\frac{1}{(n)!}\sum_{m=1}^{n}(-1)^m\left\{\binom{n+1}{m}-\binom{n}{m}\right\}(x-m)^n\sigma(x-m)\right\} + \left\{\frac{1}{(n)!}\sum_{k=0}^{n}(-1)^k\binom{n}{k}(x-k)^n\sigma(x-k)\right\}\\ &\qquad \text{by Stifel's formula Theorem 4.6 page 31}\end{aligned}$$

[79] [22], page 142, ⟨Theorem 6.1.3⟩, [25], page 84, ⟨(4.1.12)⟩

$$= \left\{ \frac{1}{(n)!} \sum_{m=1}^{n} (-1)^m \binom{n+1}{m} (x-m)^n \sigma(x-m) \right\} + \left\{ \frac{1}{(n)!} (-1)^0 \binom{n+1}{0} (x-0)^n \sigma(x-0) \right\}$$
$$= \frac{1}{(n)!} \sum_{k=0}^{n} (-1)^k \binom{n+1}{k} (x-k)^n \sigma(x-k)$$

✏

**Lemma 5.6** [80]

$$\frac{\mathrm{d}}{\mathrm{d}x}\mathsf{N}_n(x) = \mathsf{N}_{n-1}(x) - \mathsf{N}_{n-1}(x-1) \qquad \forall n \in \mathbb{W} \setminus \{1,2\},\ \forall x \in \mathbb{R}$$

✎PROOF:

(1) Proof using Fundamental Theorem of Calculus (FTC):

$$\begin{aligned}
\frac{\mathrm{d}}{\mathrm{d}x}\mathsf{N}_n(x) &= \frac{\mathrm{d}}{\mathrm{d}x}\int_0^1 \mathsf{N}_{n-1}(x-\tau)\,\mathrm{d}\tau && \text{by Lemma 5.3 page 42}\\
&= \frac{\mathrm{d}}{\mathrm{d}x}\int_{x-u=0}^{x-u=1} \mathsf{N}_{n-1}(u)(-1)\,\mathrm{d}u && \text{where } u \triangleq x-\tau \implies \tau = x-u\\
&= \frac{\mathrm{d}}{\mathrm{d}x}\int_{u=x-1}^{u=x} \mathsf{N}_{n-1}(u)\,\mathrm{d}u\\
&= \left\{ \frac{\mathrm{d}}{\mathrm{d}x}\int \mathsf{N}_{n-1}(u)\,\mathrm{d}u\Big|_{u=x} \right\} - \left\{ \frac{\mathrm{d}}{\mathrm{d}x}\int \mathsf{N}_{n-1}(u)\,\mathrm{d}u\Big|_{u=x-1} \right\} && \text{by FTC[81]}\\
&= \left\{ \mathsf{N}_{n-1}(x)\frac{\mathrm{d}}{\mathrm{d}x}(x) \right\} - \left\{ \mathsf{N}_{n-1}(x-1)\frac{\mathrm{d}}{\mathrm{d}x}(x-1) \right\} && \text{by \textit{Chain Rule}[82]}\\
&= \mathsf{N}_{n-1}(x) - \mathsf{N}_{n-1}(x-1)
\end{aligned}$$

(2) Proof by induction:

[80] 📚 [61], page 25, ⟨3.2⟩, 📚 [99], page 121, ⟨Theorem 4.16⟩
[81] 📚 [60], page 163, ⟨Theorem 4.4.3⟩
[82] 📚 [60], pages 73–74, ⟨Theorem 3.1.2⟩

(a) Proof for $n=2$ case:

$$\begin{aligned}
&\mathsf{N}_1(x) - \mathsf{N}_1(x-1)\\
&= \underbrace{x\sigma(x) - 2(x-1)\sigma(x-1) + (x-2)\sigma(x-2)}_{\mathsf{N}_1(x)}\\
&\qquad - \underbrace{[(x-1)\sigma(x-1) - 2(x-2)\sigma(x-2) + (x-3)\sigma(x-3)]}_{\mathsf{N}_1(x-1)} && \text{by Lemma 5.4 page 42}\\
&= x\sigma(x) + [-2x+2-x+1]\sigma(x-1) + [x-2+2x-4]\sigma(x-2) + [-x+3]\sigma(x-3)\\
&= x\sigma(x) + [-3x+3]\sigma(x-1) + [3x-6]\sigma(x-2) + [-x+3]\sigma(x-3)\\
&= \frac{\mathsf{d}}{\mathsf{dx}}\left\{\begin{array}{l} \frac{1}{2}x^2\sigma(x) + \left[-\frac{3}{2}x^2+3x-\frac{1}{2}\right]\sigma(x-1) + \left[\frac{3}{2}x^2-6x+3\right]\sigma(x-2)\\ \qquad +\left[-\frac{1}{2}x^2+3x-\frac{5}{2}\right]\sigma(x-3)\end{array}\right\}\\
&= \frac{\mathsf{d}}{\mathsf{dx}}\mathsf{N}_2(x) && \text{by Lemma 5.4 page 42}
\end{aligned}$$

(b) Proof that $n$ case $\implies$ $n+1$ case:

$$\begin{aligned}
\frac{\mathsf{d}}{\mathsf{dx}}\mathsf{N}_{n+1}(x) &= \frac{\mathsf{d}}{\mathsf{dx}}\int_0^1 \mathsf{N}_n(x-\tau)\,\mathsf{d}\tau && \text{by Lemma 5.3 page 42}\\
&= \int_0^1 \frac{\partial}{\partial \mathsf{t}}\mathsf{N}_n(x-\tau)\,\mathsf{d}\tau && \text{see note later}\\
&= \int_0^1 \left[\mathsf{N}_{n-1}(x-\tau) - \mathsf{N}_{n-1}(x-1-\tau)\right]\,\mathsf{d}\tau && \text{by left hypothesis}\\
&= \int_0^1 \mathsf{N}_{n-1}(x-\tau)\,\mathsf{d}\tau - \int_0^1 \mathsf{N}_{n-1}(x-1-\tau)\,\mathsf{d}\tau\\
&= \mathsf{N}_n(x) - \mathsf{N}_n(x-1) && \text{by Lemma 5.3 page 42}
\end{aligned}$$

Note: For information about differentiation of an integral, see [37], [106], [73], page 389, ⟨Chapter VII⟩

✏

**Theorem 5.7** [83] *Let* $\operatorname{supp}\mathsf{f}$ *be the* SUPPORT *of a function* $\mathsf{f}$.

[83] [22], page 140, ⟨Theorem 6.1.1⟩, [61], page 27, ⟨3.4⟩, [99], page 120, ⟨Theorem 4.15⟩, [32], page 90, ⟨B-Spline Property (i)⟩, [24], page 2, ⟨Theorem 1.1⟩, [116], page 53, ⟨Theorem 3.7⟩, [27], [31]

| | | | | | | |
|---|---|---|---|---|---|---|
| 1. | $\mathsf{N}_n(x)$ | $\geq$ | $0$ | $\forall n \in \mathbb{W},$ | $\forall x \in \mathbb{R}$ | (POSITIVE) |
| 2. | $\operatorname{supp} \mathsf{N}_n(x)$ | $=$ | $[0,\ n+1]$ | $\forall n \in \mathbb{W}$ | | (COMPACT SUPPORT) |
| 3. | $\int_{\mathbb{R}} \mathsf{N}_n(x)\,\mathrm{d}x$ | $=$ | $1$ | $\forall n \in \mathbb{W}$ | | (UNIT AREA) |
| 4. | $\sum_{k\in\mathbb{Z}} \mathsf{N}_n(x-k)$ | $=$ | $1$ | $\forall n \in \mathbb{W}\backslash 0$ | | (PARTITION OF UNITY) |
| 5. | $\mathsf{N}_n(x)$ | $=$ | $\frac{x}{n}\mathsf{N}_{n-1}(x) + \frac{n+1-x}{n}\mathsf{N}_{n-1}(x-1)$ | $\forall n \in \mathbb{W}\backslash\{1\},$ | $\forall x \in \mathbb{R}$ | |
| 6. | $\mathsf{N}_n\left(\frac{n+1}{2}+x\right)$ | $=$ | $\mathsf{N}_n\left(\frac{n+1}{2}-x\right)$ | $\forall n \in \mathbb{W}$ | $\forall x \in \mathbb{R}$ | (SYMMETRIC) |

✎PROOF:

(1) Proof that $\operatorname{supp} \mathsf{N}_n(x) \geq 0$ (proof by induction):

(a) Proof that $\mathsf{N}_0(x) \geq 0$: by Definition 5.2 page 42

(b) Proof that $\mathsf{N}_n \geq 0 \implies \mathsf{N}_{n+1} \geq 0$:

$$
\begin{aligned}
&\mathsf{N}_{n+1}(x) \\
&\qquad = \int_{\tau=0}^{\tau=1} \mathsf{N}_n(x-\tau)\,\mathrm{d}\tau && \text{by Lemma 5.3 page 42} \\
&\qquad \geq 0 && \text{by left hypothesis}
\end{aligned}
$$

(2) Proof that $\operatorname{supp} \mathsf{N}_n(x) = [0,\ n]$ (proof by induction):

(a) Proof that $\operatorname{supp} \mathsf{N}_0 = [0,\ 1]$: by Definition 5.2 page 42

(b) Proof that $\operatorname{supp} \mathsf{N}_n = [0,\ n] \implies \operatorname{supp} \mathsf{N}_{n+1} = [0,\ n+1]$:

$$
\begin{aligned}
\operatorname{supp} \mathsf{N}_{n+1}(x) &= \operatorname{supp} \int_{\tau=0}^{\tau=1} \mathsf{N}_n(x-\tau)\,\mathrm{d}\tau && \text{by Lemma 5.3 page 42} \\
&= \left\{x \in \mathbb{R} \,\middle|\, x-\tau \in [0,\ n] \text{ for some } \tau \in [0,\ 1]\right\} && \text{by left hypothesis} \\
&= [0,\ n+1]
\end{aligned}
$$

(3) Proof that $\int_{\mathbb{R}} \mathsf{N}_n(x)\,\mathrm{d}x = 1$ (proof by induction):

(a) Proof that $\int_{\mathbb{R}} \mathsf{N}_1(x) = 1$:

$$
\int_{\mathbb{R}} \mathsf{N}_0(x)\,\mathrm{d}x = 0 \qquad \text{by definition of } \mathsf{N}_1 \text{ (Definition 5.2 page 42)}
$$

(b) Proof that $\int_{\mathbb{R}} \mathsf{N}_n(x) = 1 \implies \int_{\mathbb{R}} \mathsf{N}_{n+1} = 1$:

$$\begin{aligned}
\int_{\mathbb{R}} \mathsf{N}_{n+1}(x)\,\mathrm{d}x &= \int_{\mathbb{R}} \int_0^1 \mathsf{N}_n(x-\tau)\,\mathrm{d}\tau\,\mathrm{d}x && \text{by Lemma 5.3 page 42}\\
&= \int_0^1 \int_{\mathbb{R}} \mathsf{N}_n(x-\tau)\,\mathrm{d}x\,\mathrm{d}\tau\\
&= \int_0^1 \int_{\mathbb{R}} \mathsf{N}_n(u)\,\mathrm{d}u\,\mathrm{d}\tau && \text{where } u \triangleq x-\tau \implies \tau = x-u\\
&= \int_0^1 1\,\mathrm{d}\tau && \text{by left hypothesis}\\
&= 1
\end{aligned}$$

(4) Proof that $\sum_{k\in\mathbb{Z}} \mathsf{N}_n(x-k) = 1$ for $n \in \mathbb{W}\backslash 0$ (proof by induction):

(a) Proof that $\sum_{k\in\mathbb{Z}} \mathsf{N}_1(x-k) = 1$:

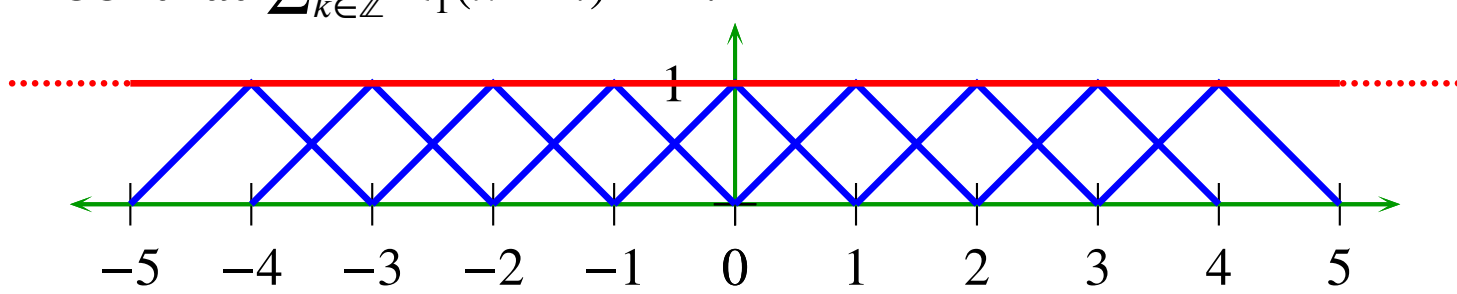

(b) Proof that $\sum_{k\in\mathbb{Z}} \mathsf{N}_n(x-k) = 1 \implies \sum_{k\in\mathbb{Z}} \mathsf{N}_{n+1}(x-k) = 1$:

$$\begin{aligned}
\sum_{k\in\mathbb{Z}} \mathsf{N}_{n+1}(x-k) &= \sum_{k\in\mathbb{Z}} \int_{\tau=0}^{\tau=1} \mathsf{N}_n(x-k-\tau)\,\mathrm{d}\tau && \text{by Lemma 5.3 page 42}\\
&= \sum_{k\in\mathbb{Z}} \int_{x-u=0}^{x-u=1} \mathsf{N}_n(u-k)(-1)\,\mathrm{d}u && \text{where } u \triangleq x-\tau \implies \tau = x-u\\
&= \sum_{k\in\mathbb{Z}} \int_{u=x-1}^{u=t} \mathsf{N}_n(u-k)\,\mathrm{d}u\\
&= \int_{u=x-1}^{u=t} \left(\sum_{k\in\mathbb{Z}} \mathsf{N}_n(u-k)\right)\mathrm{d}u\\
&= \int_{u=x-1}^{u=t} 1\,\mathrm{d}\tau && \text{by left hypothesis}\\
&= 1
\end{aligned}$$

(5) Proof for recursion equation (proof by induction):

(a) Proof for $n = 1$ case:

$$\begin{aligned}
\frac{x}{1}\mathsf{N}_0(x) + \frac{1+1-x}{1}\mathsf{N}_0(x-1) &= \frac{x}{1}\underbrace{[\sigma(x)-\sigma(x-1)]}_{\mathsf{N}_0(x)} + \frac{1+1-x}{1}\underbrace{[\sigma(x-1)-\sigma(x-2)]}_{\mathsf{N}_0(x-1)}\\
&= x\sigma(x) + [-x-x+2]\sigma(x-1) + [x-2]\sigma(x-2)\\
&= \mathsf{N}_1(x) \qquad \text{by Lemma 5.4 page 42}
\end{aligned}$$

(b) Proof that $n$ case $\implies$ $n+1$ case:

$$\begin{aligned}
&\frac{x}{n+1}\mathsf{N}_n(x) + \frac{n+2-x}{n+1}\mathsf{N}_n(x-1) + c_1 \\
&= \int \tfrac{\mathrm{d}}{\mathrm{d}x}\left\{\frac{x}{n+1}\mathsf{N}_n(x) + \frac{n+2-x}{n+1}\mathsf{N}_n(x-1)\right\}\,\mathrm{d}x \\
&= \int \underbrace{\frac{1}{n+1}\mathsf{N}_n(x) + \frac{x}{n+1}\tfrac{\mathrm{d}}{\mathrm{d}x}\mathsf{N}_n(x)}_{\frac{\mathrm{d}}{\mathrm{d}x}\frac{x}{n+1}\mathsf{N}_n(x)} + \underbrace{\frac{-1}{n+1}\mathsf{N}_n(x-1) + \frac{n+2-x}{n}\tfrac{\mathrm{d}}{\mathrm{d}x}\mathsf{N}_n(x-1)}_{\frac{\mathrm{d}}{\mathrm{d}x}\frac{n+2-x}{n+1}\mathsf{N}_n(x-1)}\,\mathrm{d}x \\
&\qquad \text{by \textit{product rule}} \\
&= \int \frac{1}{n+1}\underbrace{\left[\frac{x}{n}\mathsf{N}_{n-1}(x) + \frac{n+1-x}{n}\mathsf{N}_{n-1}(x-1)\right]}_{\text{by } n \text{ hypothesis}} + \frac{x}{n+1}\underbrace{\left[\mathsf{N}_{n-1}(x) - \mathsf{N}_{n-1}(x-1)\right]}_{\text{by Lemma 5.6 page 46}} \\
&\qquad - \underbrace{\left[\frac{x-1}{n^2+n}\mathsf{N}_{n-1}(x-1) + \frac{n-x+2}{n(n+1)}\mathsf{N}_{n-1}(x-2)\right]}_{\text{by } n \text{ hypothesis}} \\
&\qquad + \frac{n+2-x}{n+1}\underbrace{\left[\mathsf{N}_{n-1}(x-1) - \mathsf{N}_{n-1}(x-2)\right]}_{\text{by Lemma 5.6 page 46}}\,\mathrm{d}x \\
&= \int \left[\frac{x}{n(n+1)} + \frac{x}{n+1}\right]\mathsf{N}_{n-1}(x) + \left[\frac{n-x+1}{n(n+1)} - \frac{x-1}{n(n+1)} + \frac{n+2-2x}{n+1}\right]\mathsf{N}_{n-1}(x-1) \\
&\qquad + \left[\frac{-n-2+x}{n(n+1)} + \frac{-n-2+x}{n+1}\right]\mathsf{N}_{n-1}(x-2)\,\mathrm{d}x \\
&= \int \left[\frac{x+nx}{n(n+1)}\right]\mathsf{N}_{n-1}(x) + \left[\frac{n+2-2x+n(n+2-2x)}{n(n+1)}\right]\mathsf{N}_{n-1}(x-1) \\
&\qquad + \left[\frac{-n-2+x+n(-n-2+x)}{n(n+1)}\right]\mathsf{N}_{n-1}(x-2)\,\mathrm{d}x \\
&= \int \left[\frac{x}{n}\right]\mathsf{N}_{n-1}(x) + \left[\frac{n+2-2x}{n}\right]\mathsf{N}_{n-1}(x-1) + \left[\frac{-n-2+x}{n}\right]\mathsf{N}_{n-1}(x-2)\,\mathrm{d}x \\
&= \int \underbrace{\left[\frac{x}{n}\right]\mathsf{N}_{n-1}(x) + \left[\frac{n+1-x}{n}\right]\mathsf{N}_{n-1}(x-1)}_{\mathsf{N}_n(x)} \\
&\qquad \underbrace{-\left[\frac{x-1}{n}\right]\mathsf{N}_{n-1}(x-1) - \left[\frac{n+2-x}{n}\right]\mathsf{N}_{n-1}(x-2)}_{\mathsf{N}_{n-1}(x-1)}\,\mathrm{d}x \\
&= \int \mathsf{N}_n(x) - \mathsf{N}_n(x-1)\,\mathrm{d}x \qquad \text{by } n \text{ hypothesis} \\
&= \int \tfrac{\mathrm{d}}{\mathrm{d}x}\mathsf{N}_{n+1}(x)\,\mathrm{d}x \qquad \text{by Lemma 5.6 page 46}
\end{aligned}$$

$$= \mathsf{N}_{n+1}(x) + c_2$$

Proof that $c_1 = c_2$: By item 2 (page 48), $\mathsf{N}_n(x) = 0$ for $x < 0$. Therefore, $c_1 = c_2$.

(6) Proof for symmetric equation (proof by induction):
Note that it is true for $\mathsf{N}_0(x)$. Then here is the proof that $n-1$ case $\implies$ $n$ case …

$$\begin{aligned}
&\mathsf{N}_n\left(\frac{n+1}{2} + x\right) \\
&= \frac{\frac{n+1}{2} + x}{n}\mathsf{N}_{n-1}\left(\frac{n+1}{2} + x\right) + \frac{n+1-\left(\frac{n+1}{2}+x\right)}{n}\mathsf{N}_{n-1}\left(\frac{n+1}{2} + x - 1\right) && \text{by item 5 page 49} \\
&= \frac{n+1-\left(\frac{n+1}{2}-x\right)}{n}\mathsf{N}_{n-1}\left(\frac{n}{2} + \left[x + \frac{1}{2}\right]\right) + \frac{\frac{n+1}{2} - x}{n}\mathsf{N}_{n-1}\left(\frac{n}{2} + \left[x - \frac{1}{2}\right]\right) \\
&= \frac{n+1-\left(\frac{n+1}{2}-x\right)}{n}\mathsf{N}_{n-1}\left(\frac{n}{2} - \left[x + \frac{1}{2}\right]\right) + \frac{\frac{n+1}{2} - x}{n}\mathsf{N}_{n-1}\left(\frac{n}{2} - \left[x - \frac{1}{2}\right]\right) && \text{by left hypothesis} \\
&= \frac{n+1-\left(\frac{n+1}{2}-x\right)}{n}\mathsf{N}_{n-1}\left(\left[\frac{n+1}{2} - x\right] - 1\right) + \frac{\frac{n+1}{2} - x}{n}\mathsf{N}_{n-1}\left(\frac{n+1}{2} - x\right) \\
&= \mathsf{N}_n\left(\frac{n+1}{2} - x\right) && \text{by item 5 page 49}
\end{aligned}$$

✏

**Theorem 5.8** [84] *Let* $\mathsf{f}$ *be a continuous function in* $\boldsymbol{L}^2_{\mathbb{R}}$ *and* $\mathsf{f}^{(n)}$ *the* $n$*th derivative of* $\mathsf{f}$.

$$\int_{[0,1]^n} \mathsf{f}^{(n)}\left(\sum_{k=1}^{n} x_k\right) \mathrm{d}x_1\, \mathrm{d}x_2 \cdots \mathrm{d}x_n = \sum_{k=0}^{n} (-1)^{n-k}\binom{n}{k}\mathsf{f}(k) \qquad \forall n \in \mathbb{N}$$

✎PROOF: Proof by induction:

(1) Proof for $n = 1$ case:

$$\begin{aligned}
\int_{[0,1]} \mathsf{f}^{(1)}(x)\, \mathrm{d}x &= \mathsf{f}(x)\big|_0^1 \\
&= \mathsf{f}(1) - \mathsf{f}(0) \\
&= (-1)^{1+1}\binom{1}{1}\mathsf{f}(1) + (-1)^{1+0}\binom{1}{0}\mathsf{f}(0) \\
&= \sum_{k=0}^{n}(-1)^{n-k}\binom{n}{k}\mathsf{f}(k)
\end{aligned}$$

[84] 📖 [25], page 86, ⟨item (ii)⟩

(2) Proof that $n$ case $\implies$ $n+1$ case:

$$\int_{[0,1]^{n+1}} \mathsf{f}^{(n+1)}\left(\sum_{k=1}^{n} x_k\right) \mathrm{d}x_1\,\mathrm{d}x_2\cdots\mathrm{d}x_{n+1}$$
$$= \int_{[0,1]^{n}} \left\{ \mathsf{f}^{(n)}\left(x_{n+1} + \sum_{k=1}^{n} x_k\right)\Bigg|_0^1 \right\} \mathrm{d}x_1\,\mathrm{d}x_2\cdots\mathrm{d}x_{n}$$
$$= \int_{[0,1]^{n}} \mathsf{f}^{(n)}\left(1 + \sum_{k=1}^{n} x_k\right) - \mathsf{f}^{(n)}\left(0 + \sum_{k=1}^{n} x_k\right) \mathrm{d}x_1\,\mathrm{d}x_2\cdots\mathrm{d}x_{n}$$
$$= \sum_{k=0}^{n}(-1)^{n-k}\binom{n}{k}\mathsf{f}(k+1) - \sum_{k=0}^{n}(-1)^{n-k}\binom{n}{k}\mathsf{f}(k) \qquad \text{by } n \text{ case hypothesis}$$
$$= \sum_{k=1}^{n+1}(-1)^{n-k+1}\binom{n}{k-1}\mathsf{f}(k) + \sum_{k=0}^{n}(-1)(-1)^{n-k}\binom{n}{k}\mathsf{f}(k)$$
$$= \left\{ \mathsf{f}(n+1) + \sum_{k=1}^{n}(-1)^{n-k+1}\binom{n}{k-1}\mathsf{f}(k) \right\} + \left\{ (-1)^{n+1}\mathsf{f}(0) + \sum_{k=1}^{n}(-1)^{n-k+1}\binom{n}{k}\mathsf{f}(k) \right\}$$
$$= \mathsf{f}(n+1) + (-1)^{n+1}\mathsf{f}(0) + \sum_{k=1}^{n}(-1)^{n-k+1}\left[\binom{n}{k-1} + \binom{n}{k}\right]\mathsf{f}(k)$$
$$= \mathsf{f}(n+1) + (-1)^{n+1}\mathsf{f}(0) + \sum_{k=1}^{n}(-1)^{n-k+1}\binom{n+1}{k}\mathsf{f}(k) \qquad \text{by Pascal's Recursion Theorem 4.6 page 31}$$
$$= \sum_{k=0}^{n+1}(-1)^{n-k+1}\binom{n+1}{k}\mathsf{f}(k)$$

✏

**Theorem 5.9** [85] *Let* $\mathsf{f}$ *be a continuous function in* $\boldsymbol{L}^2_{\mathbb{R}}$.

1. $$\int_{\mathbb{R}} \mathsf{f}(x)\mathsf{N}_n(x)\,\mathrm{d}x = \int_{[0,1]^{n+1}} \mathsf{f}(x_0 + x_1 + \cdots + x_n)\,\mathrm{d}x_0\,\mathrm{d}x_1\cdots\mathrm{d}x_n$$
2. $$\int_{\mathbb{R}} \mathsf{f}^{(n)}(x)\mathsf{N}_n(x)\,\mathrm{d}x = \sum_{k=0}^{n}(-1)^{n-k}\binom{n}{k}\mathsf{f}(k)$$

✎PROOF:

(1) Proof for (1) (proof by induction):

(a) Proof for $n=0$ case:

$$\int_{\mathbb{R}} \mathsf{N}_0(x)\mathsf{f}(x)\,\mathrm{d}x = \int_{[0,1]} \mathsf{f}(x)\,\mathrm{d}x$$

[85] [25], page 85, ⟨(4.2.2), (4.2.3)⟩, [22], page 140, ⟨Theorem 6.1.1⟩

(b) Proof that $\mathsf{N}_n$ case $\implies$ $\mathsf{N}_{n+1}$ case:

$$\begin{aligned}
&\int_{\mathbb{R}} \mathsf{N}_{n+1}(x)\mathsf{f}(x)\,\mathrm{d}x \\
&= \int_{\mathbb{R}} \left(\int_{[0,\,1]} \mathsf{N}_n(x-\tau)\,\mathrm{d}\tau\right)\mathsf{f}(x)\,\mathrm{d}x && \text{by Lemma 5.3 page 42} \\
&= \int_{[0,\,1]}\int_{\mathbb{R}} \mathsf{N}_n(x-\tau)\mathsf{f}(x)\,\mathrm{d}x\,\mathrm{d}\tau \\
&= \int_{[0,\,1]}\int_{\mathbb{R}} \mathsf{N}_n(u)\mathsf{f}(u+\tau)\,\mathrm{d}u\,\mathrm{d}\tau && \text{where } u \triangleq x-\tau \implies x = u+\tau \\
&= \int_{[0,\,1]}\int_{[0,\,1]^n} \mathsf{f}(u_0+u_1+\cdots+u_n+\tau)\,\mathrm{d}u_0\,\mathrm{d}u_1\cdots\mathrm{d}u_n\,\mathrm{d}\tau && \text{by left hypothesis} \\
&= \int_{[0,\,1]^{n+1}} \mathsf{f}(u_0+u_1+\cdots+u_n)\,\mathrm{d}u_0\,\mathrm{d}u_1\cdots\mathrm{d}u_n\,\mathrm{d}\tau \\
&= \int_{[0,\,1]^{n+1}} \mathsf{f}(x_0+x_1+\cdots+x_n+x_{n+1})\,\mathrm{d}x_0\,\mathrm{d}x_1\cdots\mathrm{d}x_n\,\mathrm{d}x_{n+1} && \text{by change of variables}
\end{aligned}$$

(2) Proof for (2):

$$\begin{aligned}
\int_{\mathbb{R}} \mathsf{f}^{(n)}(x)\mathsf{N}_n(x)\,\mathrm{d}x &= \int_{[0,\,1]^{n+1}} \mathsf{f}^{(n)}\left(\sum_{k=0}^{n} x_k\right)\mathrm{d}x_0\,\mathrm{d}x_1\cdots\mathrm{d}x_n && \text{by item 1} \\
&= \sum_{k=0}^{n}(-1)^{n-k}\binom{n}{k}\mathsf{f}(k) && \text{by Theorem 5.8 page 51}
\end{aligned}$$

✏

**Theorem 5.10** [86] *Let* $\tilde{\mathbf{F}}$ *be the* FOURIER TRANSFORM *operator (Definition 1.26 page 8).*

$$\tilde{\mathbf{F}}\mathsf{N}_n(\omega) = \frac{1}{\sqrt{2\pi}}\left(\frac{1-e^{-i\omega}}{i\omega}\right)^{n+1} = \frac{1}{\sqrt{2\pi}}e^{-i\frac{(n+1)\omega}{2}}\underbrace{\left(\frac{\sin\left(\frac{\omega}{2}\right)}{\frac{\omega}{2}}\right)}_{\operatorname{sinc}\left(\frac{\omega}{2}\right)}^{n+1}$$

[86] [22], page 142, ⟨Corollary 6.1.2⟩

✎PROOF:

$$\begin{aligned}
\tilde{\mathbf{F}}\mathsf{N}_n(\omega) &= \frac{1}{\sqrt{2\pi}} \int_{\mathbb{R}} \mathsf{N}_n(x) e^{-i\omega x}\,\mathrm{d}x && \text{by definition of } \tilde{\mathbf{F}} \text{ Definition 1.26 page 8}\\
&= \frac{1}{\sqrt{2\pi}} \int_{[0,\,1]^{n+1}} e^{-i\omega(x_0+x_1+\cdots x_n)}\,\mathrm{d}x_0\,\mathrm{d}x_1\,\cdots,\,\mathrm{d}x_n && \text{by Theorem 5.9 page 52}\\
&= \frac{1}{\sqrt{2\pi}} \prod_{k=0}^{n} \left( \int_{[0,\,1]} e^{-i\omega x_k}\,\mathrm{d}x_k \right) && \text{by Theorem 5.9 page 52}\\
&= \frac{1}{\sqrt{2\pi}} \left( \int_{[0,\,1]} e^{-i\omega x}\,\mathrm{d}x \right)^{n+1}\\
&= \frac{1}{\sqrt{2\pi}} \left( \left.\frac{e^{-i\omega x}}{-i\omega}\right|_0^1 \right)^{n+1}\\
&= \frac{1}{\sqrt{2\pi}} \left( \frac{1-e^{-i\omega}}{i\omega} \right)^{n+1}\\
&= \frac{1}{\sqrt{2\pi}} \left[ e^{-i\frac{\omega}{2}} \left( \frac{e^{i\frac{\omega}{2}} - e^{-i\frac{\omega}{2}}}{i\omega} \right) \right]^{n+1}\\
&= \frac{1}{\sqrt{2\pi}} \left[ e^{-i\frac{\omega}{2}} \left( \frac{2i\sin\left(\frac{\omega}{2}\right)}{\frac{2i\omega}{2}} \right) \right]^{n+1} && \text{by } \textit{Euler formulas} \text{ (Corollary 1.14 page 5)}\\
&= \frac{1}{\sqrt{2\pi}} e^{-i\frac{(n+1)\omega}{2}} \left( \frac{\sin\left(\frac{\omega}{2}\right)}{\frac{\omega}{2}} \right)^{n+1}
\end{aligned}$$

✏

## 5.3 Spline function spaces

**Definition 5.11** [87] Let $\mathsf{N}_n(x)$ be an *n th order cardinal B-spline* (Definition 5.2 page 42). The **space of all splines of order n** is denoted $\boldsymbol{S}^n(a\mathbb{Z})$ and is defined as

$$\boldsymbol{S}^n(a\mathbb{Z}) \triangleq \operatorname{span}\left\{ \mathbf{T}^m \mathsf{N}_n(ax) \,|\, m\in\mathbb{Z} \right\}.$$

**Theorem 5.12** [88] *Let $\boldsymbol{S}^n(\mathbb{Z})$ be the* SPACE OF ALL SPLINES OF ORDER N *(Definition 5.11 page 54).*

$$\underbrace{\left\{ \mathsf{f}(x) = \sum_{k\in\mathbb{Z}} \alpha_k \mathbf{T}^n \mathsf{N}_n(x-k) = \sum_{k\in\mathbb{Z}} \beta_k \mathbf{T}^n \mathsf{N}_n(x-k) \right\} \implies \left\{ (\!(\alpha_k)\!)_{k\in\mathbb{Z}} = (\!(\beta_k)\!)_{k\in\mathbb{Z}} \right\}}_{\textit{coefficients are } \text{UNIQUE}}$$

---

87 [116], page 52, ⟨Definition 3.5⟩
88 [116], page 55, ⟨Theorem 3.11⟩

**Lemma 5.13** [89] *Let $\boldsymbol{S}^n(\mathbb{Z})$ be the* SPACE OF ALL SPLINES OF ORDER N *(Definition 5.11 page 54). For each $n \in \mathbb{W}$,*
$\left(\!\left(\mathbf{T}^n \mathsf{N}_n(x)\right)\!\right)_{n\in\mathbb{Z}}$ *is a* RIESZ BASIS *in* $\boldsymbol{L}^2_{\mathbb{R}}$.

## 5.4 Examples

**Example 5.14** (Square pulse)

The B-Spline $\mathsf{N}_0(x)$ is calculated in Lemma 5.4 page 42 and illustrated to the right.

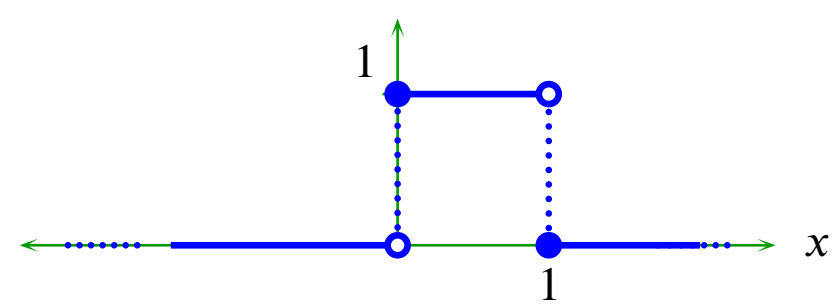

The B-spline $\mathsf{N}_0(x)$ forms a *partition of unity* (cross reference Theorem 5.7 page 47).

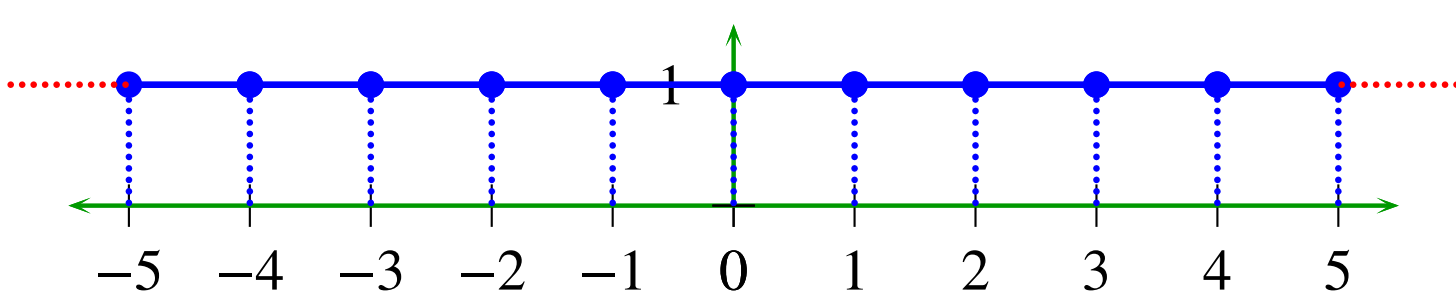

Here is the Fourier transform $\left[\tilde{\mathbf{F}}\mathsf{f}\right](\omega)$ of $\mathsf{N}_0(x)$:

$$
\begin{aligned}
\tilde{\mathbf{F}}\mathsf{f}(x) &\triangleq \frac{1}{\sqrt{2\pi}} \int_0^1 e^{-i\omega x}\,\mathrm{d}x && \text{by definition of } \tilde{\mathbf{F}} \text{ page 8}\\
&= \frac{1}{-i\omega}\frac{1}{\sqrt{2\pi}} e^{-i\omega x}\Big|_0^1 \\
&= \frac{1}{-i\omega}\frac{1}{\sqrt{2\pi}}\left(e^{-i\omega\frac{1}{2}} - e^{i\omega\frac{1}{2}}\right)e^{-i\omega\frac{1}{2}} \\
&= \frac{1}{-i\omega}\frac{1}{\sqrt{2\pi}}\left[-2i\sin\left(\frac{\omega}{2}\right)\right]e^{-i\omega\frac{1}{2}} && \text{by Corollary 1.14 page 5}\\
&= \frac{2}{2}\frac{1}{\sqrt{2\pi}}\frac{\sin\left(\frac{\omega}{2}\right)}{\omega\frac{1}{2}}e^{-i\omega\frac{1}{2}} \\
&= \frac{1}{\sqrt{2\pi}}\frac{\sin\left(\frac{\omega}{2}\right)}{\frac{\omega}{2}}e^{-i\omega\frac{1}{2}}
\end{aligned}
$$

Note that $\tilde{\mathbf{F}}\mathsf{N}_0(0) = \frac{1}{\sqrt{2\pi}}$, which agrees with the result demonstrated in Theorem 5.10 page 53.

[89] [116], page 56, ⟨Proposition 3.12⟩

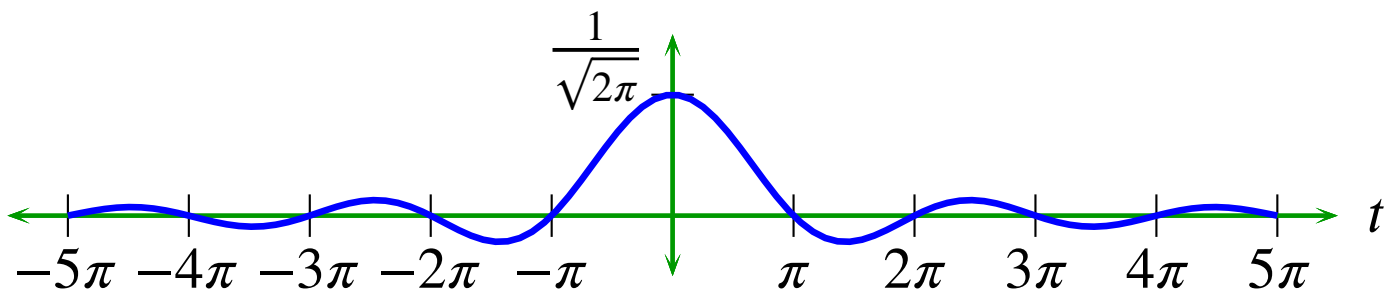

**Example 5.15** [90]

The B-Spline $\mathsf{N}_1(x)$ is calculated in Lemma 5.4 page 42 and illustrated to the right.

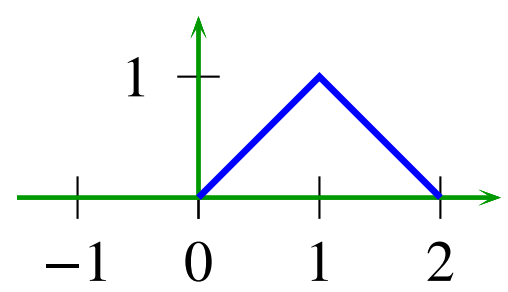

B-spline $\mathsf{N}_1(x)$ forms a *partition of unity* (cross reference Theorem 5.7 page 47).

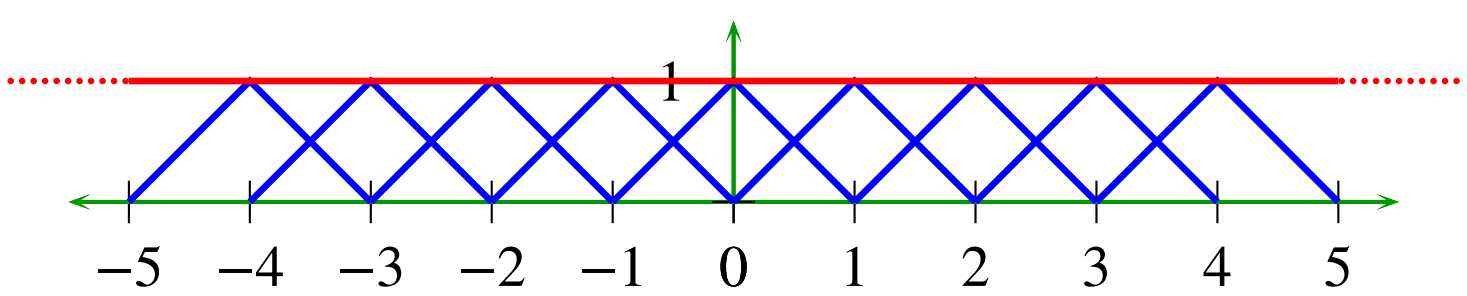

Here is the Fourier transform $\left[\tilde{\mathbf{F}}\mathsf{N}_1\right](\omega)$ of the function $\mathsf{N}_1(x)$:

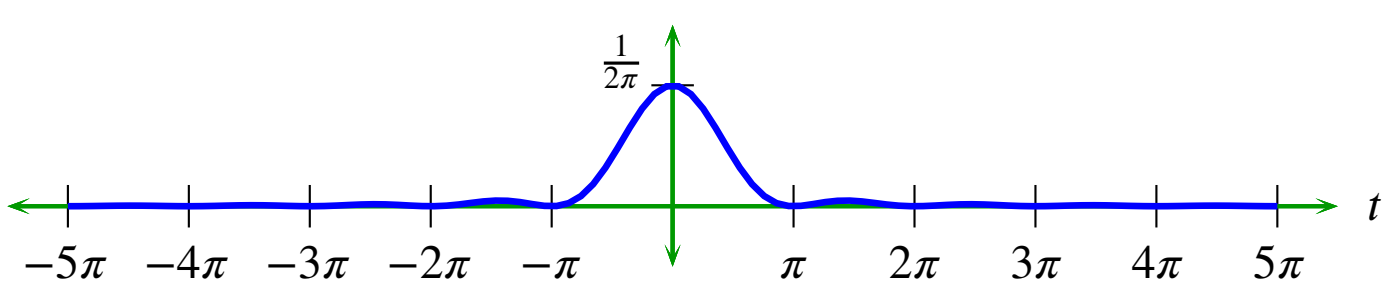

**Example 5.16** (centered cubic B-spline) [91] Let a function $\mathsf{f}$ be the *centered cubic B-spline* defined as follows:

$$\mathsf{f}(x) \triangleq \begin{cases} \frac{2}{3} - \frac{1}{2}|x|^2(2-|x|) & \text{for } |x| < 1 \\ \frac{1}{6}(2-|x|)^3 & \text{for } 1 \le |x| < 2 \\ 0 & \text{otherwise} \end{cases}$$

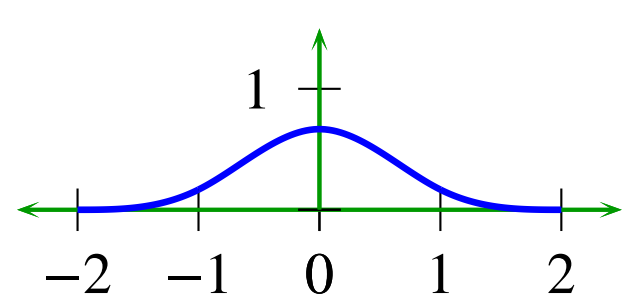

Then $\mathsf{f}$ forms a *partition of unity* because $\sum_{n\in\mathbb{Z}} \mathsf{f}(x-n) = 1$.

[90] [22], pages 146–147, ⟨Corollary 6.2.1⟩
[91] [22], page 146, ⟨Corollary 6.2.1⟩, [13], page 479, [32]

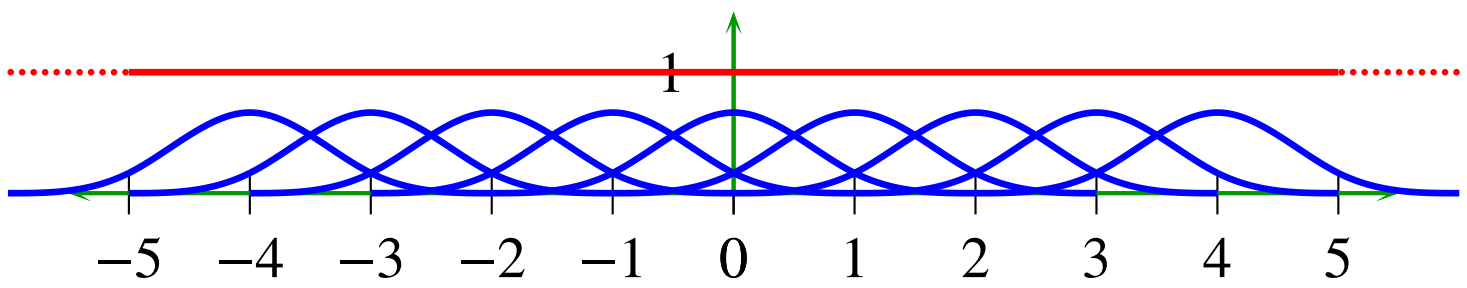

✎PROOF: Note that the function $\mathsf{h}(x) \triangleq \sum_{n\in\mathbb{Z}} \mathsf{f}(x-n)$ is periodic with period 1 (Proposition 2.3 page 15). So it is only necessary to examine a single interval of length one. Here we use the interval $[0,\ 1]$. In this interval, there are four functions contributing to the sum $\sum_{n\in\mathbb{Z}} \mathsf{f}(x-n)$ (see previous illustration).

$$\begin{aligned}
\sum_{n=-1}^{n=2} \mathsf{f}(x-n) &= \underbrace{\frac{1}{6}(2-|x+1|)^3}_{\mathsf{f}(x+1)} + \underbrace{\frac{2}{3}-\frac{1}{2}|x|^2(2-|x|)}_{\mathsf{f}(x)} + \underbrace{\frac{2}{3}-\frac{1}{2}|x-1|^2(2-|x-1|)}_{\mathsf{f}(x-1)} + \underbrace{\frac{1}{6}(2-|x-2|)^3}_{\mathsf{f}(x-2)} \\
&= \underbrace{\frac{1}{6}(2-(x+1))^3}_{\mathsf{f}(x+1)} + \underbrace{\frac{2}{3}-\frac{1}{2}x^2(2-x)}_{\mathsf{f}(x)} + \underbrace{\frac{2}{3}-\frac{1}{2}(1-x)^2(2-(1-x))}_{\mathsf{f}(x-1)} + \underbrace{\frac{1}{6}(2-(2-x))^3}_{\mathsf{f}(x-2)} \\
&= \underbrace{\frac{1}{6}\left(-x^3+3x^2-3x+1\right)}_{\mathsf{f}(x+1)} + \underbrace{\frac{2}{3}-\frac{1}{2}\left(-x^3+2x^2\right)}_{\mathsf{f}(x)} + \underbrace{\frac{2}{3}-\frac{1}{2}\left(x^2-2x+1\right)(x+1)}_{\mathsf{f}(x-1)} + \underbrace{\frac{1}{6}x^3}_{\mathsf{f}(x-2)} \\
&= \underbrace{\frac{1}{6}\left(-x^3+3x^2-3x+1\right)}_{\mathsf{f}(x+1)} + \underbrace{\frac{2}{3}-\frac{1}{2}\left(-x^3+2x^2\right)}_{\mathsf{f}(x)} + \underbrace{\frac{2}{3}-\frac{1}{2}\left(x^3-x^2-x+1\right)}_{\mathsf{f}(x-1)} + \underbrace{\frac{1}{6}x^3}_{\mathsf{f}(x-2)} \\
&= x^3\left(-\frac{1}{6}+\frac{1}{2}-\frac{1}{2}+\frac{1}{6}\right) + x^2\left(\frac{3}{6}-\frac{2}{2}+\frac{1}{2}\right) + t\left(-\frac{3}{6}+\frac{1}{2}\right) + \left(\frac{1}{6}+\frac{2}{3}+\frac{2}{3}-\frac{1}{2}\right) \\
&= 1
\end{aligned}$$

✎

# 6 Partition of unity

## 6.1 Motivation

A very common property of scaling functions (Definition 3.1 page 21) is the *partition of unity* property (Definition 6.2 page 58). The partition of unity is a kind of generalization of *orthonormality*; that is, *all* orthonormal scaling functions form a partition of unity. But the partition of unity property is not just a consequence of orthonormality, but also a generalization of orthonormality, in that if you remove the orthonormality constraint, the partition of unity is still a reasonable constraint in and of itself.

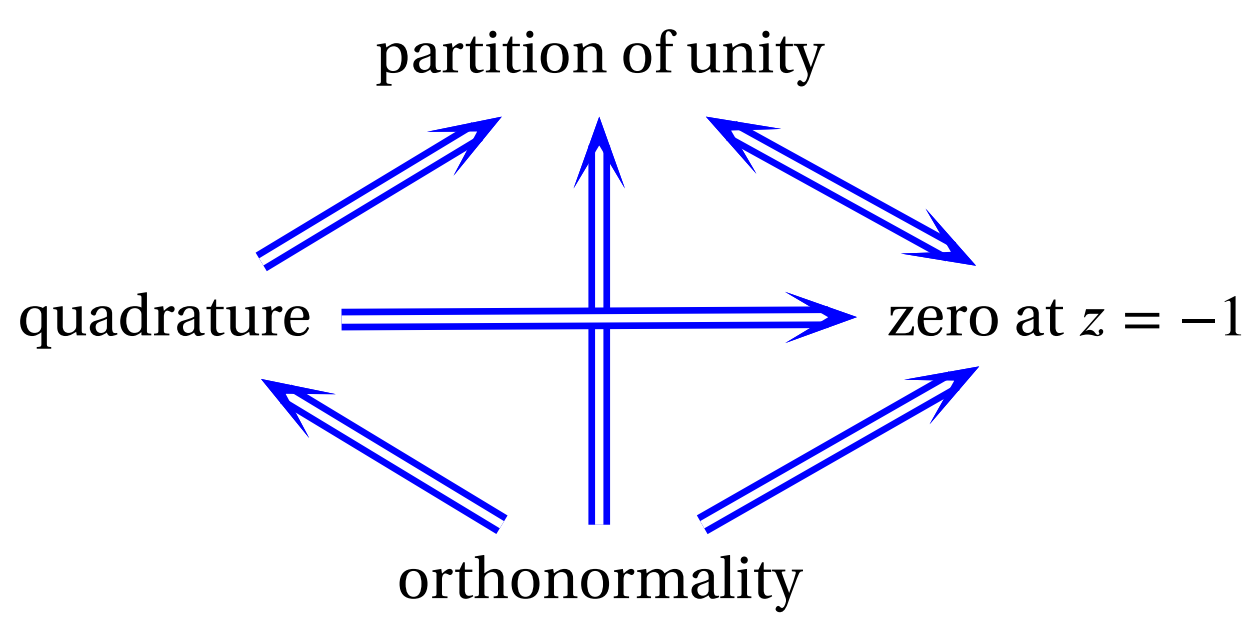

Figure 3: Implications of scaling function properties

There are two reasons why the partition of unity property is a reasonable constraint on its own:

- Without a partition of unity, it is difficult to represent a function as simple as a constant.[92]
- For a multiresolution system $\left(\boldsymbol{L}^2_{\mathbb{R}}, \left(\!\left(\boldsymbol{V}_j\right)\!\right), \phi, \left(\!\left(h_n\right)\!\right)\right)$, the partition of unity property is equivalent to $\sum_{n\in\mathbb{Z}}(-1)^n h_n = 0$ (Theorem 6.8 page 61). As viewed from the perspective of discrete time signal processing, this implies that the scaling coefficients form a "*lowpass filter*"; lowpass filters provide a kind of "coarse approximation" of a function. And that is what the scaling function is "supposed" to do—to provide a coarse approximation at some resolution or "scale" (Definition 3.1 page 21).

## 6.2 Definition and results

**Definition 6.1** The **Kronecker delta function** $\bar{\delta}_n$ is defined as

$$\bar{\delta}_n \triangleq \begin{cases} 1 & \text{for } n=0 \\ 0 & \text{for } n\neq 0 \end{cases} \quad \text{and} \qquad \forall n\in\mathbb{Z}$$

**Definition 6.2** [93] A function $\mathsf{f} \in \mathbb{R}^{\mathbb{R}}$ forms a **partition of unity** if

$$\sum_{n\in\mathbb{Z}} \mathbf{T}^n \mathsf{f}(x) = 1 \qquad \forall x\in\mathbb{R}.$$

**Theorem 6.3** [94] *Let* $\left(\boldsymbol{L}^2_{\mathbb{R}}, \left(\!\left(\boldsymbol{V}_j\right)\!\right), \phi, \left(\!\left(h_n\right)\!\right)\right)$ *be a multiresolution system (Definition 3.6 page 23). Let* $\tilde{\mathbf{F}}\mathsf{f}(\omega)$ *be the* FOURIER TRANSFORM *(Definition 1.26 page 8) of a function* $\mathsf{f} \in \boldsymbol{L}^2_{\mathbb{R}}$*. Let* $\bar{\delta}_n$ *be the*

[92] [65], page 8
[93] [71], page 171, [89], page 225, [63], page 116, [114], page 152, ⟨item 20C⟩, [115], page 152, ⟨item 20C⟩
[94] [65], page 8

KRONECKER DELTA FUNCTION.

$$\underbrace{\sum_{n\in\mathbb{Z}} \mathbf{T}^n \mathsf{f} = c}_{\text{PARTITION OF UNITY in "time"}} \iff \underbrace{\left[\tilde{\mathbf{F}}\mathsf{f}\right](2\pi n) = \bar{\delta}_n}_{\text{PARTITION OF UNITY in "frequency"}}$$

✎PROOF: Let $\mathbb{Z}_e$ be the set of even integers and $\mathbb{Z}_o$ the set of odd integers.

(1) Proof for ( $\implies$ ) case:

$$\begin{aligned}
c &= \sum_{m\in\mathbb{Z}} \mathbf{T}^m \mathsf{f}(x) && \text{by left hypothesis}\\
&= \sum_{m\in\mathbb{Z}} \mathsf{f}(x-m) && \text{by definition of } \mathbf{T} \text{ (Definition 2.1 page 14)}\\
&= \sqrt{2\pi}\sum_{m\in\mathbb{Z}} \tilde{\mathsf{f}}(2\pi m)e^{i2\pi mx} && \text{by } \textit{PSF} \text{ (Theorem 2.22 page 19)}\\
&= \underbrace{\sqrt{2\pi}\,\tilde{\mathsf{f}}(2\pi n)e^{i2\pi nx}}_{\text{real and constant for } n=0} + \underbrace{\sqrt{2\pi}\sum_{m\in\mathbb{Z}\setminus n} \tilde{\mathsf{f}}(2\pi m)e^{i2\pi mx}}_{\text{complex and non-constant}}\\
\implies \sqrt{2\pi}\tilde{\mathsf{f}}(2\pi n) &= c\bar{\delta}_n && \text{because } c \text{ is real and constant for all } t
\end{aligned}$$

(2) Proof for ( $\impliedby$ ) case:

$$\begin{aligned}
\sum_{n\in\mathbb{Z}} \mathbf{T}^n \mathsf{f}(x) &= \sum_{n\in\mathbb{Z}} \mathsf{f}(x-n) && \text{by definition of } \mathbf{T} \text{ (Definition 2.1 page 14)}\\
&= \sqrt{2\pi}\sum_{n\in\mathbb{Z}} \tilde{\mathsf{f}}(2\pi n)e^{-i2\pi nx} && \text{by } \textit{PSF} \text{ (Theorem 2.22 page 19)}\\
&= \sqrt{2\pi}\sum_{n\in\mathbb{Z}} \frac{c}{\sqrt{2\pi}}\bar{\delta}_n e^{-i2\pi nx} && \text{by right hypothesis}\\
&= \sqrt{2\pi}\frac{c}{\sqrt{2\pi}} e^{-i2\pi 0x} && \text{by definition of } \bar{\delta}_n \text{ (Definition 6.1 page 58)}\\
&= c
\end{aligned}$$

✎

**Corollary 6.4**

$$\left\{\begin{array}{l} \exists \mathsf{g} \in \boldsymbol{L}^2_{\mathbb{R}} \textit{ such that} \\ \mathsf{f}(x) = \mathbb{1}_{[-1,1)}(x) \star \mathsf{g}(x) \end{array}\right\} \implies \left\{\begin{array}{l} \mathsf{f}(x) \textit{ generates} \\ \textit{a } \text{PARTITION OF UNITY} \end{array}\right\}$$

**Example 6.5** All *B-splines* form a partition of unity. All B-splines of order $n = 1$ or greater can be generated by convolution with a *pulse* function, similar to that specified in Corollary 6.4 (page 59).

**Example 6.6** Let a function $\mathsf{f}$ be defined in terms of the cosine function (Definition 1.5 page 4) as follows:

$$\mathsf{f}(x) \triangleq \begin{cases} \cos^2\left(\frac{\pi}{2}x\right) & \text{for } |x| \le 1 \\ 0 & \text{otherwise} \end{cases}$$

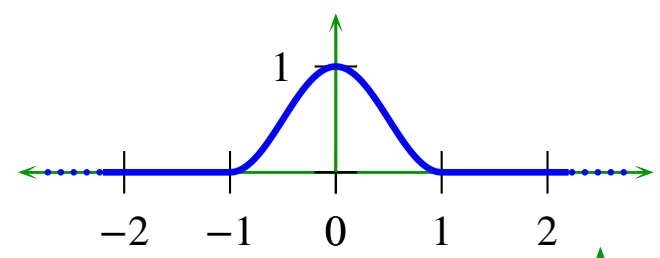

Then $\mathsf{f}$ forms a *partition of unity*:

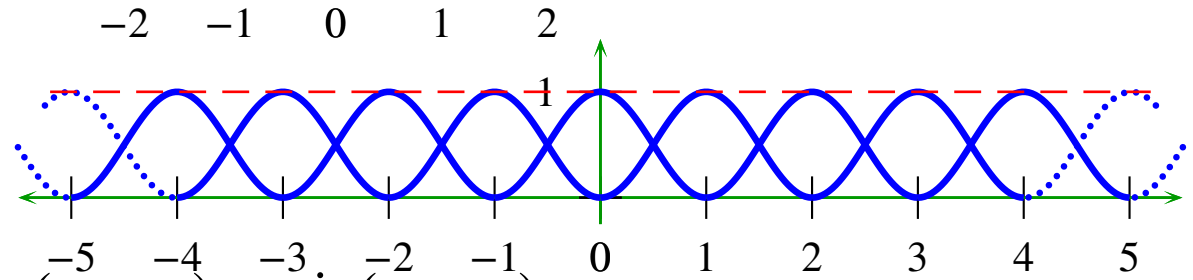

Note that $$\tilde{\mathsf{f}}(\omega) = \frac{1}{2\sqrt{2\pi}}\left[\underbrace{\frac{2\sin\omega}{\omega}}_{2\,\mathrm{sinc}(\omega)} + \underbrace{\frac{\sin(\omega-\pi)}{(\omega-\pi)}}_{\mathrm{sinc}(\omega-\pi)} + \underbrace{\frac{\sin(\omega+\pi)}{(\omega+\pi)}}_{\mathrm{sinc}(\omega+\pi)}\right]$$

and so $\tilde{\mathsf{f}}(2\pi n) = \frac{1}{\sqrt{2\pi}}\bar{\delta}_n$:

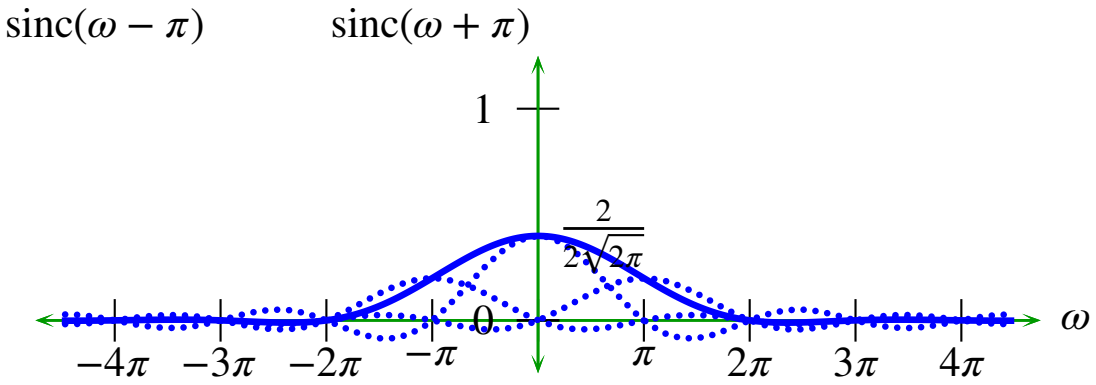

**Example 6.7** (raised cosine) [95] Let a function $\mathsf{f}$ be defined in terms of the cosine function (Definition 1.5 page 4) as follows:

$$\text{Let } \mathsf{f}(x) \triangleq \begin{cases} 1 & \text{for } \quad 0 \le |x| < \frac{1-\beta}{2} \\ \frac{1}{2}\left\{1 + \cos\left[\frac{\pi}{\beta}\left(|x| - \frac{1-\beta}{2}\right)\right]\right\} & \text{for } \quad \frac{1-\beta}{2} \le |x| < \frac{1+\beta}{2} \\ 0 & \text{otherwise} \end{cases}$$

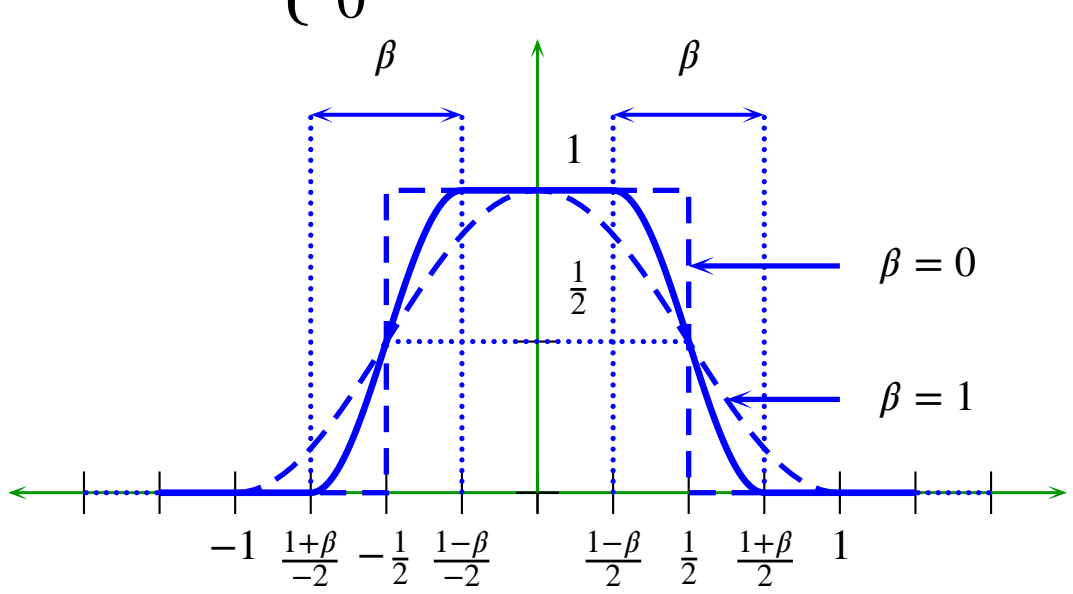

Then $\mathsf{f}$ forms a *partition of unity*:

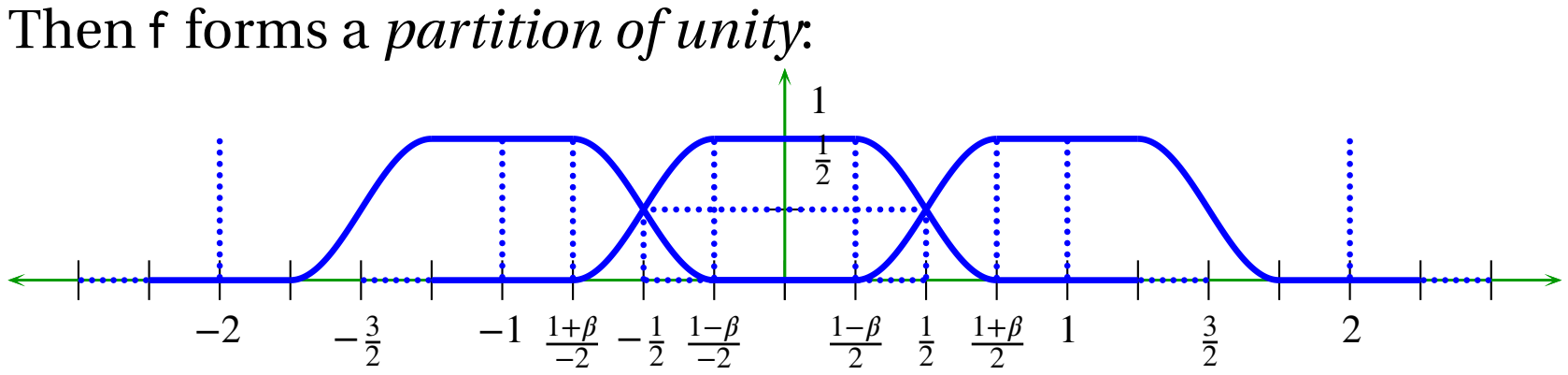

95 [94], pages 560–561

## 6.3 Scaling functions with partition of unity

The $Z$ transform (Definition 1.39 page 11) of a sequence $(\!(h_n)\!)$ with sum $\sum_{n\in\mathbb{Z}}(-1)^n h_n = 0$ has a zero at $z=-1$. Somewhat surprisingly, the *partition of unity* and *zero at* $z=-1$ properties are actually equivalent (next theorem).

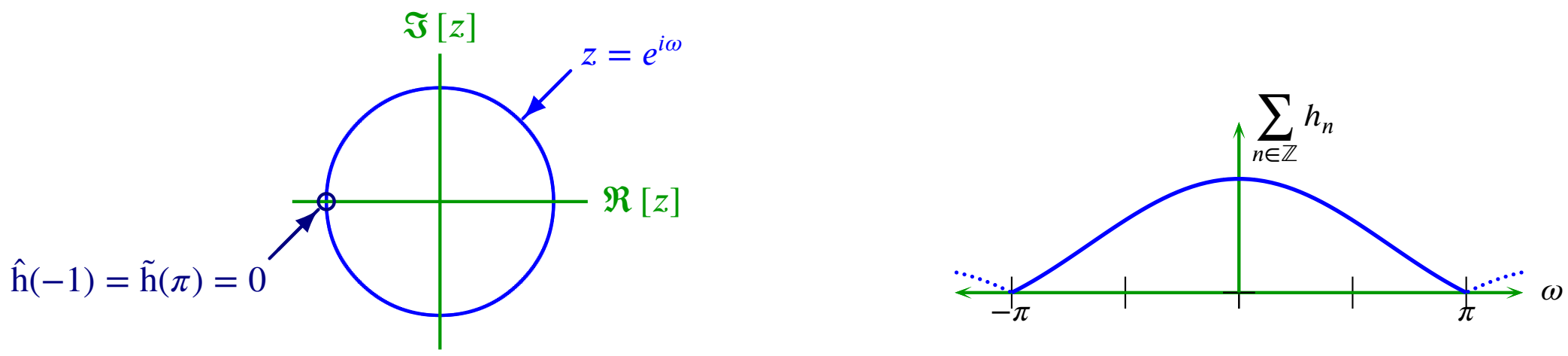

**Theorem 6.8** [96] *Let* $\left(L^2_{\mathbb{R}}, (\!(V_j)\!), \phi, (\!(h_n)\!)\right)$ *be a multiresolution system* (Definition 3.6 page 23). *Let* $\tilde{\mathbf{F}}\mathsf{f}(\omega)$ *be the* FOURIER TRANSFORM (Definition 1.26 page 8) *of a function* $\mathsf{f} \in L^2_{\mathbb{R}}$. *Let* $\bar{\delta}_n$ *be the* KRONECKER DELTA FUNCTION.

$$\underbrace{\sum_{n\in\mathbb{Z}} \mathbf{T}^n\phi = c \quad \text{for some } c\in\mathbb{R}\backslash 0}_{(1)\text{ PARTITION OF UNITY}} \iff \underbrace{\sum_{n\in\mathbb{Z}}(-1)^n h_n = 0}_{(2)\text{ ZERO AT } z=-1} \iff \underbrace{\sum_{n\in\mathbb{Z}} h_{2n} = \sum_{n\in\mathbb{Z}} h_{2n+1} = \frac{\sqrt{2}}{2}}_{(3)\text{ sum of even = sum of odd = } \frac{\sqrt{2}}{2}}$$

✎PROOF: Let $\mathbb{Z}_e$ be the set of even integers and $\mathbb{Z}_o$ the set of odd integers.

(1) Proof that (1) $\Longleftarrow$ (2):

$$\begin{aligned}
\sum_{n\in\mathbb{Z}} \mathbf{T}^n\phi &= \sum_{n\in\mathbb{Z}} \mathbf{T}^n \left[\sum_{m\in\mathbb{Z}} h_m \mathbf{D}\mathbf{T}^m\phi\right] && \text{by dilation equ. (Theorem 3.4 page 23)}\\
&= \sum_{m\in\mathbb{Z}} h_m \sum_{n\in\mathbb{Z}} \mathbf{T}^n\mathbf{D}\mathbf{T}^m\phi\\
&= \sum_{m\in\mathbb{Z}} h_m \sum_{n\in\mathbb{Z}} \mathbf{D}\mathbf{T}^{2n}\mathbf{T}^m\phi && \text{by Proposition 2.9 page 16}\\
&= \mathbf{D}\sum_{m\in\mathbb{Z}} h_m \sum_{n\in\mathbb{Z}} \mathbf{T}^{2n}\mathbf{T}^m\phi\\
&= \mathbf{D}\sum_{m\in\mathbb{Z}} h_m \left[\sqrt{\frac{2\pi}{2}}\hat{\mathbf{F}}^{-1}\mathbf{S}_2\tilde{\mathbf{F}}(\mathbf{T}^m\phi)\right] && \text{by PSF (Theorem 2.22 page 19)}\\
&= \sqrt{\pi}\mathbf{D}\sum_{m\in\mathbb{Z}} h_m \hat{\mathbf{F}}^{-1}\mathbf{S}_2 e^{-i\omega m}\tilde{\mathbf{F}}\phi && \text{by Corollary 2.19 page 18}\\
&= \sqrt{\pi}\mathbf{D}\sum_{m\in\mathbb{Z}} h_m \hat{\mathbf{F}}^{-1} e^{-i\frac{2\pi}{2}km}\mathbf{S}_2\tilde{\mathbf{F}}\phi && \text{by definition of } \mathbf{S} \text{ (Theorem 2.22 page 19)}
\end{aligned}$$

[96] [65], page 8, [25], page 123

$$
\begin{aligned}
&= \sqrt{\pi}\mathbf{D}\sum_{m\in\mathbb{Z}} h_m \hat{\mathbf{F}}^{-1}(-1)^{km}\mathbf{S}_2\tilde{\mathbf{F}}\phi \\
&= \sqrt{\pi}\mathbf{D}\sum_{m\in\mathbb{Z}} h_m \left[\tfrac{\sqrt{2}}{2}\sum_{k\in\mathbb{Z}}(-1)^{km}\left(\mathbf{S}_2\tilde{\mathbf{F}}\phi\right)e^{i\frac{2\pi}{2}kx}\right] && \text{by def. of } \hat{\mathbf{F}}^{-1} \text{ (Theorem 1.21 page 6)} \\
&= \frac{\sqrt{2\pi}}{2}\mathbf{D}\sum_{k\in\mathbb{Z}}\left(\mathbf{S}_2\tilde{\mathbf{F}}\phi\right)e^{i\pi kx}\sum_{m\in\mathbb{Z}}(-1)^{km}h_m \\
&= \frac{\sqrt{2\pi}}{2}\mathbf{D}\sum_{k\in\mathbb{Z}_e}\left(\mathbf{S}_2\tilde{\mathbf{F}}\phi\right)e^{i\pi kx}\sum_{m\in\mathbb{Z}}(-1)^{km}h_m \\
&\quad + \frac{\sqrt{2\pi}}{2}\mathbf{D}\sum_{k\in\mathbb{Z}_o}\left(\mathbf{S}_2\tilde{\mathbf{F}}\phi\right)e^{i\pi kx}\sum_{m\in\mathbb{Z}}(-1)^{km}h_m \\
&= \frac{\sqrt{2\pi}}{2}\mathbf{D}\sum_{k\in\mathbb{Z}_e}\left(\mathbf{S}_2\tilde{\mathbf{F}}\phi\right)e^{i\pi kx}\underbrace{\sum_{m\in\mathbb{Z}}h_m}_{\sqrt{2}} \\
&\quad + \frac{\sqrt{2\pi}}{2}\mathbf{D}\sum_{k\in\mathbb{Z}_o}\left(\mathbf{S}_2\tilde{\mathbf{F}}\phi\right)e^{i\pi kx}\underbrace{\sum_{m\in\mathbb{Z}}(-1)^m h_m}_{0} \\
&= \sqrt{\pi}\mathbf{D}\sum_{k\in\mathbb{Z}_e}\left(\mathbf{S}_2\tilde{\mathbf{F}}\phi\right)e^{i\pi kx} && \text{by Theorem 3.9 (page 24) and right hyp.} \\
&= \sqrt{\pi}\mathbf{D}\sum_{k\in\mathbb{Z}_e}\tilde{\phi}\left(\frac{2\pi}{2}k\right)e^{i\pi kx} && \text{by definitions of } \tilde{\mathbf{F}} \text{ and } \mathbf{S}_2 \\
&= \sqrt{\pi}\mathbf{D}\sum_{k\in\mathbb{Z}}\tilde{\phi}(2\pi k)e^{i2\pi kx} && \text{by definition of } \mathbb{Z}_e \\
&= \frac{1}{\sqrt{2}}\mathbf{D}\left\{\sqrt{2\pi}\sum_{k\in\mathbb{Z}}\tilde{\phi}(2\pi k)e^{i2\pi kx}\right\} \\
&= \frac{1}{\sqrt{2}}\mathbf{D}\sum_{n\in\mathbb{Z}}\phi(x+n) && \text{by PSF (Theorem 2.22 page 19)} \\
&= \frac{1}{\sqrt{2}}\mathbf{D}\sum_{n}\mathbf{T}^n\phi && \text{by definition of } \mathbf{T} \text{ (Definition 2.1 page 14)}
\end{aligned}
$$

The above equation sequence demonstrates that

$$\mathbf{D}\sum_n \mathbf{T}^n\phi = \sqrt{2}\sum_n \mathbf{T}^n\phi$$

(essentially that $\sum_n \mathbf{T}^n\phi$ is equal to it's own dilation). This implies that $\sum_n \mathbf{T}^n\phi$ is a constant (Proposition 2.12 page 17).

(2) Proof that (1) $\implies$ (2):

$$
\begin{aligned}
c &= \sum_{n\in\mathbb{Z}} \mathbf{T}^n \phi && \text{by left hypothesis} \\
&= \sqrt{2\pi}\, \hat{\mathbf{F}}^{-1} \mathbf{S} \tilde{\mathbf{F}} \phi && \text{by PSF (Theorem 2.22 page 19)} \\
&= \sqrt{2\pi}\, \hat{\mathbf{F}}^{-1} \mathbf{S} \underbrace{\sqrt{2}\left( \mathbf{D}^{-1} \sum_{n\in\mathbb{Z}} h_n e^{-i\omega n} \right) \left( \mathbf{D}^{-1} \tilde{\mathbf{F}} \phi \right)}_{\tilde{\mathbf{F}}\phi} && \text{by Lemma 3.5 page 23} \\
&= 2\sqrt{\pi}\, \hat{\mathbf{F}}^{-1} \left( \mathbf{S}\mathbf{D}^{-1} \sum_{n\in\mathbb{Z}} h_n e^{-i\omega n} \right) \left( \mathbf{S} \tilde{\mathbf{F}} \mathbf{D} \phi \right) && \text{by Corollary 2.19 page 18} \\
&= 2\sqrt{\pi}\, \hat{\mathbf{F}}^{-1} \left( \mathbf{S} \frac{1}{\sqrt{2}} \sum_{n\in\mathbb{Z}} h_n e^{-i\frac{\omega}{2} n} \right) \left( \mathbf{S} \tilde{\mathbf{F}} \mathbf{D} \phi \right) && \text{by Proposition 2.4 page 15} \\
&= \sqrt{2\pi}\, \hat{\mathbf{F}}^{-1} \left( \sum_{n\in\mathbb{Z}} h_n e^{-i\frac{2\pi k}{2} n} \right) \left( \mathbf{S} \tilde{\mathbf{F}} \mathbf{D} \phi \right) && \text{by def. of } \mathbf{S} \text{ (Theorem 2.22 page 19)} \\
&= \sqrt{2\pi}\, \hat{\mathbf{F}}^{-1} \left( \sum_{n\in\mathbb{Z}} h_n (-1)^{kn} \right) \left( \mathbf{S} \mathbf{D}^{-1} \mathbf{F} \phi \right) \\
&= \sqrt{2\pi}\, \hat{\mathbf{F}}^{-1} \left( \sum_{n\in\mathbb{Z}} h_n (-1)^{kn} \right) \left( \mathbf{S} \frac{1}{\sqrt{2}} \tilde{\phi}\left(\frac{\omega}{2}\right) \right) && \text{by def. of } \mathbf{S} \text{ (Theorem 2.22 page 19)} \\
&= \sqrt{2\pi}\, \hat{\mathbf{F}}^{-1} \left( \sum_{n\in\mathbb{Z}} h_n (-1)^{kn} \right) \left( \frac{1}{\sqrt{2}} \tilde{\phi}\left(\frac{2\pi k}{2}\right) \right) \\
&= \sqrt{\pi} \sum_{k\in\mathbb{Z}} \sum_{n\in\mathbb{Z}} h_n (-1)^{kn}\, \tilde{\phi}(\pi k) e^{i2\pi kx} && \text{by Theorem 1.21 page 6} \\
&= \sqrt{\pi} \sum_{k \text{ even}} \sum_{n\in\mathbb{Z}} h_n (-1)^{kn}\, \tilde{\phi}(\pi k) e^{i2\pi kx} \\
&\quad + \sqrt{\pi} \sum_{k \text{ odd}} \sum_{n\in\mathbb{Z}} h_n (-1)^{kn}\, \tilde{\phi}(\pi k) e^{i2\pi kx} \\
&= \sqrt{\pi} \sum_{k \text{ even}} \underbrace{\left( \sum_{n\in\mathbb{Z}} h_n \right)}_{\sqrt{2}} \tilde{\phi}(\pi k) e^{i2\pi kx} \\
&\quad + \sqrt{\pi} \sum_{k \text{ odd}} \left( \sum_{n\in\mathbb{Z}} h_n (-1)^n \right) \tilde{\phi}(\pi k) e^{i2\pi kx} \\
&= \sqrt{\pi} \sum_{k\in\mathbb{Z}} \sqrt{2}\, \tilde{\phi}(\pi 2k) e^{i2\pi 2kx}
\end{aligned}
$$

$$+\sqrt{\pi}\sum_{k\in\mathbb{Z}}\left(\sum_{n\in\mathbb{Z}} h_n(-1)^n\right)\tilde{\phi}(\pi[2k+1])e^{i2\pi[2k+1]x}$$ by Theorem 3.9 page 24

$$=\frac{\sqrt{2\pi}}{\sqrt{2\pi}}\tilde{\phi}(0)+\sqrt{\pi}e^{i2\pi x}\sum_{n\in\mathbb{Z}} h_n(-1)^n\sum_{k\in\mathbb{Z}}\tilde{\phi}(\pi[2k+1])e^{i4\pi kx}$$ by left hyp. and Theorem 6.3 page 58

$$\implies \left(\sum_{n\in\mathbb{Z}} h_n(-1)^n\right)=0$$ because the right side must equal $c$

(3) Proof that (2) $\implies$ (3):

$$\sum_{n\in\mathbb{Z}_e} h_n=\sum_{n\in\mathbb{Z}_o} h_n=\frac{1}{2}\sum_{n\in\mathbb{Z}} h_n$$ by (2) and Proposition 1.45 page 13

$$=\frac{\sqrt{2}}{2}$$ by *admissibility condition* (Theorem 3.9 page 24)

(4) Proof that (2) $\impliedby$ (3):

$$\frac{\sqrt{2}}{2}=\underbrace{\sum_{n\in\mathbb{Z}_e}(-1)^n h_n}_{\text{even terms}}+\underbrace{\sum_{n\in\mathbb{Z}_o}(-1)^n h_n}_{\text{odd terms}}$$ by (3)

$$\implies \sum_{n\in\mathbb{Z}}(-1)^n h_n=0$$ by Proposition 1.45 page 13

✏

**Proposition 6.9**
$\phi(x)$ *generates a* PARTITION OF UNITY $\not\Longrightarrow$ $\phi(x)$ *generates an* MRA *system.*

## 6.4 Spline wavelet systems

**Theorem 6.10** [97] *Let* $\boldsymbol{S}^n(\mathbb{Z})$ *be the* SPACE OF ALL SPLINES OF ORDER N *(Definition 5.11 page 54). For each* $n\in\mathbb{W}$,
$\boldsymbol{S}^n(2^k\mathbb{Z})$ *is a* MULTIRESOLUTION ANALYSIS *(an* MRA*).*

**Theorem 6.11** (B-spline wavelet coefficients) *Let* $\left(\boldsymbol{L}^2_{\mathbb{R}}, (\!(\boldsymbol{V}_j)\!), \phi, (\!(h_n)\!)\right)$ *be an* MRA SYSTEM *(Definition 3.6 page 23). Let* $\mathsf{N}_n(x)$ *be a* $n$TH ORDER B-SPLINE.

$$\underbrace{\phi(x)\triangleq \mathsf{N}_n(x)}_{\text{(1) B-spline scaling function}} \implies (\!(h_k)\!)=\begin{cases}\frac{\sqrt{2}}{2^{n+1}}\binom{n}{k} & \text{for } k=0,1,\ldots,n\\ 0 & \text{otherwise}\end{cases}$$ *(2) scaling sequence in "time"*

$$\iff \left.\breve{\mathsf{h}}(z)\right|_{z\triangleq e^{i\omega}}=\left.\frac{\sqrt{2}}{2^n}\left(1+z^{-1}\right)^{n+1}\right|_{z\triangleq e^{i\omega}}$$ *(3) scaling sequence in "z domain"*

$$\iff \tilde{\mathsf{h}}(\omega)=2\sqrt{2}e^{-i\frac{n+1}{2}\omega}\left[\cos\left(\frac{\omega}{2}\right)\right]^{n+1}$$ *(4) scaling sequence in "frequency"*

[97] 📚 [116], page 57, ⟨Theorem 3.13⟩

✎PROOF:

(1) Proof that (1) $\Longrightarrow$ (3): By Theorem 6.10 page 64 we know that $\mathsf{N}_n(x)$ is a *scaling function* (Definition 3.1 page 21). So then we know that we can use Lemma 3.5 page 23.

$$\begin{aligned}
\tilde{\mathrm{h}}(\omega) &= \sqrt{2}\frac{\tilde{\phi}(2\omega)}{\tilde{\phi}(\omega)} && \text{by Lemma 3.5 page 23}\\
&= \sqrt{2}\frac{\tilde{\mathsf{N}}_n(2\omega)}{\tilde{\mathsf{N}}_n(\omega)} && \text{by (1)}\\
&= \sqrt{2}\frac{\frac{1}{\sqrt{2\pi}}\left(\frac{1-e^{-i2\omega}}{2i\omega}\right)^{n+1}}{\frac{1}{\sqrt{2\pi}}\left(\frac{1-e^{-i\omega}}{i\omega}\right)^{n+1}} && \text{by Theorem 5.10 page 53}\\
&= \frac{\sqrt{2}}{2^{n+1}}\left(\frac{1-z^{-2}}{1-z^{-1}}\right)^{n+1}\Bigg|_{z=e^{i\omega}}\\
&= \frac{\sqrt{2}}{2^{n+1}}\left[\left(\frac{1-z^{-2}}{1-z^{-1}}\right)\left(\frac{1+z^{-1}}{1+z^{-1}}\right)\right]^{n+1}\Bigg|_{z=e^{i\omega}}\\
&= \frac{\sqrt{2}}{2^{n+1}}\left(\frac{\left(1-z^{-2}\right)\left(1+z^{-1}\right)}{1-z^{-2}}\right)^{n+1}\Bigg|_{z=e^{i\omega}}\\
&= \frac{\sqrt{2}}{2^{n}}\left(1+z^{-1}\right)^{n+1}\Big|_{z\triangleq e^{i\omega}}
\end{aligned}$$

(2) Proof that (3) $\Longleftrightarrow$ (2):

$$\begin{aligned}
\hat{\mathrm{h}}(z)\Big|_{z\triangleq e^{i\omega}} &= \frac{\sqrt{2}}{2^{n}}\left(1+z^{-1}\right)^{n+1}\Bigg|_{z\triangleq e^{i\omega}} && \text{by (3)}\\
&= \frac{\sqrt{2}}{2^{n}}\left(\sum_{k=0}^{n+1}\binom{n}{k}z^{-k}\right)\Bigg|_{z\triangleq e^{i\omega}} && \text{by \textit{binomial theorem}}\\
\Longleftrightarrow\qquad h_k &= \frac{\sqrt{2}}{2^{n+1}}\binom{n}{k} && \text{by definition of \textit{Z transform} (Definition 1.39 page 11)}
\end{aligned}$$

(3) Proof that (3) $\Longrightarrow$ (4):

$$\begin{aligned}
\tilde{\mathrm{h}}(\omega) &= \hat{\mathrm{h}}(z)\Big|_{z\triangleq e^{i\omega}} && \text{by definition of \textit{DTFT} (Definition 1.42 page 12)}\\
&= \frac{\sqrt{2}}{2^{n}}\left(1+z^{-1}\right)^{n+1}\Bigg|_{z\triangleq e^{i\omega}} && \text{by (3)}\\
&= \frac{\sqrt{2}}{2^{n}}\left(1+e^{-i\omega}\right)^{n+1} && \text{by definition of } z
\end{aligned}$$

$$= \frac{\sqrt{2}}{2^n}\left[e^{-i\frac{1}{2}\omega}\left(e^{i\frac{\omega}{2}} + e^{-i\frac{\omega}{2}}\right)\right]^{n+1}$$
$$= \frac{\sqrt{2}}{2^n}e^{-i\frac{n+1}{2}\omega}\left[2\cos\left(\frac{\omega}{2}\right)\right]^{n+1}$$
$$= 2\sqrt{2}e^{-i\frac{n+1}{2}\omega}\left[\cos\left(\frac{\omega}{2}\right)\right]^{n+1}$$

(4) Proof that (3) $\Longleftarrow$ (4):

$$\begin{aligned}
\hat{\mathrm{h}}(z)\Big|_{z\triangleq e^{i\omega}} &= \hat{\mathrm{h}}\left(e^{i\omega}\right) \\
&= \tilde{\mathrm{h}}(\omega) \\
&= 2\sqrt{2}e^{-i\frac{n+1}{2}\omega}\left[\cos\left(\frac{\omega}{2}\right)\right]^{n+1} && \text{by (4)} \\
&= \frac{\sqrt{2}}{2^n}e^{-i\frac{n+1}{2}\omega}\left[2\cos\left(\frac{\omega}{2}\right)\right]^{n+1} \\
&= \frac{\sqrt{2}}{2^n}\left[e^{-i\frac{1}{2}\omega}\left(e^{i\frac{\omega}{2}} + e^{-i\frac{\omega}{2}}\right)\right]^{n+1} \\
&= \frac{\sqrt{2}}{2^n}\left(1 + e^{-i\omega}\right)^{n+1} \\
&= \frac{\sqrt{2}}{2^n}\left(1 + z^{-1}\right)^{n+1}\Big|_{z\triangleq e^{i\omega}}
\end{aligned}$$

✏

## 6.5 Examples

**Example 6.12** (2 coefficient case/Haar wavelet system/order 0 B-spline wavelet system) [98]

Let $\left(\boldsymbol{L}^2_{\mathbb{R}}, (\!(\boldsymbol{V}_j)\!), (\!(\boldsymbol{W}_j)\!), \phi, \psi, (\!(h_n)\!), (\!(g_n)\!)\right)$ be an *orthogonal* wavelet system with two non-zero scaling coefficients.

$$\left\{\begin{array}{lll}
1. & \operatorname{supp}\phi(x) = [0,\ 1] & \text{(Theorem 3.20 page 28) and} \\
2. & \textit{admissibility condition} & \text{(Theorem 3.9 page 24) and} \\
3. & \textit{partition of unity} & \text{(Theorem 6.8 page 61) and} \\
4. & \mathrm{g}_n = \pm(-1)^n h^*_{N-n}\ \forall n\in\mathbb{Z} & \text{(Theorem 3.18 page 27)}
\end{array}\right\} \Longrightarrow \left\{\begin{array}{r|c|c}
n & h_n & g_n \\ \hline
0 & \frac{\sqrt{2}}{2} & \frac{\sqrt{2}}{2} \\
1 & \frac{\sqrt{2}}{2} & -\frac{\sqrt{2}}{2} \\
\text{other} & 0 & 0
\end{array}\right\}$$

[98] [54], [116], pages 14–15, ⟨"Sources and comments"⟩

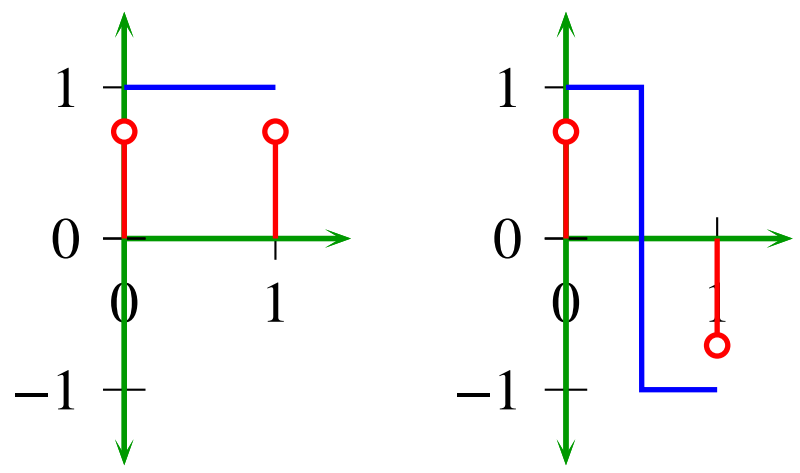

✎PROOF:

(1) Proof that (1) $\Longrightarrow$ that only $h_0$ and $h_1$ are non-zero: by Theorem 3.20 page 28.

(2) Proof for values of $h_0$ and $h_1$:

(a) Method 1: Under the constraint of two non-zero scaling coefficients, a scaling function design is fully constrained using the *admissibility equation* (Theorem 3.9 page 24) and the *partition of unity* constraint (Definition 6.2 page 58). The partition of unity formed by $\phi(x)$ is illustrated in Example 5.14 page 55.
Here are the equations:

$$\begin{array}{lcll} h_0 + h_1 &=& \sqrt{2} & \text{(admissibility equation Theorem 3.9 page 24)} \\ h_0 - h_1 &=& 0 & \text{(partition of unity/zero at } -1 \text{ Theorem 6.8 page 61)} \end{array}$$

Here are the calculations for the coefficients:

$$\begin{aligned} (h_0 + h_1) + (h_0 - h_1) = 2h_0 &= \sqrt{2} && \text{(add two equations together)} \\ (h_0 + h_1) - (h_0 - h_1) = 2h_1 &= \sqrt{2} && \text{(subtract second from first)} \end{aligned}$$

$$\begin{aligned} g_0 &= h_1 \\ g_1 &= -h_0 \end{aligned}$$

(b) Method 2: By Theorem 6.11 page 64.

(3) Note: $h_0$ and $h_1$ can also be produced using other systems of equations including the following:

(a) Admissibility condition and *orthonormality*

(b) *Daubechies-p*1 wavelets computed using spectral techniques

(4) Proof for values of $g_0$ and $g_1$: by (4) and Theorem 3.18 page 27.

✏

**Example 6.13** (order 1 B-spline wavelet system) [99] The following figures illustrate scaling and wavelet coefficients and functions for the *B-Spline* $B_2$, or *tent function*. The partition of unity formed by the scaling function $\phi(x)$ is illustrated in Example 5.15 page 56.

[99] 📚 [104], page 616, 📚 [30], pages 146–148, ⟨§5.4⟩

| $n$ | $h_n$ | $g_n$ |
|---|---|---|
| 0 | $\left(\frac{\sqrt{2}}{4}\right)$ | $\left(\frac{\sqrt{2}}{4}\right)$ |
| 1 | $2\left(\frac{\sqrt{2}}{4}\right)$ | $-2\left(\frac{\sqrt{2}}{4}\right)$ |
| 2 | $\left(\frac{\sqrt{2}}{4}\right)$ | $\left(\frac{\sqrt{2}}{4}\right)$ |

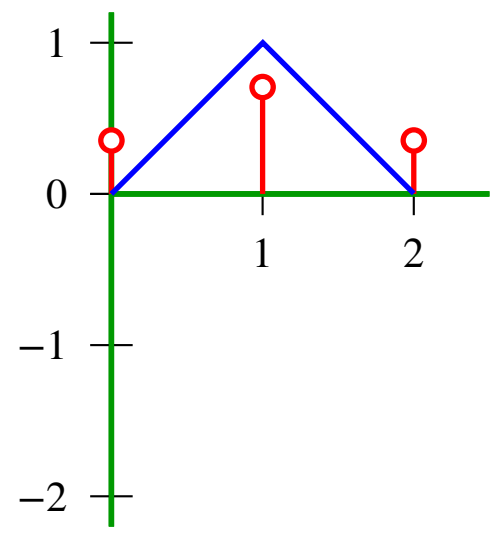

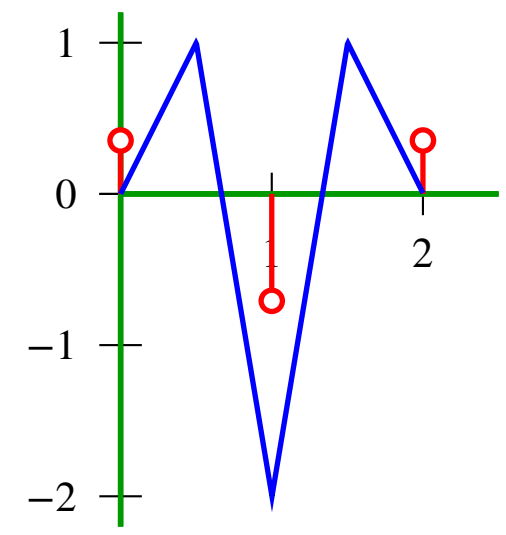

✎PROOF: These results follow from Theorem 6.11 page 64.

$$\left(\begin{array}{ccccc} & & 1 & & \\ & 1 & & 1 & \\ 1 & & 2 & & 1 \end{array}\right)$$

✏

**Example 6.14** (order 3 B-spline wavelet system) [100] The following figures illustrate scaling and wavelet coefficients and functions for a *B-spline*.

| $n$ | $h_n$ | $g_n$ |
|---|---|---|
| 0 | $\left(\frac{\sqrt{2}}{16}\right)$ | $\left(\frac{\sqrt{2}}{16}\right)$ |
| 1 | $4\left(\frac{\sqrt{2}}{16}\right)$ | $-4\left(\frac{\sqrt{2}}{16}\right)$ |
| 2 | $6\left(\frac{\sqrt{2}}{16}\right)$ | $6\left(\frac{\sqrt{2}}{16}\right)$ |
| 3 | $4\left(\frac{\sqrt{2}}{16}\right)$ | $-4\left(\frac{\sqrt{2}}{16}\right)$ |
| 4 | $\left(\frac{\sqrt{2}}{16}\right)$ | $\left(\frac{\sqrt{2}}{16}\right)$ |

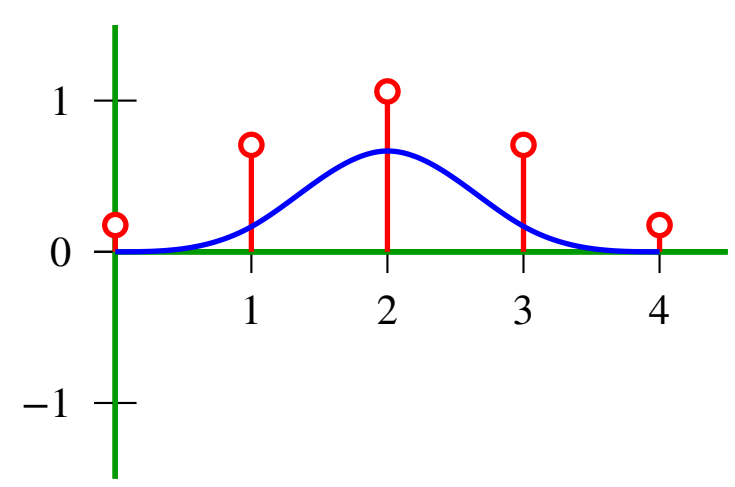

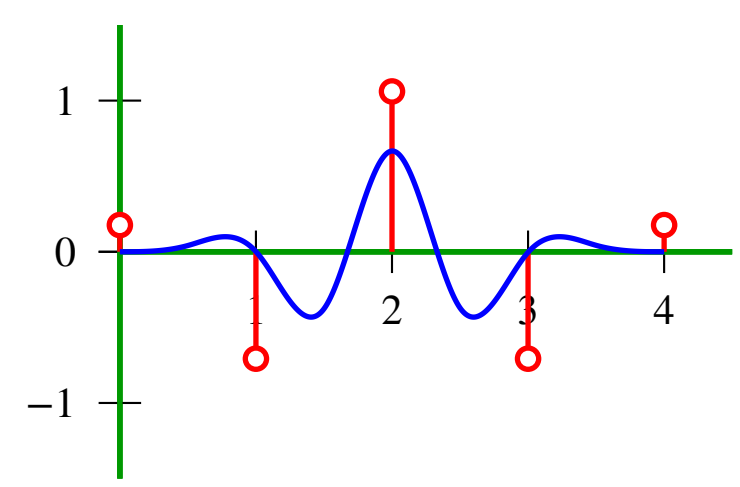

✎PROOF: These results follow from Theorem 6.11 page 64.

$$\left(\begin{array}{ccccccccc} & & & & 1 & & & & \\ & & & 1 & & 1 & & & \\ & & 1 & & 2 & & 1 & & \\ & 1 & & 3 & & 3 & & 1 & \\ 1 & & 4 & & 6 & & 4 & & 1 \end{array}\right)$$

✏

Not all functions that form a *partition of unity* are a bases for an *MRA*. Counterexample 6.15 (next) and Counterexample 6.16 (page 71) provide two counterexamples.

**Counterexample 6.15** Let a function `f` be defined in terms of the sine function (Definition 1.5 page 4) as follows:

[100] [104], page 616

$$\phi(x) \triangleq \begin{cases} \sin^2\left(\frac{\pi}{2}x\right) & \text{for } x \in [0,\ 2] \\ 0 & \text{otherwise} \end{cases}$$

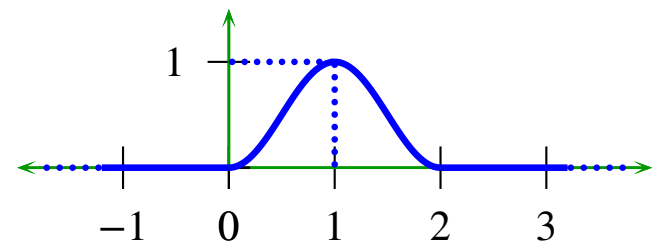

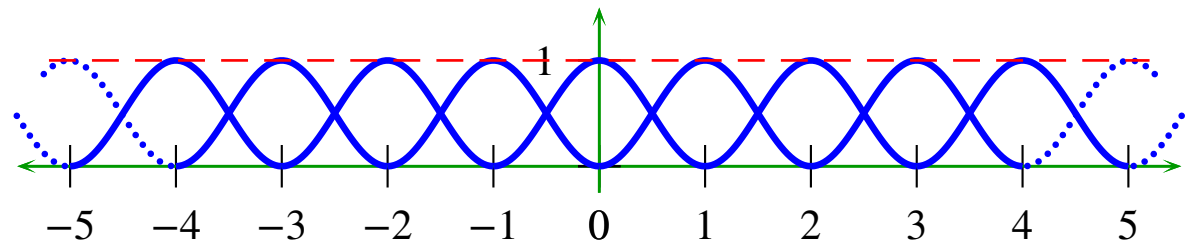

Then $\int_{\mathbb{R}} \phi(x)\, \mathrm{d}x = 1$ and $\phi$ forms a *partition of unity*

**but** $\left\{\mathbf{T}^n\phi|_{n\in\mathbb{Z}}\right\}$ does **not** generate an *MRA*.

✎PROOF: Let $\mathbb{1}_A(x)$ be the *set indicator function* (Definition 1.3 page 3) on a set $A$.

(1) Proof that $\int_{\mathbb{R}} \phi(x)\, \mathrm{d}x = 1$:

$$\begin{aligned}
\int_{\mathbb{R}} \phi(x)\, \mathrm{d}x &= \int_{\mathbb{R}} \sin^2\left(\frac{\pi}{2}x\right)\mathbb{1}_{[0,\ 2]}(x)\, \mathrm{d}x && \text{by definition of } \phi(x) \\
&= \int_0^2 \sin^2\left(\frac{\pi}{2}x\right) \mathrm{d}x && \text{by definition of } \mathbb{1}_A(x) \text{ (Definition 1.3 page 3)} \\
&= \int_0^2 \frac{1}{2}[1 - \cos(\pi x)]\, \mathrm{d}x && \text{by Theorem 1.19 page 6} \\
&= \frac{1}{2}\left[x - \frac{1}{\pi}\sin(\pi x)\right]_0^2 \\
&= \frac{1}{2}[2 - 0 - 0 - 0] \\
&= 1
\end{aligned}$$

(2) Proof that $\phi(x)$ forms a *partition of unity*:

$$\begin{aligned}
\sum_{n\in\mathbb{Z}} \mathbf{T}^n\phi(x) &= \sum_{n\in\mathbb{Z}} \mathbf{T}^n \sin^2\left(\frac{\pi}{2}x\right)\mathbb{1}_{[0,\ 2]}(x) && \text{by definition of } \phi(x) \\
&= \sum_{n\in\mathbb{Z}} \mathbf{T}^n \sin^2\left(\frac{\pi}{2}x\right)\mathbb{1}_{[0,\ 2)}(x) && \text{because } \sin^2\left(\frac{\pi}{2}x\right) = 0 \text{ when } x = 2 \\
&= \sum_{m\in\mathbb{Z}} \mathbf{T}^{m-1} \sin^2\left(\frac{\pi}{2}x\right)\mathbb{1}_{[0,\ 2)}(x) && \text{where } m \triangleq n+1 \implies n = m-1 \\
&= \sum_{m\in\mathbb{Z}} \sin^2\left(\frac{\pi}{2}(x-m+1)\right)\mathbb{1}_{[0,\ 2)}(x-m+1) && \text{by definition of } \mathbf{T} \text{ (Definition 2.1 page 14)} \\
&= \sum_{m\in\mathbb{Z}} \sin^2\left(\frac{\pi}{2}(x-m) + \frac{\pi}{2}\right)\mathbb{1}_{[-1,\ 1)}(x-m) \\
&= \sum_{m\in\mathbb{Z}} \cos^2\left(\frac{\pi}{2}(x-m)\right)\mathbb{1}_{[-1,\ 1)}(x-m) && \text{by Theorem 1.19 page 6} \\
&= \sum_{m\in\mathbb{Z}} \mathbf{T}^m \cos^2\left(\frac{\pi}{2}x\right)\mathbb{1}_{[-1,\ 1)}(x) && \text{by definition of } \mathbf{T} \text{ (Definition 2.1 page 14)}
\end{aligned}$$

$$= \sum_{m\in\mathbb{Z}} \mathbf{T}^m \cos^2\left(\frac{\pi}{2}x\right)\mathbb{1}_{[-1,\,1]}(x) \qquad \text{because } \cos^2\left(\frac{\pi}{2}x\right) = 0 \text{ when } x = 1$$

$$= 1 \qquad \text{by Example 6.6 page 60}$$

(3) Proof that $\phi(x) \notin \operatorname{span}\left\{\mathbf{DT}^n\phi(x)|_{n\in\mathbb{Z}}\right\}$ (and so does not generate an *MRA*):

(a) Note that the *support* (Definition 3.19 page 28) of $\phi$ is $\operatorname{supp}\phi = [0,\,2]$.

(b) Therefore, the *support* of $(\!(h_n)\!)$ is $\operatorname{supp}(\!(h_n)\!) = \{0, 1, 2\}$ (Theorem 3.20 page 28).

(c) So if $\phi(x)$ *is* an *MRA*, we only need to compute $\{h_0, h_1, h_2\}$ (the rest would be 0).

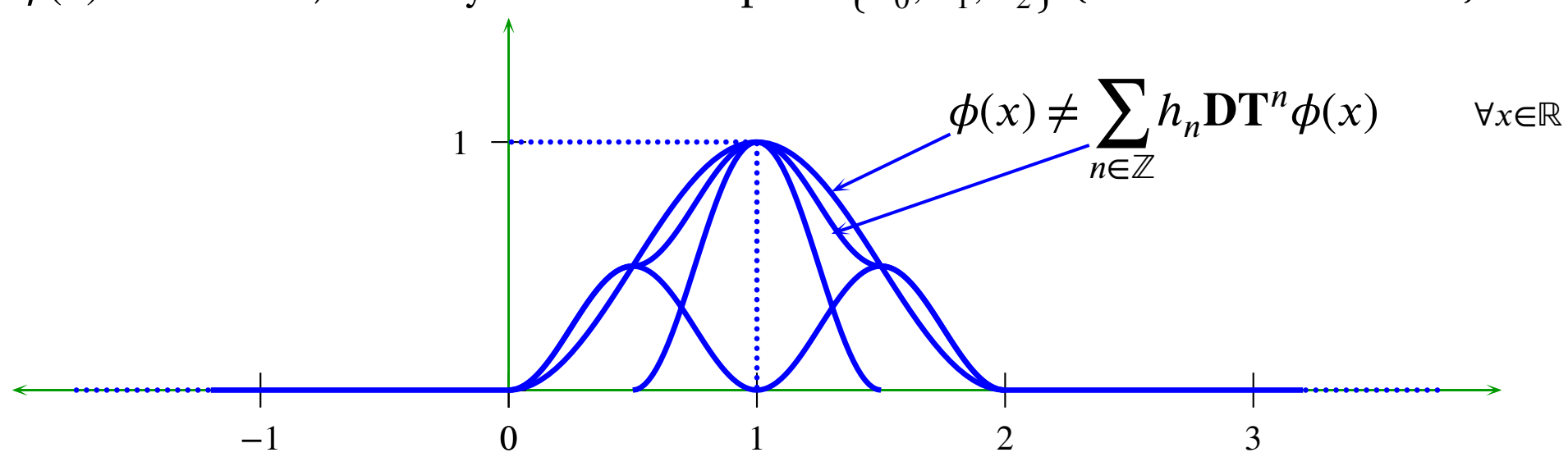

Here would be the values of $\{h_1, h_2, h_3\}$:

$$\begin{aligned}
\phi(x) &= \sum_{n\in\mathbb{Z}} h_n \mathbf{DT}^n \phi(x) \\
&= \sum_{n\in\mathbb{Z}} h_n \mathbf{DT}^n \sin^2\left(\frac{\pi}{2}x\right)\mathbb{1}_{[0,\,2]}(x) \\
&= \sum_{n\in\mathbb{Z}} h_n \sin^2\left(\frac{\pi}{2}(2x-n)\right)\mathbb{1}_{[0,\,2]}(2x-n) \\
&= \sum_{n=0}^{2} h_n \sin^2\left(\frac{\pi}{2}(2x-n)\right)\mathbb{1}_{[0,\,2]}(2x-n) && \text{by Theorem 3.20}
\end{aligned}$$

(d) The values of $(\!(h_0, h_1, h_2)\!)$ can be conveniently calculated at the knot locations $x = \frac{1}{2}$, $x = 1$, and $x = \frac{3}{2}$ (see the diagram in item 3c page 70):

$$\begin{aligned}
\tfrac{\sqrt{2}}{2}\cdot\frac{1}{2} &= \tfrac{\sqrt{2}}{2}\left(\frac{1}{\sqrt{2}}\right)^2 \\
&= \tfrac{\sqrt{2}}{2}\sin^2\left(\frac{\pi}{4}\right) \\
&\triangleq \tfrac{\sqrt{2}}{2}\phi\left(\frac{1}{2}\right) \\
&= \tfrac{\sqrt{2}}{2}\sqrt{2}\sum_{n\in\mathbb{Z}} h_n \sin^2\left(\frac{\pi}{2}(1-n)\right)\mathbb{1}_{[0,\,2]}(1-n) \\
&= h_0 \sin^2\left(\frac{\pi}{2}(1-0)\right)\mathbb{1}_{[0,\,2]}(1-0) + h_1 \sin^2\left(\frac{\pi}{2}(1-1)\right)\mathbb{1}_{[0,\,2]}(1-1) \\
&\quad + h_2 \sin^2\left(\frac{\pi}{2}(1-2)\right)\mathbb{1}_{[0,\,2]}(1-2)
\end{aligned}$$

$$
\begin{aligned}
&= h_0 \cdot 1 \cdot 1 + h_1 \cdot 0 \cdot 1 + h_2(-1) \cdot 0 \\
&= h_0
\end{aligned}
$$

$$
\begin{aligned}
\tfrac{\sqrt{2}}{2} \cdot 1 &= \tfrac{\sqrt{2}}{2}(1)^2 \\
&= \tfrac{\sqrt{2}}{2} \sin^2\left(\frac{\pi}{2}\right) \\
&\triangleq \tfrac{\sqrt{2}}{2}\phi(1) \\
&= \tfrac{\sqrt{2}}{2}\sqrt{2}\sum_{n\in\mathbb{Z}} h_n \sin^2\left(\frac{\pi}{2}(2-n)\right)\mathbb{1}_{[0,2]}(2-n) \\
&= h_0 \sin^2\left(\frac{\pi}{2}(2-0)\right)\mathbb{1}_{[0,2]}(2-0) + h_1 \sin^2\left(\frac{\pi}{2}(2-1)\right)\mathbb{1}_{[0,2]}(2-1) \\
&\quad + h_2 \sin^2\left(\frac{\pi}{2}(2-2)\right)\mathbb{1}_{[0,2]}(2-2) \\
&= h_0 \cdot 0 \cdot 1 + h_1 \cdot 1 \cdot 1 + h_2 \cdot 0 \cdot 1 \\
&= h_1
\end{aligned}
$$

$$
\begin{aligned}
\tfrac{\sqrt{2}}{2} \cdot \frac{1}{2} &= \tfrac{\sqrt{2}}{2}\left(\frac{1}{-\sqrt{2}}\right)^2 \\
&= \tfrac{\sqrt{2}}{2} \sin^2\left(\frac{3\pi}{4}\right) \\
&\triangleq \tfrac{\sqrt{2}}{2}\phi\left(\frac{3}{2}\right) \\
&= \tfrac{\sqrt{2}}{2}\sqrt{2}\sum_{n\in\mathbb{Z}} h_n \sin^2\left(\frac{\pi}{2}(3-n)\right)\mathbb{1}_{[0,2]}(3-n) \\
&= h_0 \sin^2\left(\frac{\pi}{2}(3-0)\right)\mathbb{1}_{[0,2]}(3-0) + h_1 \sin^2\left(\frac{\pi}{2}(3-1)\right)\mathbb{1}_{[0,2]}(3-1) \\
&\quad + h_2 \sin^2\left(\frac{\pi}{2}(3-2)\right)\mathbb{1}_{[0,2]}(3-2) \\
&= h_0 \cdot (-1) \cdot 0 + h_1 \cdot 0 \cdot 1 + h_2 1 \cdot 1 \\
&= h_2
\end{aligned}
$$

(e) These values for $\langle\!\langle h_0, h_1, h_2 \rangle\!\rangle$ are valid for the knot locations $x = \frac{1}{2}$, $x = 1$, and $x = \frac{3}{2}$, **but** they don't satisfy the *dilation equation* (Theorem 3.4 page 23). In particular,

$$\phi(x) \neq \sum_{n\in\mathbb{Z}} h_n \mathbf{D}\mathbf{T}^n \phi(x)$$

(see the diagram in item 3c page 70)

✏

**Counterexample 6.16** (raised sine) [101] Let a function f be defined in terms of a shifted

[101] 📚 [94], pages 560–561

cosine function (Definition 1.5 page 4) as follows:

$$\phi(x) \triangleq \begin{cases} \frac{1}{2}\left\{1+\cos\left[\pi(|x-1|)\right]\right\} & \text{for } 0 \le x < 2 \\ 0 & \text{otherwise} \end{cases}$$

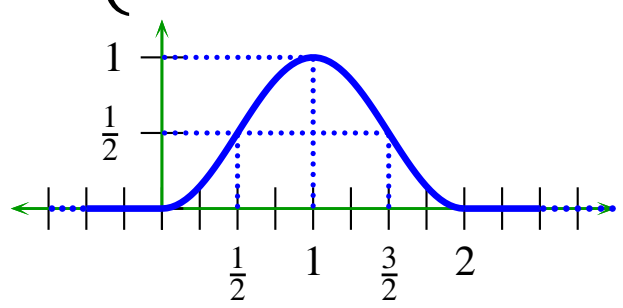

Then $\phi$ forms a *partition of unity*:

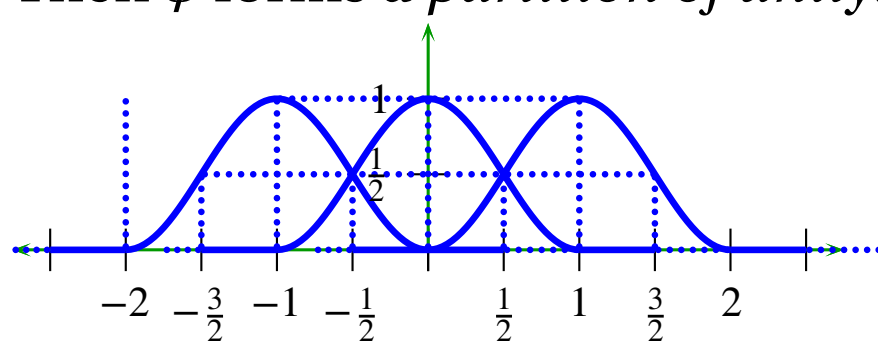

**but** $\left\{\mathbf{T}^n\phi|_{n\in\mathbb{Z}}\right\}$ does **not** generate an *MRA*.

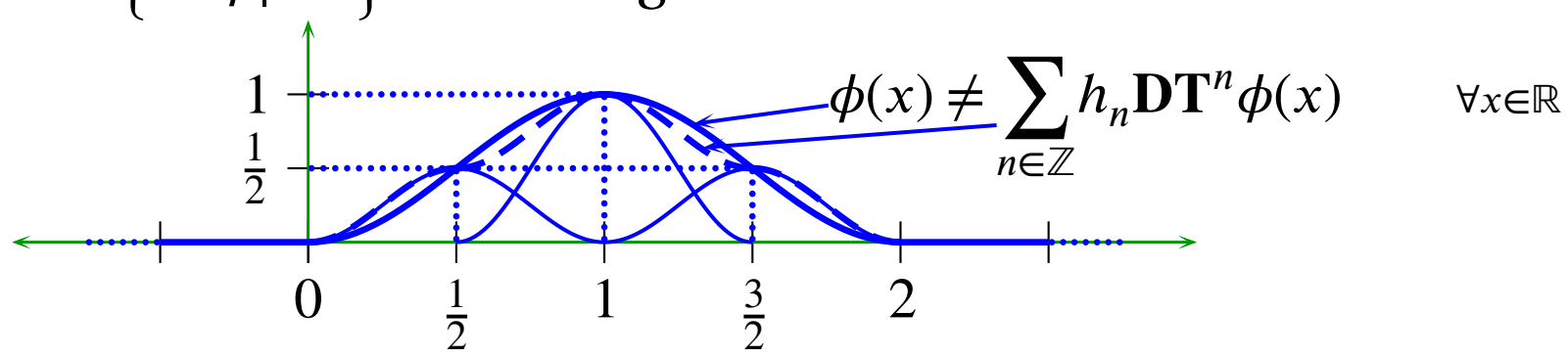

✎PROOF: Let $\mathbb{1}_A(x)$ be the *set indicator function* (Definition 1.3 page 3) on a set $A$.

(1) Proof that $\phi(x)$ forms a *partition of unity*:

$$\begin{aligned}
\sum_{n\in\mathbb{Z}} \mathbf{T}^n\phi(x) &= \sum_{n\in\mathbb{Z}} \mathbf{T}^n\phi(x+1) && \text{by Proposition 2.3 page 15} \\
&= \sum_{n\in\mathbb{Z}} \phi(x+1-n) && \text{by Definition 2.1 page 14} \\
&= \sum_{n\in\mathbb{Z}} \frac{1}{2}\{1+\cos\left[\pi(|x-1+1-n|)\right]\}\mathbb{1}_{[0,\,2)}(x+1-n) && \text{by definition of } \phi(x) \\
&= \sum_{n\in\mathbb{Z}} \frac{1}{2}\{1+\cos\left[\pi(|x-n|)\right]\}\mathbb{1}_{[-1,\,1)}(x-n) && \text{by Definition 1.3 page 3} \\
&= \sum_{n\in\mathbb{Z}} \underbrace{\frac{1}{2}\left\{1+\cos\left[\frac{\pi}{\beta}\left(|x-n|-\frac{1-\beta}{2}\right)\right]\right\}\mathbb{1}_{[-1,\,1)}(x-n)\bigg|_{\beta=1}}_{\textit{raised cosine} \text{ (Example 6.7 page 60) with } \beta=1} \\
&= 1 && \text{by Example 6.7 page 60}
\end{aligned}$$

(2) Proof that $\phi(x) \notin \operatorname{span}\left\{\mathbf{D}\mathbf{T}^n\phi(x)|_{n\in\mathbb{Z}}\right\}$ (and so does not generate an *MRA*):

(a) Note that the *support* (Definition 3.19 page 28) of $\phi$ is $\operatorname{supp}\phi = [0,\,2]$.

(b) Therefore, the *support* of $(\!(h_n)\!)$ is $\operatorname{supp}(\!(h_n)\!) = \{0,1,2\}$ (Theorem 3.20 page 28).

(c) So if $\phi(x)$ *is* an *MRA*, we only need to compute $\{h_0, h_1, h_2\}$ (the rest would be 0). Here would be the values of $\{h_1, h_2, h_3\}$:

$$\begin{aligned}
\phi(x) &= \sum_{n\in\mathbb{Z}} h_n \mathbf{D}\mathbf{T}^n \phi(x) \\
&= \sum_{n\in\mathbb{Z}} h_n \mathbf{D}\mathbf{T}^n \ \frac{1}{2}\left\{ 1+\cos\left[\pi(|x-1|)\right] \right\} \mathbb{1}_{[0,\,2]}(x) && \text{by definition of } \phi(x) \\
&= \sum_{n\in\mathbb{Z}} h_n \ \tfrac{\sqrt{2}}{2}\left\{ 1+\cos\left[\pi(|2x-1-n|)\right] \right\} \mathbb{1}_{[0,\,2]}(2x-n) && \text{by Definition 2.1 page 14} \\
&= \sum_{n=0}^{2} h_n \ \tfrac{\sqrt{2}}{2}\left\{ 1+\cos\left[\pi(|2x-1-n|)\right] \right\} \mathbb{1}_{[0,\,2]}(2x-n) && \text{by Theorem 3.20}
\end{aligned}$$

(d) The values of $(\!(h_0, h_1, h_2)\!)$ can be conveniently calculated at the knot locations $x = \frac{1}{2}$, $x = 1$, and $x = \frac{3}{2}$ (see the diagram in item 3c page 70):

$$\begin{aligned}
\frac{1}{2} &= \left.\sum_{n=0}^{2} h_n \ \tfrac{\sqrt{2}}{2}\left\{ 1+\cos\left[\pi(|2x-1-n|)\right] \right\} \mathbb{1}_{[0,\,2]}(2x-n)\right|_{x=\frac{1}{2}} \\
&= h_0 \ \tfrac{\sqrt{2}}{2}\left\{ 1+\cos\left[1-1-0\right] \right\} \\
&= h_0\sqrt{2} \\
&\implies h_0 = \frac{\sqrt{2}}{4}
\end{aligned}$$

$$\begin{aligned}
1 &= \left.\sum_{n=0}^{2} h_n \ \tfrac{\sqrt{2}}{2}\left\{ 1+\cos\left[\pi(|2x-1-n|)\right] \right\} \mathbb{1}_{[0,\,2]}(2x-n)\right|_{x=1} \\
&= h_1 \ \tfrac{\sqrt{2}}{2}\left\{ 1+\cos\left[2-1-1\right] \right\} \\
&= h_1\sqrt{2} \\
&\implies h_1 = \tfrac{\sqrt{2}}{2}
\end{aligned}$$

$$\begin{aligned}
\frac{1}{2} &= \left.\sum_{n=0}^{2} h_n \ \tfrac{\sqrt{2}}{2}\left\{ 1+\cos\left[\pi(|2x-1-n|)\right] \right\} \mathbb{1}_{[0,\,2]}(2x-n)\right|_{x=\frac{3}{2}} \\
&= h_2 \ \tfrac{\sqrt{2}}{2}\left\{ 1+\cos\left[1-1-0\right] \right\} \\
&= h_2\sqrt{2} \\
&\implies h_2 = \frac{\sqrt{2}}{4}
\end{aligned}$$

(e) These values for $\langle\!\langle h_0, h_1, h_2 \rangle\!\rangle$ are valid for the knot locations $x = \frac{1}{2}$, $x = 1$, and $x = \frac{3}{2}$, **but** they don't satisfy the *dilation equation* (Theorem 3.4 page 23). In particular (see diagram),

$$\phi(x) \neq \sum_{n\in\mathbb{Z}} h_n \mathbf{D}\mathbf{T}^n \phi(x) .$$

✏

*Telecommunications Engineering Department, National Chiao-Tung University, Hsinchu, Taiwan*

dgreenhoe@gmail.com

[102] This document was typeset using XƎLATEX, which is part of the TEX family of typesetting engines, which is arguably the greatest development since the Gutenberg Press. Graphics were rendered using the pstricks and related packages and LATEX graphics support. The main roman, *italic*, and **bold** font typefaces used are all from the Heuristica family of typefaces (based on the Utopia typeface, released by Adobe Systems Incorporated). The math font is XITS from the XITS font project. The font used for the title in the footers is Adventor (similar to Avant-Garde) from the TEX-Gyre Project. The font used for the version number in the footers is LIQUID CRYSTAL (Liquid Crystal) from FontLab Studio. This handwriting font is Lavi from the Free Software Foundation. Chinese glyphs appearing in the text are from the font 王漢宗中明體繁 (pinyin: Wáng Hàn Zōng Zhōng Mińg Tǐ Fán; translation: Hàn Zōng Wáng's Medium-weight Mińg-style Traditional Characters; literal: 王漢宗~font designer's name; 中~medium; 明~Mińg (a dynasty); 體~style; 繁~traditional).